\documentclass[12pt]{elsarticleNoFoot}
\usepackage[a4paper,text={16cm,24cm},centering]{geometry}
\usepackage[latin1]{inputenc}
\usepackage[pdftex]{color}
\usepackage{graphicx}
\usepackage{amssymb}
\usepackage[hidelinks]{hyperref}
\usepackage{karn}

\renewcommand{\arraystretch}{1.0}

\makeatletter

\def\newlog#1#2{%
	\def#1{\mathop{\operator@font #2}\nolimits}%
}
\def\Ite#1#2#3{\setbox\@tempboxa\hbox{$#1$}%
\ifdim\wd\@tempboxa=0cm#3\else#2\fi\relax}

\makeatother

\definecolor{coLabel}	{rgb}{0.50,0.00,0.00}	% label (debug only)
\definecolor{coRemark}	{rgb}{0.99,0.00,0.00}	% label (debug only)

\definecolor{coRef}	{rgb}{0.00,0.00,0.99}	% reference
\definecolor{coThmName}	{rgb}{0.00,0.00,0.99}
\definecolor{coQNeg}	{rgb}{0.40,0.00,0.00}	% qmc negation
\definecolor{coQDef}	{rgb}{0.00,0.00,0.60}	% qmc definition
\definecolor{coQDer}	{rgb}{0.60,0.60,0.99}	% qmc derived
\definecolor{coQRef}	{rgb}{0.00,0.00,0.50}	% qmc reference
\definecolor{coDpd}	{rgb}{0.00,0.50,0.00}	% proof depth at def.
\definecolor{coDpr}	{rgb}{0.50,0.99,0.50}	% proof depth at ref.

\definecolor{coQii}	{rgb}{0.00,0.00,0.50}	% level 2 lemma
\definecolor{coQiii}	{rgb}{0.30,0.30,0.75}	% level 3 lemma
\definecolor{coQiv}	{rgb}{0.60,0.60,0.80}	% level 4 lemma
\definecolor{coQv}	{rgb}{0.70,0.70,0.80}	% level 5 lemma
\definecolor{coQvi}	{rgb}{0.80,0.80,0.80}	% level 6 lemma
\newcommand{\LABEL}[1]{\label{#1}}

\newcommand{\REF}[1]{\textcolor{coRef}{\ref{#1}}}
\newcommand{\REFF}[2]{\textcolor{coRef}{\ref{#1}.\ref{#1 #2}}}

\newcommand{\La}{\Leftarrow}
\newcommand{\Lra}{\Leftrightarrow}
\newcommand{\ra}{\rightarrow}
\newcommand{\lra}{\leftrightarrow}
\newcommand{\Ra}{\Rightarrow}
\renewcommand{\leq}{\leqslant}
\renewcommand{\geq}{\geqslant}

\newlog{\dom}{dom}
\newlog{\ran}{ran}

\newcommand{\tpl}[1]{\langle #1 \rangle}	% tuple
\newcommand{\set}[1]{\{ #1 \}}			% set
\newcommand{\disjUnion}{\mathop{\stackrel{.}{\cup}}}
\newcommand{\false}{{\it false}}
\newcommand{\true}{{\it true}}
\newcommand{\N}{I\!\!N}
\newcommand{\Z}{Z\!\!\!Z}

\newcommand{\dpd}[1]{$^{\textcolor{coDpd}{[#1]}}$}	% depth definition
\newcommand{\dpr}[1]{}	% depth reference

\newcommand{\QDef}[1]{\hypertarget{QL#1}{\textcolor{coQDef}{\sf #1}}}

\newcommand{\QDer}[1]{\hypertarget{QL#1}{$^*$\textcolor{coQDer}{\sf #1}}}% follows from qmc output
\newcommand{\QDeR}[1]{\hypertarget{QL#1}{\textcolor{coQDef}{\sf #1}}}% follows from manually proved lemma only
\newcommand{\QDee}[1]{\hypertarget{QL#1}{$^e$\textcolor{coQDer}{\sf #1}}}% example only

\newcommand{\QRef}[1]{\hyperlink{QL#1}{\textcolor{coQRef}{\sf #1}}}%

\newcommand{\QNeg}[1]{\textcolor{coQNeg}{$\lnot$#1}}

\newcommand{\restrict}[2]{#1 \!\!\mid_{#2}}

\newcommand{\eqc}[2]{[#1]_{#2}}

\renewcommand{\:}[4]{%
	{%
	\renewcommand{\:}[4]{%
		{%
		\renewcommand{\:}[4]{error\error}%
		\renewcommand{\j}{{##2}}%
		{##4}%
		##1...##1%
		\renewcommand{\j}{{##3}}%
		{##4}%
		}%
	}%
	\renewcommand{\i}{{#2}}%
	{#4}%
	#1\ldots#1%
	\renewcommand{\i}{{#3}}%
	{#4}%
	}%
}

\newproof{pf}{Proof}
\newcommand{\THM}[3]{%
	\vspace{0ex plus 1ex}%
	\begin{#1}%
	\Ite{#2}{ \textcolor{coThmName}{(#2)} }{ }%
	{\rm #3}%
	\end{#1}%
	}
\newcommand{\PRF}[2]{%
	\vspace{-2ex plus 1ex}%
	\begin{pf}%
	{#2}%
	\end{pf}%
	\vspace{0ex plus 1ex}%
	}
\newcommand{\LEMMA}[2]		{\THM{lemma}		{#1}{#2}}
\newcommand{\DEFINITION}[2]	{\THM{definition}	{#1}{#2}}
\newcommand{\EXAMPLE}[2]	{\THM{example}		{#1}{#2}}

\newcommand{\PROOF}[1]		{\PRF{\em Proof. }	{#1}}

\newcommand{\CITE}[1]{\citet{#1}}
\newcommand{\CITEPAGE}[2]{\citet[][#1]{#2}}

\renewcommand{\.}[1]{\!#1\!}

\begin{document}

\newsavebox{\atpic}
\savebox{\atpic}{\includegraphics[scale=0.08]{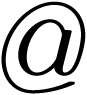}}

\begin{frontmatter}

\title{Simple Laws about\\Nonprominent Properties of Binary Relations}
\author{Jochen Burghardt}
\address{jochen.burghardt\usebox{\atpic}alumni.tu-berlin.de}
\address{\rm Nov 2018}

\begin{abstract}
We checked each binary relation on a 5-element set for a given set of
properties, including usual ones like asymmetry and less known ones
like Euclideanness.
%~
Using a poor man's Quine-McCluskey algorithm, we computed prime
implicants of non-occurring property combinations, like
\textit{``not irreflexive, but asymmetric''}.
%~
We considered the laws obtained this way, and manually proved
them true for binary relations on arbitrary sets, thus contributing to
the encyclopedic knowledge about less known properties.
\end{abstract}

\begin{keyword}
Binary relation;
Quine-McCluskey algorithm;
Hypotheses generation
%reflexivity;
%irreflexivity;
%co-reflexivity;
%symmetry;
%asymmetry;
%anti-symmetry;
%transitivity;
%anti-transitivity;
%quasi-transitivity;
%Euclideanness;
%semi-connex relation;
%connex relation;
%semi-order
\end{keyword}

\end{frontmatter}

\vfill

{

\setlength{\parskip}{0cm}

\clearpage
\tableofcontents

%\vfill

%\clearpage
\listoffigures

%\vfill
}

\clearpage

\section{Introduction}
\LABEL{Introduction}

In order to flesh out encyclopedic articles\footnote{
	at \url{https://en.wikipedia.org}
}
about less common
properties (like e.g.\ anti-transitivity) of binary relations,
we implemented a simple \texttt{C} program to iterate over all
relations on a small
finite set and to check each of them for given
properties.
%~
We implemented checks for the properties given in
Def.~\REF{def} below.
%~
Figure~\REF{Source code for transitivity check}
shows the \texttt{C} source code to check a relation
\texttt{R} for transitivity, where \texttt{card} is the universe size
and \texttt{elemT} is the type encoding a universe element.

This way, we could, in a first stage,
(attempt to) falsify intuitively found hypotheses
about laws involving such properties, and search for illustrative
counter-examples to known, or intuitively guessed, non-laws.
%~
For example,
Fig.~\REF{Source code to search for right Euclidean non-transitive relations}
shows the source code to search for right Euclidean non-transitive
relations over a $2$-element universe, where \texttt{printRel} prints
its argument relation in a human-readable form.
%~
For a universe of $n$ elements, $n^2$ \texttt{for} loops are nested.
%~
In Sect.~\REF{Improved relation enumeration} we describe an improved
way to iterate over all relations.

Relations on a set of up to $6$ elements could be dealt with in
reasonable time on a 2.3 GHz CPU.
%~
Figure~\REF{Timing vs. universe cardinality}
gives an overview, where all times are wall clock
times in seconds,
and ``tr$\Ra$qt'' indicates the task of validating
that each transitive binary relation is also quasi-transitive.
%~
Note the considerable amount of compile time,\footnote{
	We used \texttt{gcc} version 7.3.0 with the highest
	optimization level.
}
presumably
caused by excessive use
of inlining, deeply nested loops, and abuse of array elements as loop
variables.

In a second stage, we aimed at supporting the generation of
law hypotheses, rather than their validation.

We used a 5-element universe set, and
checked each binary relation for each of the properties.\footnote{
	For this run, we hadn't provided checks for left and right
	quasi-reflexivity
	(Def.~\REFF{def}{LfQuasiRefl}+\REF{def RgQuasiRefl}),
	but only for the conjunction of both,
	viz.\ quasi-reflexivity (Def.~\REFF{def}{QuasiRefl}).
	%~
	As additional properties, we provided a check for the empty
	relation ($\forall x, y \in X. \;\; \lnot xRy$)
	and for the universal relation
	($\forall x, y \in X. \;\; xRy$).
}
%~
The latter were encoded by bits of a 64-bit word.
%~
After that, we applied
a poor-man's Quine-McCluskey algorithm\footnote{
	See \CITE{Quine.1952} and \CITE{McCluskey.1956}
	for the original algorithm.
}
(denoted ``QMc'' in Fig.~\ref{Timing vs. universe cardinality})
to obtain a short description
of property combinations that didn't occur at all.
%~
For example, an output line ``\verb|~Irrefl ASym|'' indicated that the
program didn't find any relation that was asymmetric but not
irreflexive, i.e.\ that each asymmetric relation on a 5-element set is
irreflexive.
%~
Section~\REF{Reported law suggestions} shows the complete output
on a 5-element universe.

\begin{figure}
\begin{center}
\small
\begin{verbatim}
bool isTrans(const bool R[card][card]) {
  elemT x,y,z;
  for (x=0; x<card; ++x)
    for (y=0; y<card; ++y)
      if (R[x][y])
        for (z=0; z<card; ++z)
          if (R[y][z] && ! R[x][z])
            return false;
  return true;
}
\end{verbatim}
\caption{Source code for transitivity check}
\LABEL{Source code for transitivity check}
\end{center}
\end{figure}

\begin{figure}
\begin{center}
\small
\begin{verbatim}
void check03(void) {
  bool R[card][card];
  for (R[0][0]=false; R[0][0]<=true; ++R[0][0])
    for (R[0][1]=false; R[0][1]<=true; ++R[0][1])
      for (R[1][0]=false; R[1][0]<=true; ++R[1][0])
        for (R[1][1]=false; R[1][1]<=true; ++R[1][1])
          if (isRgEucl(R) && ! isTrans(R))
            printRel(R);
}
\end{verbatim}
\caption{Source code to search for right Euclidean non-transitive relations}
\LABEL{Source code to search for right Euclidean non-transitive relations}
\end{center}
\end{figure}

\begin{figure}
\begin{center}
\begin{tabular}{|l|r|r|r|r|r|}
\hline
Universe card
	& 2
	& 3
	& 4
	& 5
	& 6
%	& 7
	\\
\hline
Relation count
	& 16				% 2
	& 512				% 3
	& 6.55e04			% 4
	& 3.35e07			% 5
	& 6.87e10			% 6
%	& 5.62e14			% 7
	\\
\cline{2-6}
	& 10				% 2
	& 140				% 3
	& 6.17e03			% 4
	& 9.07e05			% 5
	& 4.60e08			% 6
%	& 8.27e11			% 7
	\\
\hline
Compile time
	& 7.123				% 2
	& 14.254			% 3
	& 20.868			% 4
	& 27.923			% 5
	& 41.965			% 6
%	& 30.939			% 7
	\\
Run time tr$\Ra$qt
	& 0.007				% 2
	& 0.007				% 3
	& 0.009				% 4
	& 0.132				% 5
	& 50.386			% 6
%	&				% 7	>1d estimated
	\\
%Run time QMc
%	&				% 2
%	&				% 3
%	&				% 4
%	&				% 5
%	& --				% 6
%	& --				% 7
%	\\
\hline
\end{tabular}
\caption{Timing vs. universe cardinality}
\LABEL{Timing vs. universe cardinality}
\end{center}
\end{figure}

We took each printed law as a suggestion to be proven for all binary
relations (on arbitrary sets).
%~
Many of the considered laws were trivial, in particular those involving
co-reflexivity, as this property applies only to a relatively small
number of relations (32 on a 5-element set).

A couple of laws appeared to be interesting, and we could prove them
fairly easily by hand for the general case\footnote{
	We needed to require a minimum cardinality of the universe
	set in some lemmas, e.g.\ Lem.~\REF{conn 2} and~\REF{eucl 7}.
}.
%~
For those laws involving less usual properties (like anti-transitivity,
quasi-transitivity, Euclideanness) there is good chance that they
haven't been stated in the literature before.
%~
However, while they may contribute to the completeness of an encyclopedia,
it is not clear whether they may serve any other purpose.

Disregarding the particular area of binary relations, the method of
computing law suggestions by the Quine-McCluskey algorithm might be
used as a source of fresh exercises whose solutions are unlikely to be
found on web pages.

Some of the laws, e.g.\ Lem.~\REF{eucl 5}, appeared
surprising, but turned out during the proof to be vacuously true.
%~
The proof attempt to some laws gave rise to the assertion of other
lemmas that weren't directly obtained from the computed
output:
%~
Lemma~\REF{corefl 1} was needed for the proof of
Lem.~\REF{quasiTrans 5}, and
Lem.~\REF{conn 3} was needed for Lem.~\REF{eucl 7}.

Our Quine-McCluskey approach restricts law suggestions to formulas of
the form $\forall R. \; \:{\lor}1n{{\it prop}_\i(R)}$,
where the quantification is over all binary relations,
and ${\it prop}_i$ is one of the considered properties or a negation
thereof.

\begin{figure}
\begin{center}
$$\begin{array}{l|rrr|}
\multicolumn{1}{l}{_x \setminus ^y} & 0 & 1 & \multicolumn{1}{r}{2}	\\
\cline{2-4}
0 & 1 & 2 & 4	\\
1 & 8 & 16 & 32	\\
2 & 64 & 128 & 256	\\
\cline{2-4}
\end{array}$$
\caption{Encoding scheme for $3 \times 3$ relations
	for a \CITE{Burghardt.2002b} approach}
\LABEL{Encoding scheme for 3 times 3 relations}
\end{center}
\end{figure}

\begin{figure}
\begin{center}
$$\begin{array}{cr c@{\;}c@{\;}c@{\;}c@{\;}c@{\;}c@{\;}c@{\;}c@{\;}c}
R_0
	& ::=  & {\it symCls}(R_0)
	& \mid & R_0 \cup R_0	\\
\vdots	\\
R_{10}
	& ::=  & {\it symCls}(R_2)
	& \mid & {\it symCls}(R_8)
	& \mid & {\it symCls}(R_{10})	\\
	& \mid & R_0 \cup R_{10}
	& \mid & R_2 \cup R_8
	& \mid & R_2 \cup R_{10}
	& \mid & R_8 \cup R_2
	& \mid & R_8 \cup R_{10}	\\
	& \mid & R_{10} \cup R_0
	& \mid & R_{10} \cup R_2
	& \mid & R_{10} \cup R_8
	& \mid & R_{10} \cup R_{10}	\\
\vdots	\\
R_{511}
	& ::=  & {\it symCls}(R_{311})
	& \mid & \ldots
	& \mid & {\it symCls}(R_{511})	\\
	& \mid & R_0 \cup R_{511}
	& \mid & \ldots
	& \mid & R_{511} \cup R_{511}	\\
[2ex]
\true
	& ::=  & {\it isRefl}(R_{273})
	& \mid & \ldots
	& \mid & {\it isRefl}(R_{511})	\\
	& \mid & {\it isSym}(R_0)
	& \mid & \ldots
	& \mid & {\it isSym}(R_{511})	\\
	& \mid & \false \lor \true
	& \mid & \true \lor \false
	& \mid & \true \lor \true
	& \mid & \true \land \true
	& \mid & \lnot \false	\\
\false
	& ::=  & {\it isRefl}(R_0)
	& \mid & \ldots
	& \mid & {\it isSym}(R_2)
	& \mid & \ldots
	& \mid & \lnot \true	\\
\end{array}$$
\caption{Tree grammar sketch for \CITE{Burghardt.2002b} approach}
\LABEL{Tree grammar sketch for Burghardt.2002b approach}
\end{center}
\end{figure}

For an approach to compute more general forms of law suggestions,
see \CITE{Burghardt.2002b};
however, due to its run-time complexity
this approach is feasible only for even smaller universe
sets.
%~
In order to handle all relations on a 3-element set,
a regular tree grammar of 512 nonterminals, one for each
relation, plus 2 nonterminals, one for each truth value,
would be needed.
%~
Using the encoding scheme from
Fig.~\REF{Encoding scheme for 3 times 3 relations},
the original grammar would consist of rules as sketched\footnote{
	For sake of simplicity, only one unary and one binary
	operation on relations is considered, viz.\ symmetric closure
	${\it symCls}$ and union $\cup$.
	%~
	Only two properties of relations are considered, viz.\
	reflexivity ${\it isRefl}$ and symmetry ${\it isSym}$.
	%~
	It should be obvious how to incorporate more operators and
	predicates on relations.
	%~
	By additionally providing a sort for sets, operations like
	$\dom$, $\ran$, restriction, etc.\ could be considered also.
}
in Fig.~\REF{Tree grammar sketch for Burghardt.2002b approach}.
%~
However, this grammar grows very large, and its $n$-fold
product would be needed if all laws in $n$ variables were to be computed.

The rest of this paper is organized as follows.
%~
In Sect.~\REF{Definitions}, we formally define each considered
property, and introduce some other notions.
%~
In Sect.~\REF{Reported law suggestions}, we show the annotated output
for a run of our algorithm on a 5-element set, also indicating which
law suggestions gave rise to which lemmas.
%~
The latter are stated and proven in
Sect.~\REF{Formal proofs of property laws}, which is the main part of
this paper.
%~
In addition, we state the proofs of some laws that weren't of the form
admitted by our approach; some of them were, however, obtained using
the assistance of the counter-example search in our \texttt{C}
program.
%~
In Sect.~\REF{Examples},
we discuss those computed law suggestions that lead
to single examples, rather than to general laws.
%~
In Sect.~\REF{Implementation issues},
we comment on some program details.

This paper is a follow-up version of
\url{https://arxiv.org/abs/1806.05036v1}.
%~
Compared to the previous version,
we considered $9$ more properties (see Def.~\REF{def}), including
being the empty and being the universal relation, to avoid
circumscriptions like ``Irrefl$\lor$CoRefl$\lor\lnot$ASym''
in favor of ``Empty$\lor\lnot$ASym'';
in the new setting,
we found a total of $274$ law suggestions, and proved or disproved
all of them.
%~
I am thankful to all people who have helped with their comments and
corrections.

\clearpage

\section{Definitions}
\LABEL{Definitions}

\DEFINITION{Binary relation properties}{%
\LABEL{def}%
Let $X$ be a set.
%~
A (homogeneous) binary relation $R$ on $X$ is a subset of $X \times X$.
%~
The relation $R$ is called
\begin{enumerate}
\item\LABEL{def Refl}
	reflexive
	(``Refl'', ``rf'')
	~ if ~ $\forall x \in X. \;\; xRx$;
\item\LABEL{def Irrefl}
	irreflexive
	(``Irrefl'', ``ir'')
	~ if ~ $\forall x \in X. \;\; \lnot xRx$;
\item\LABEL{def CoRefl}
	co-reflexive
	(``CoRefl'', ``cr'')
	~ if ~ $\forall x,y \in X. \;\; xRy \ra x=y$;
\item\LABEL{def LfQuasiRefl}
	left quasi-reflexive
	(``lq'')
	~ if ~ $\forall x,y \in X. \;\; xRy \ra xRx$;
\item\LABEL{def RgQuasiRefl}
	right quasi-reflexive
	(``rq'')
	~ if ~ $\forall x,y \in X. \;\; xRy \ra yRy$;
\item\LABEL{def QuasiRefl}
	quasi-reflexive
	(``QuasiRefl'')
	~ if ~ it is both left and right quasi-reflexive;
\item\LABEL{def Sym}
	symmetric
	(``Sym'', ``sy'')
	~ if ~ $\forall x,y \in X. \;\; xRy \ra yRx$;
\item\LABEL{def ASym}
	asymmetric
	(``ASym'', ``as'')
	~ if ~ $\forall x,y \in X. \;\; xRy \ra \lnot yRx$;
\item\LABEL{def AntiSym}
	anti-symmetric
	(``AntiSym'', ``an'')
	~ if ~ $\forall x,y \in X. \;\; xRy \land x \neq y \ra \lnot yRx$;
\item\LABEL{def SemiConnex}
	semi-connex
	(``SemiConnex'', ``sc'')
	 if ~ $\forall x,y \in X. \;\; xRy \lor yRx \lor x=y$;
\item\LABEL{def Connex}
	connex
	(``Connex'', ``co'')
	 if ~ $\forall x,y \in X. \;\; xRy \lor yRx$;
\item\LABEL{def Trans}
	transitive
	(``Trans'', ``tr'')
	~ if ~ $\forall x,y,z \in X. \;\; xRy \land yRz \ra xRz$;
\item\LABEL{def AntiTrans}
	anti-transitive
	(``AntiTrans'', ``at'')
	~ if ~ $\forall x,y,z \in X. \;\; xRy \land yRz \ra \lnot xRz$;
\item\LABEL{def QuasiTrans}
	quasi-transitive
	(``QuasiTrans'', ``qt'')
	~ if ~ $\forall x,y,z \in X. \;\;
	xRy \land \lnot yRx \land yRz \land \lnot zRy
	\ra xRz \land \lnot zRx$;
\item\LABEL{def RgEucl}
	right Euclidean
	(``RgEucl'', ``re'')
	~ if ~ $\forall x,y,z \in X. \;\; xRy \land xRz \ra yRz$;
\item\LABEL{def LfEucl}
	left Euclidean
	(``LfEucl'', ``le'')
	~ if ~ $\forall x,y,z \in X. \;\; yRx \land zRx \ra yRz$;
\item\LABEL{def SemiOrd1}
	semi-order property 1
	(``SemiOrd1'', ``s1'')
	~ if ~ $\forall w,x,y,z \in X. \;\;
	wRx \land \lnot xRy \land \lnot yRx \land yRz \ra wRz$;
\item\LABEL{def SemiOrd2}
	semi-order property 2
	(``SemiOrd2'', ``s2'')
	~ if ~ $\forall w,x,y,z \in X. \;\;
	xRy \land yRz
	\ra wRx \lor xRw \lor wRy \lor yRw \lor wRz \lor zRw$.
\item\LABEL{def RgSerial}
	right serial
	(``RgSerial'', ``rs'')
	% formerly: serial
	~ if ~ $\forall x \in X \; \exists y \in X. \;\; xRy$
\item\LABEL{def LfSerial}
	left serial
	(``LfSerial'', ``ls'')
	% formerly: inverse serial
	~ if ~ $\forall y \in X \; \exists x \in X. \;\; xRy$
\item\LABEL{def Dense}
	dense
	(``Dense'', ``de'')
	~ if ~ $\forall x, z \in X \; \exists y \in X. \;\;
	xRz \ra xRy \land yRz$.
	% inspired by "dense order", needed for "idempotent relation"
\item\LABEL{def IncTrans}
	incomparability-transitive
	(``IncTrans'', ``it'')
	~ if ~ $\forall x, y, z \in X. \;\;
	\lnot xRy \land \lnot yRx \land \lnot yRz \land \lnot zRy
	\ra \lnot xRz \land \lnot zRx$.
\item\LABEL{def LfUnique}
	left unique
	(``LfUnique'', ``lu'')
	~ if ~ $\forall x_1, x_2, y \in X \;\;
	x_1 R y \land x_2 R y \ra x_1 = x_2$.
\item\LABEL{def RgUnique}
	right unique
	(``RgUnique'', ``ru'')
	~ if ~ $\forall x, y_1, y_2 \in X \;\;
	x R y_1 \land x R y_2 \ra y_1 = y_2$.
\end{enumerate}
The capitalized abbreviations in parentheses are used by our
algorithm; the two-letter codes are used in tables and pictures when
space is scarce.

The ``left'' and ``right'' properties are dual to each other.
%~
All other properties are self-dual.
%~
For example, a relation $R$ is left unique iff its converse, $R^{-1}$,
is right unique;
a relation $R$ is dense iff its converse is dense.

We say that $x,y$ are incomparable w.r.t.\ $R$,
if $\lnot xRy \land \lnot yRx$ holds.
\qed
}

\DEFINITION{Kinds of binary relations}{%
\LABEL{Kinds of binary relations}%
A binary relation $R$ on a set $X$ is called
\begin{enumerate}
\item an equivalence
	~ if ~
	it is reflexive, symmetric, and transitive;
\item a partial equivalence
	~ if ~
	it is symmetric and transitive;

\item a tolerance relation
	~ if ~
	it is reflexive and symmetric;
\item idempotent
	~ if ~
	it is dense and transitive;
\item trichotomous
	~ if ~
	it is irreflexive, asymmetric, and semi-connex;

\item a non-strict partial order
	~ if ~
	it is reflexive, anti-symmetric, and transitive;
\item a strict partial order
	~ if ~
	it is irreflexive, asymmetric, and transitive;
\item\LABEL{semi-order}
	a semi-order
	~ if ~
	it is asymmetric and satisfies semi-order properties 1 and 2;
\item a preorder
	~ if ~
	it is reflexive and transitive;
\item a weak ordering
	~ if ~
	it is irreflexive, asymmetric, transitive, and
	incomparability-transitive;

\item a partial function
	~ if ~
	it is right unique;
\item a total function
	~ if ~
	it is right unique and right serial;
\item an injective function
	~ if ~
	it is left unique, right unique, and right serial;
\item a surjective function
	~ if ~
	it is right unique and and left and right serial;
\item a bijective function
	~ if ~
	it is left and right unique and left and right serial.
\qed
\end{enumerate}
}

\DEFINITION{Operations on relations}{%
\begin{enumerate}
\item
	For a relation $R$ on a set $X$ and a subset $Y \subseteq X$,
	we write $\restrict{R}{Y}$ for the restriction of $R$ to $Y$.
	%~
	Formally, $\restrict{R}{Y}$ is the relation on $Y \times Y$
	defined by $x (\restrict{R}{Y}) y :\Lra xRy$ for each $x,y \in Y$.
\item
	For an equivalence relation $R$ on a set $X$, we write $\eqc{x}{R}$
	for the equivalence class of $x \in X$ w.r.t.\ $R$.
	%~
	Formally, $\eqc{x}{R} := \set{ y \in X \mid xRy }$.
\item
	For a relation $R$ on a set $X$ and $x, y \in X$, we write
	$xR$ for the set of elements $x$ is related to,
	and $Ry$ for the set of elements that are related to $y$.
	%~
	Formally,
	$xR := \set{ y \in X \mid xRy }$
	and $Ry := \set{ x \in X \mid xRy }$.
\qed
\end{enumerate}
}

\clearpage
\section{Reported law suggestions}
\LABEL{Reported law suggestions}

\newcommand{\x}{\textcolor{black}{\sf X}}

In this section,
we show the complete output produced by our Quine-McCluskey algorithm run.

In the Fig.~\REF{Reported laws for level 2}
to~\REF{Reported laws for level 8},
we list the computed prime implicants for missing relation property
combinations on a 5-element universe set.
%~
We took each prime implicant as a suggested law about all binary
relations.
%~
These suggestions are grouped by the number of their
literals (``level'').

In the leftmost column, we provide a consecutive law number for
referencing.
%~
In the middle column, the law is given in textual
representation, ``\QNeg{P}'' denoting the negation of P,
and juxtaposition used for conjunction.
%~
The property names correspond to those used
by the \texttt{C} program; they should be understandable without
further explanation, but can also be looked up via
Fig.~\REF{Number of relations on a $5$-element set}, if necessary.
%~
In the rightmost column, we annotated a reference to the lemma
(in Sect.~\REF{Formal proofs of property laws})
where the law has been formally proven or to the example
(in Sect.~\REF{Examples}) where it is discussed.

For example, line \QRef{039}, in level 2
(Fig.~\REF{Reported laws for level 2} left),
reports that no relation was found to be
asymmetric (property~\REFF{def}{ASym})
and non-irreflexive (negation of property~\REFF{def}{Irrefl});
we show the formal proof that every asymmetric relation is irreflexive in
Lem.~\REFF{asym 1b}{1}.\footnote{
	A warning about possible confusion appears advisable here:
	%~
	In the setting of the Quine-McCluskey algorithm, a prime
	implicant is a conjunction of negated and/or unnegated variables.
	%~
	However, its corresponding law suggestion is its complement,
	and hence a disjunction, as should be clear from the example.
	%~
	Where possible, we used the term ``literal'' in favor of
	``conjunct'' or ``disjunct''.
}

Laws that could be derived from others by purely propositional
reasoning and without referring to the
property definitions in Def.~\REF{def} are considered redundant;
they are marked with a star ``$^*$''.\footnote{
	We marked all redundancies we became aware of; we don't claim
	that no undetected ones exist.
}
%~
For example, law \QRef{044} (``no relation is asymmetric and reflexive'')
is marked since it
follows immediately from \QRef{046}
(``no relation is irreflexive and reflexive'')
and \QRef{039}.

No laws were reported for level 1 and level 9 and beyond.
~
A text version of these tables is available in the ancillary file
\texttt{reportedLaws.txt} at \texttt{arxiv.org}.

In Fig.~\REF{Law index + - lf} to~\REF{Law index - - / + + rg},
we summarize the found laws.
%~
We omitted suggestions that couldn't be
manually verified as laws, and suggestions marked as redundant.

Figure~\REF{Law index + - lf} and~\REF{Law index + - rg}
shows the left and right half of an implication table,
respectively.
%~
Every field lists all law numbers that can possibly be used
to derive the column property from the row property.

For example, law \QRef{129} appears in line ``tr''
(transitive) and column ``as'' (asymmetric)
in Fig.~\REF{Law index + - lf}
because that law
(well-known, and proven in Lem.~\REFF{asym 1a}{2})
allows one to infer a relation's asymmetry from its transitivity,
provided that it is also known to be irreflexive.

Fields belonging to the table's diagonal are marked by ``\x''.
%~
Law numbers are colored by number of literals,
deeply-colored and pale-colored numbers indicating
few and many literals, respectively.

Similarly, the table consisting of
Fig.~\REF{Law index - - / + + lf}
and~\REF{Law index - - / + + rg}
shows below and above
its diagonal laws about required disjunctions and
impossible conjunctions, respectively.

For example, law \QRef{223} appears below the diagonal in
line ``co'' (connex) and column ``em'' (empty) of
Fig.~\REF{Law index - - / + + lf}, since the law
(proven in Lem.~\REF{incTrans 6})
requires every relation to be connex
or empty, provided it is quasi-reflexive and incomparability-transitive.

Law \QRef{145} appears above the diagonal in line ``le''
(left Euclidean) and column ``lu'' (left unique), since the law
(proven in Lem.~\REF{eucl 10})
ensures that no relation
can be left Euclidean and
left unique, provided it isn't anti-symmetric.

Figure~\REF{Implications and incompatibilities between properties}
shows all proper implications (black) and incompatibilities
(red) from level 2, except for the empty and the universal relation.
%~
Vertex labels use the abbreviations from Fig.~\REF{Law index + - lf},
edge labels refer to law numbers in Fig.~\REF{Reported laws for level 2}.

\renewcommand{\arraystretch}{0.9}

\begin{figure}
\footnotesize
\begin{tabular}{|r|l|c|}
\hline
\QDee{001} & Empty Univ & \REFF{empty}{10}	\\\hline
\QDee{002} & Empty \QNeg{CoRefl} & \REFF{empty}{1}	\\\hline
\QDee{003} & Univ CoRefl & \REFF{univ}{10}	\\\hline
\QDee{004} & Empty \QNeg{LfEucl} & \REFF{empty}{2}	\\\hline
\QDee{005} & Univ \QNeg{LfEucl} & \REFF{univ}{1}	\\\hline
\QDef{006} & CoRefl \QNeg{LfEucl} & \REFF{corefl 5}{1}	\\\hline
\QDee{007} & Empty \QNeg{RgEucl} & \REFF{empty}{2}	\\\hline
\QDee{008} & Univ \QNeg{RgEucl} & \REFF{univ}{1}	\\\hline
\QDef{009} & CoRefl \QNeg{RgEucl} & \REFF{corefl 5}{2}	\\\hline
\QDee{010} & Empty \QNeg{LfUnique} & \REFF{empty}{3}	\\\hline
\QDee{011} & Univ LfUnique & \REFF{univ}{12}	\\\hline
\QDef{012} & CoRefl \QNeg{LfUnique} & \REFF{corefl 5}{3}	\\\hline
\QDee{013} & Empty \QNeg{RgUnique} & \REFF{empty}{3}	\\\hline
\QDee{014} & Univ RgUnique & \REFF{univ}{12}	\\\hline
\QDef{015} & CoRefl \QNeg{RgUnique} & \REFF{corefl 5}{4}	\\\hline
\QDee{016} & Empty \QNeg{Sym} & \REFF{empty}{4}	\\\hline
\QDee{017} & Univ \QNeg{Sym} & \REFF{univ}{2}	\\\hline
\QDef{018} & CoRefl \QNeg{Sym} & \REFF{corefl 5}{5}	\\\hline
\QDee{019} & Empty \QNeg{AntiTrans} & \REFF{empty}{5}	\\\hline
\QDee{020} & Univ AntiTrans & \REFF{univ}{9}	\\\hline
\QDee{021} & Empty \QNeg{ASym} & \REFF{empty}{6}	\\\hline
\QDee{022} & Univ ASym & \REFF{univ}{9}	\\\hline
\QDee{023} & Empty Connex & \REFF{empty}{12}	\\\hline
\QDee{024} & Univ \QNeg{Connex} & \REFF{univ}{3}	\\\hline
\QDer{025} & CoRefl Connex & \REFF{corefl 5}{9}	\\\hline
\QDer{026} & LfUnique Connex & \REF{conn 2}	\\\hline
\QDer{027} & RgUnique Connex & \REF{conn 2}	\\\hline
\QDer{028} & AntiTrans Connex & \REF{conn 2}	\\\hline
\QDer{029} & ASym Connex & \REF{irrefl 2}	\\\hline
\QDee{030} & Empty \QNeg{Trans} & \REFF{empty}{7}	\\\hline
\QDee{031} & Univ \QNeg{Trans} & \REFF{univ}{4}	\\\hline
\QDef{032} & CoRefl \QNeg{Trans} & \REFF{corefl 5}{7}	\\\hline
\QDee{033} & Empty \QNeg{SemiOrd1} & \REFF{empty}{8}	\\\hline
\QDee{034} & Univ \QNeg{SemiOrd1} & \REFF{univ}{5}	\\\hline
\QDef{035} & Connex \QNeg{SemiOrd1} & \REF{semiOrd1 1a}	\\\hline
\QDee{036} & Empty \QNeg{Irrefl} & \REFF{empty}{6}	\\\hline
\QDee{037} & Univ Irrefl & \REFF{univ}{9}	\\\hline
\QDef{038} & AntiTrans \QNeg{Irrefl} & \REF{antitrans 1}	\\\hline
\QDef{039} & ASym \QNeg{Irrefl} & \REFF{asym 1b}{1}	\\\hline
\QDer{040} & Connex Irrefl & \REF{irrefl 2}	\\\hline
\QDee{041} & Empty Refl & \REFF{empty}{11}	\\\hline
\QDee{042} & Univ \QNeg{Refl} & \REFF{univ}{6}	\\\hline
\QDer{043} & AntiTrans Refl & \REF{irrefl 2}	\\\hline
\QDer{044} & ASym Refl & \REF{irrefl 2}	\\\hline
\QDef{045} & Connex \QNeg{Refl} & \REF{conn 1}	\\\hline
\QDef{046} & Irrefl Refl & \REF{irrefl 2}	\\\hline
\QDee{047} & Empty \QNeg{QuasiRefl} & \REFF{empty}{1}	\\\hline
\end{tabular}
\hfill
\begin{tabular}{|r|l|c|}
\hline
\QDee{048} & Univ \QNeg{QuasiRefl} & \REFF{univ}{6}	\\\hline
\QDef{049} & CoRefl \QNeg{QuasiRefl} & \REFF{corefl 5}{10}	\\\hline
\QDer{050} & Connex \QNeg{QuasiRefl} & \REF{conn 1}	\\\hline
\QDef{051} & Refl \QNeg{QuasiRefl} & \REF{quasiRefl 2}	\\\hline
\QDee{052} & Empty \QNeg{AntiSym} & \REFF{empty}{6}	\\\hline
\QDee{053} & Univ AntiSym & \REFF{univ}{11}	\\\hline
\QDef{054} & CoRefl \QNeg{AntiSym} & \REFF{corefl 5}{6}	\\\hline
\QDef{055} & ASym \QNeg{AntiSym} & \REFF{asym 1b}{2}	\\\hline
\QDee{056} & Empty SemiConnex & \REFF{empty}{12}	\\\hline
\QDee{057} & Univ \QNeg{SemiConnex} & \REFF{univ}{3}	\\\hline
\QDef{058} & CoRefl SemiConnex & \REFF{corefl 5}{9}	\\\hline
\QDef{059} & LfUnique SemiConnex & \REF{conn 2}	\\\hline
\QDef{060} & RgUnique SemiConnex & \REF{conn 2}	\\\hline
\QDef{061} & AntiTrans SemiConnex & \REF{conn 2}	\\\hline
\QDef{062} & Connex \QNeg{SemiConnex} & \REF{conn 1}	\\\hline
\QDee{063} & Empty \QNeg{IncTrans} & \REFF{empty}{8b}	\\\hline
\QDee{064} & Univ \QNeg{IncTrans} & \REFF{univ}{7}	\\\hline
\QDer{065} & Connex \QNeg{IncTrans} & \REF{incTrans 1}	\\\hline
\QDef{066} & SemiConnex \QNeg{IncTrans} & \REF{incTrans 1}	\\\hline
\QDee{067} & Empty \QNeg{SemiOrd2} & \REFF{empty}{8}	\\\hline
\QDee{068} & Univ \QNeg{SemiOrd2} & \REFF{univ}{5}	\\\hline
\QDer{069} & Connex \QNeg{SemiOrd2} & \REF{semiOrd1 1a}	\\\hline
\QDer{070} & SemiConnex \QNeg{SemiOrd2} & \REF{incTrans 1}	\\\hline
\QDef{071} & IncTrans \QNeg{SemiOrd2} & \REF{incTrans 13}	\\\hline
\QDee{072} & Empty \QNeg{QuasiTrans} & \REFF{empty}{7b}	\\\hline
\QDee{073} & Univ \QNeg{QuasiTrans} & \REFF{univ}{4}	\\\hline
\QDer{074} & CoRefl \QNeg{QuasiTrans} & \REFF{corefl 5}{7}	\\\hline
\QDef{075} & LfEucl \QNeg{QuasiTrans} & \REF{eucl 5}	\\\hline
\QDef{076} & RgEucl \QNeg{QuasiTrans} & \REF{eucl 5}	\\\hline
\QDef{077} & Sym \QNeg{QuasiTrans} & \REF{quasiTrans 4}	\\\hline
\QDef{078} & Trans \QNeg{QuasiTrans} & \REF{quasiTrans 4}	\\\hline
\QDee{079} & Empty \QNeg{Dense} & \REFF{empty}{9}	\\\hline
\QDee{080} & Univ \QNeg{Dense} & \REFF{univ}{6}	\\\hline
\QDer{081} & CoRefl \QNeg{Dense} & \REFF{corefl 5}{8}	\\\hline
\QDeR{082} & LfEucl \QNeg{Dense} & \REFF{den 1}{5}	\\\hline
\QDeR{083} & RgEucl \QNeg{Dense} & \REFF{den 1}{6}	\\\hline
\QDer{084} & Connex \QNeg{Dense} & \REFF{den 1}{8}	\\\hline
\QDer{085} & Refl \QNeg{Dense} & \REFF{den 1}{1}	\\\hline
\QDef{086} & QuasiRefl \QNeg{Dense} & \REFF{den 1}{3}	\\\hline
\QDee{087} & Empty LfSerial & \REFF{empty}{11}	\\\hline
\QDee{088} & Univ \QNeg{LfSerial} & \REFF{univ}{6}	\\\hline
\QDer{089} & Connex \QNeg{LfSerial} & \REF{ser 1}	\\\hline
\QDef{090} & Refl \QNeg{LfSerial} & \REF{ser 1}	\\\hline
\QDee{091} & Empty RgSerial & \REFF{empty}{11}	\\\hline
\QDee{092} & Univ \QNeg{RgSerial} & \REFF{univ}{6}	\\\hline
\QDer{093} & Connex \QNeg{RgSerial} & \REF{ser 1}	\\\hline
\QDef{094} & Refl \QNeg{RgSerial} & \REF{ser 1}	\\\hline
\end{tabular}
\caption{Reported laws for level 2}
\LABEL{Reported laws for level 2}
\end{figure}

\begin{figure}
\footnotesize
\begin{tabular}{|r|l|c|}
\hline
\QDeR{095} & \QNeg{CoRefl} RgEucl LfUnique & \REFF{corefl 6}{3}	\\\hline
\QDeR{096} & \QNeg{CoRefl} LfEucl RgUnique & \REFF{corefl 6}{4}	\\\hline
\QDef{097} & LfEucl RgEucl \QNeg{Sym} & \REF{eucl 4}	\\\hline
\QDef{098} & LfEucl \QNeg{RgEucl} Sym & \REFF{sym 1}{2}	\\\hline
\QDef{099} & \QNeg{LfEucl} RgEucl Sym & \REFF{sym 1}{2}	\\\hline
\QDef{100} & LfUnique \QNeg{RgUnique} Sym & \REFF{sym 1}{4}	\\\hline
\QDef{101} & \QNeg{LfUnique} RgUnique Sym & \REFF{sym 1}{4}	\\\hline
\QDer{102} & \QNeg{Empty} CoRefl AntiTrans & \REF{irrefl 1}	\\\hline
\QDer{103} & \QNeg{Empty} LfEucl AntiTrans & \REF{irrefl 1}	\\\hline
\QDer{104} & \QNeg{Empty} RgEucl AntiTrans & \REF{irrefl 1}	\\\hline
\QDer{105} & \QNeg{Empty} CoRefl ASym & \REF{irrefl 1}	\\\hline
\QDer{106} & \QNeg{Empty} LfEucl ASym & \REF{irrefl 1}	\\\hline
\QDer{107} & \QNeg{Empty} RgEucl ASym & \REF{irrefl 1}	\\\hline
\QDef{108} & \QNeg{Empty} Sym ASym & \REF{sym 2}	\\\hline
\QDer{109} & LfUnique \QNeg{AntiTrans} ASym & \REF{antitrans 3}	\\\hline
\QDer{110} & RgUnique \QNeg{AntiTrans} ASym & \REF{antitrans 3}	\\\hline
\QDer{111} & \QNeg{Univ} LfEucl Connex & \REFF{conn 4}{2}	\\\hline
\QDer{112} & \QNeg{Univ} RgEucl Connex & \REFF{conn 4}{3}	\\\hline
\QDef{113} & \QNeg{Univ} Sym Connex & \REFF{conn 4}{1}	\\\hline
\QDef{114} & LfEucl RgEucl \QNeg{Trans} & \REF{eucl 4}	\\\hline
\QDef{115} & LfEucl LfUnique \QNeg{Trans} & \REF{eucl 9}	\\\hline
\QDef{116} & RgEucl RgUnique \QNeg{Trans} & \REF{eucl 9}	\\\hline
\QDef{117} & \QNeg{LfEucl} Sym Trans & \REF{eucl 1}	\\\hline
\QDer{118} & AntiTrans \QNeg{ASym} Trans & \REFF{asym 1a}{6}	\\\hline
\QDer{119} & AntiTrans \QNeg{ASym} SemiOrd1 & \REFF{asym 1a}{5}	\\\hline
\QDef{120} & LfUnique \QNeg{Trans} SemiOrd1 & \REFF{semiOrd1 12}{1}	\\\hline
\QDef{121} & RgUnique \QNeg{Trans} SemiOrd1 & \REFF{semiOrd1 12}{2}	\\\hline
\QDer{122} & AntiTrans \QNeg{Trans} SemiOrd1 & \REFF{semiOrd1 12}{4}	\\\hline
\QDer{123} & ASym \QNeg{Trans} SemiOrd1 & \REFF{semiOrd1 12}{5}	\\\hline
\QDef{124} & \QNeg{Empty} CoRefl Irrefl & \REFF{irrefl 1}{1}	\\\hline
\QDef{125} & \QNeg{Empty} LfEucl Irrefl & \REFF{irrefl 1}{4}	\\\hline
\QDef{126} & \QNeg{Empty} RgEucl Irrefl & \REFF{irrefl 1}{5}	\\\hline
\QDef{127} & LfUnique \QNeg{AntiTrans} Irrefl & \REF{antitrans 3}	\\\hline
\QDef{128} & RgUnique \QNeg{AntiTrans} Irrefl & \REF{antitrans 3}	\\\hline
\QDef{129} & \QNeg{ASym} Trans Irrefl & \REFF{asym 1a}{2}	\\\hline
\QDer{130} & \QNeg{ASym} SemiOrd1 Irrefl & \REFF{asym 1a}{4}	\\\hline
\QDef{131} & LfEucl \QNeg{RgEucl} Refl & \REF{eucl 2}	\\\hline
\QDef{132} & \QNeg{LfEucl} RgEucl Refl & \REF{eucl 2}	\\\hline
\QDer{133} & \QNeg{CoRefl} LfUnique Refl & \REFF{corefl 6}{4a}	\\\hline
\QDer{134} & \QNeg{CoRefl} RgUnique Refl & \REFF{corefl 6}{4b}	\\\hline
\QDef{135} & CoRefl SemiOrd1 Refl & \REFF{corefl 2}{4}	\\\hline
\QDef{136} & \QNeg{Connex} SemiOrd1 Refl & \REF{semiOrd1 1a}	\\\hline
\QDef{137} & LfEucl RgEucl \QNeg{QuasiRefl} & \REF{eucl 4}	\\\hline
\QDef{138} & LfEucl \QNeg{RgEucl} QuasiRefl & \REF{eucl 2}	\\\hline
\QDef{139} & \QNeg{LfEucl} RgEucl QuasiRefl & \REF{eucl 2}	\\\hline
\QDef{140} & \QNeg{CoRefl} LfUnique QuasiRefl & \REFF{corefl 6}{1}	\\\hline
\QDef{141} & \QNeg{CoRefl} RgUnique QuasiRefl & \REFF{corefl 6}{2}	\\\hline
\end{tabular}
\hfill
\begin{tabular}{|r|l|c|}
\hline
\QDer{142} & \QNeg{Empty} AntiTrans QuasiRefl & \REF{irrefl 1}	\\\hline
\QDer{143} & \QNeg{Empty} ASym QuasiRefl & \REF{irrefl 1}	\\\hline
\QDef{144} & \QNeg{Empty} Irrefl QuasiRefl & \REFF{irrefl 1}{2}	\\\hline
\QDef{145} & LfEucl LfUnique \QNeg{AntiSym} & \REF{eucl 10}	\\\hline
\QDef{146} & LfEucl \QNeg{LfUnique} AntiSym & \REF{eucl 10}	\\\hline
\QDef{147} & RgEucl RgUnique \QNeg{AntiSym} & \REF{eucl 10}	\\\hline
\QDef{148} & RgEucl \QNeg{RgUnique} AntiSym & \REF{eucl 10}	\\\hline
\QDef{149} & \QNeg{CoRefl} Sym AntiSym & \REFF{corefl 6}{5}	\\\hline
\QDer{150} & AntiTrans \QNeg{ASym} AntiSym & \REFF{asym 1a}{3}	\\\hline
\QDef{151} & LfUnique Trans \QNeg{AntiSym} & \REF{uniq 1}	\\\hline
\QDef{152} & RgUnique Trans \QNeg{AntiSym} & \REF{uniq 1}	\\\hline
\QDef{153} & \QNeg{ASym} Irrefl AntiSym & \REFF{asym 1a}{1}	\\\hline
\QDef{154} & LfEucl \QNeg{Trans} SemiConnex & \REF{eucl 6}	\\\hline
\QDef{155} & RgEucl \QNeg{Trans} SemiConnex & \REF{eucl 6}	\\\hline
\QDer{156} & LfEucl \QNeg{SemiOrd1} SemiConnex & \REFF{semiOrd1 3a}{2}	\\\hline
\QDer{157} & RgEucl \QNeg{SemiOrd1} SemiConnex & \REFF{semiOrd1 3a}{3}	\\\hline
\QDer{158} & Trans \QNeg{SemiOrd1} SemiConnex & \REFF{semiOrd1 3a}{1}	\\\hline
\QDer{159} & \QNeg{Connex} Refl SemiConnex & \REF{conn 1}	\\\hline
\QDef{160} & \QNeg{Connex} QuasiRefl SemiConnex & \REF{conn 1}	\\\hline
\QDef{161} & \QNeg{Empty} CoRefl IncTrans & \REFF{corefl 3}{2}	\\\hline
\QDef{162} & LfEucl \QNeg{Trans} IncTrans & \REF{incTrans 3}	\\\hline
\QDef{163} & RgEucl \QNeg{Trans} IncTrans & \REF{incTrans 3}	\\\hline
\QDef{164} & LfEucl \QNeg{SemiOrd1} IncTrans & \REF{incTrans 3}	\\\hline
\QDef{165} & RgEucl \QNeg{SemiOrd1} IncTrans & \REF{incTrans 3}	\\\hline
\QDef{166} & Trans \QNeg{SemiOrd1} IncTrans & \REFF{semiOrd1 3a}{6}	\\\hline
\QDef{167} & \QNeg{Connex} Refl IncTrans & \REF{incTrans 2}	\\\hline
\QDef{168} & \QNeg{SemiOrd1} QuasiRefl IncTrans & \REFF{semiOrd1 3a}{5}	\\\hline
\QDef{169} & \QNeg{Empty} CoRefl SemiOrd2 & \REFF{corefl 3}{1}	\\\hline
\QDef{170} & LfEucl \QNeg{Trans} SemiOrd2 & \REF{semiOrd2 4}	\\\hline
\QDef{171} & RgEucl \QNeg{Trans} SemiOrd2 & \REF{semiOrd2 4}	\\\hline
\QDef{172} & AntiTrans Trans \QNeg{SemiOrd2} & \REF{antitrans 4}	\\\hline
\QDef{173} & LfEucl \QNeg{SemiOrd1} SemiOrd2 & \REF{semiOrd2 5}	\\\hline
\QDef{174} & RgEucl \QNeg{SemiOrd1} SemiOrd2 & \REF{semiOrd2 5}	\\\hline
\QDef{175} & \QNeg{Connex} Refl SemiOrd2 & \REF{semiOrd1 1a}	\\\hline
\QDef{176} & \QNeg{SemiOrd1} QuasiRefl SemiOrd2 & \REF{semiOrd2 10}	\\\hline
\QDeR{177} & LfEucl \QNeg{IncTrans} SemiOrd2 & \REFF{incTrans 14}{5}	\\\hline
\QDeR{178} & RgEucl \QNeg{IncTrans} SemiOrd2 & \REFF{incTrans 14}{6}	\\\hline
\QDef{179} & Sym \QNeg{IncTrans} SemiOrd2 & \REFF{incTrans 14}{3}	\\\hline
\QDef{180} & QuasiRefl \QNeg{IncTrans} SemiOrd2 & \REFF{incTrans 14}{1}	\\\hline
\QDer{181} & ASym \QNeg{Trans} QuasiTrans & \REF{quasiTrans 5}	\\\hline
\QDef{182} & \QNeg{Trans} AntiSym QuasiTrans & \REF{quasiTrans 5}	\\\hline
\QDef{183} & \QNeg{LfEucl} LfUnique Dense & \REF{eucl 11b}	\\\hline
\QDef{184} & \QNeg{RgEucl} RgUnique Dense & \REF{eucl 11b}	\\\hline
\QDef{185} & \QNeg{Empty} AntiTrans Dense & \REF{den 3}	\\\hline
\QDee{186} & \QNeg{Empty} ASym Dense & \REF{exm 186}	\\\hline
\QDef{187} & Sym SemiOrd1 \QNeg{Dense} & \REFF{den 1}{7}	\\\hline
\QDef{188} & Sym SemiConnex \QNeg{Dense} & \REFF{den 1}{9}	\\\hline
\end{tabular}
\caption{Reported laws for level 3 (a)}
\LABEL{Reported laws for level 3 (a)}
\end{figure}

\begin{figure}
\footnotesize
\begin{tabular}{|r|l|c|}
\hline
\QDef{189} & \QNeg{LfEucl} RgEucl LfSerial & \REF{eucl 3}	\\\hline
\QDef{190} & \QNeg{LfUnique} RgUnique LfSerial & \REF{uniq 2}	\\\hline
\QDef{191} & AntiTrans Trans LfSerial & \REF{antitrans 4}	\\\hline
\QDef{192} & ASym Trans LfSerial & \REF{asym 3}	\\\hline
\QDef{193} & CoRefl SemiOrd1 LfSerial & \REFF{corefl 2}{4}	\\\hline
\QDef{194} & RgUnique SemiOrd1 LfSerial & \REFF{semiOrd1 7}{2}	\\\hline
\QDef{195} & CoRefl \QNeg{Refl} LfSerial & \REF{ser 3}	\\\hline
\QDeR{196} & RgEucl \QNeg{Refl} LfSerial & \REF{ser 3}	\\\hline
\QDef{197} & \QNeg{Refl} QuasiRefl LfSerial & \REF{ser 3}	\\\hline
\QDef{198} & LfEucl SemiConnex \QNeg{LfSerial} & \REF{ser 4}	\\\hline
\QDef{199} & Sym SemiConnex \QNeg{LfSerial} & \REF{ser 4}	\\\hline
\QDef{200} & RgUnique IncTrans LfSerial & \REF{incTrans 12}	\\\hline
\QDef{201} & LfEucl \QNeg{RgEucl} RgSerial & \REF{eucl 3}	\\\hline
\QDef{202} & LfUnique \QNeg{RgUnique} RgSerial & \REF{uniq 2}	\\\hline
\end{tabular}
\hfill
\begin{tabular}{|r|l|c|}
\hline
\QDef{203} & AntiTrans Trans RgSerial & \REF{antitrans 4}	\\\hline
\QDef{204} & ASym Trans RgSerial & \REF{asym 3}	\\\hline
\QDef{205} & CoRefl SemiOrd1 RgSerial & \REFF{corefl 2}{4}	\\\hline
\QDef{206} & LfUnique SemiOrd1 RgSerial & \REFF{semiOrd1 7}{1}	\\\hline
\QDef{207} & CoRefl \QNeg{Refl} RgSerial & \REF{ser 3}	\\\hline
\QDeR{208} & LfEucl \QNeg{Refl} RgSerial & \REF{ser 3}	\\\hline
\QDef{209} & \QNeg{Refl} QuasiRefl RgSerial & \REF{ser 3}	\\\hline
\QDef{210} & RgEucl SemiConnex \QNeg{RgSerial} & \REF{ser 4}	\\\hline
\QDef{211} & Sym SemiConnex \QNeg{RgSerial} & \REF{ser 4}	\\\hline
\QDef{212} & LfUnique IncTrans RgSerial & \REF{incTrans 12}	\\\hline
\QDef{213} & LfUnique \QNeg{LfSerial} RgSerial & \REF{uniq 2}	\\\hline
\QDef{214} & RgUnique LfSerial \QNeg{RgSerial} & \REF{uniq 2}	\\\hline
\QDef{215} & Sym LfSerial \QNeg{RgSerial} & \REFF{sym 1}{3}	\\\hline
\QDef{216} & Sym \QNeg{LfSerial} RgSerial & \REFF{sym 1}{3}	\\\hline
\end{tabular}
\caption{Reported laws for level 3 (b)}
\LABEL{Reported laws for level 3 (b)}
\end{figure}

\begin{figure}
\begin{center}
\footnotesize
\begin{tabular}{|r|l|c|}
\hline
\QDef{217} & \QNeg{LfEucl} LfUnique \QNeg{AntiTrans} SemiOrd1 & \REF{eucl 11b}	\\\hline
\QDef{218} & \QNeg{RgEucl} RgUnique \QNeg{AntiTrans} SemiOrd1 & \REF{eucl 11b}	\\\hline
\QDef{219} & \QNeg{LfEucl} Sym SemiOrd1 QuasiRefl & \REFF{semiOrd1 4}{1}	\\\hline
\QDef{220} & \QNeg{Empty} LfUnique RgUnique IncTrans & \REFF{incTrans 9}{1}	\\\hline
\QDer{221} & \QNeg{LfEucl} LfUnique \QNeg{AntiTrans} IncTrans & \REFF{eucl 11b}{3}	\\\hline
\QDer{222} & \QNeg{RgEucl} RgUnique \QNeg{AntiTrans} IncTrans & \REF{eucl 11b}	\\\hline
\QDef{223} & \QNeg{Empty} \QNeg{Connex} QuasiRefl IncTrans & \REF{incTrans 6}	\\\hline
\QDef{224} & \QNeg{LfEucl} LfUnique \QNeg{AntiTrans} SemiOrd2 & \REFF{eucl 11b}{5}	\\\hline
\QDef{225} & \QNeg{RgEucl} RgUnique \QNeg{AntiTrans} SemiOrd2 & \REF{eucl 11b}	\\\hline
\QDef{226} & LfUnique RgUnique \QNeg{ASym} SemiOrd2 & \REF{semiOrd2 8}	\\\hline
\QDee{227} & Trans \QNeg{SemiOrd1} \QNeg{AntiSym} SemiOrd2 & \REF{exm 227}	\\\hline
\QDer{228} & LfUnique \QNeg{ASym} IncTrans \QNeg{QuasiTrans} & \REFF{incTrans 9}{3}	\\\hline
\QDer{229} & RgUnique \QNeg{ASym} IncTrans \QNeg{QuasiTrans} & \REFF{incTrans 9}{3}	\\\hline
\QDef{230} & LfUnique \QNeg{ASym} SemiOrd2 \QNeg{QuasiTrans} & \REFF{semiOrd2 7b}{3}	\\\hline
\QDef{231} & RgUnique \QNeg{ASym} SemiOrd2 \QNeg{QuasiTrans} & \REF{semiOrd2 7b}	\\\hline
\QDef{232} & Sym \QNeg{AntiTrans} IncTrans \QNeg{Dense} & \REF{incTrans 7}	\\\hline
\QDef{233} & Trans \QNeg{SemiOrd1} SemiOrd2 Dense & \REFF{semiOrd1 3a}{4}	\\\hline
\QDef{234} & Trans \QNeg{IncTrans} SemiOrd2 Dense & \REFF{incTrans 14}{4}	\\\hline
\QDef{235} & LfUnique Sym AntiTrans LfSerial & \REF{uniq 3}	\\\hline
\QDef{236} & \QNeg{LfEucl} LfUnique Trans LfSerial & \REFF{eucl 11b}{2}	\\\hline
\QDef{237} & \QNeg{Empty} LfEucl IncTrans \QNeg{LfSerial} & \REF{incTrans 3}	\\\hline
\QDef{238} & \QNeg{Empty} Sym IncTrans \QNeg{LfSerial} & \REF{incTrans 8}	\\\hline
\QDer{239} & LfUnique \QNeg{ASym} IncTrans \QNeg{LfSerial} & \REFF{incTrans 9}{4}	\\\hline
\QDef{240} & LfUnique ASym IncTrans LfSerial & \REFF{incTrans 9}{5}	\\\hline
\QDef{241} & LfUnique SemiOrd1 \QNeg{IncTrans} LfSerial & \REFF{semiOrd1 6}{1}	\\\hline
\QDef{242} & LfUnique \QNeg{ASym} SemiOrd2 \QNeg{LfSerial} & \REFF{semiOrd2 7b}{4}	\\\hline
\QDee{243} & RgUnique ASym \QNeg{SemiOrd2} LfSerial & \REF{exm 243,251}	\\\hline
\QDef{244} & RgUnique \QNeg{Sym} QuasiTrans LfSerial & \REFF{quasiTrans 8}{1}	\\\hline
\QDef{245} & \QNeg{RgEucl} RgUnique Trans RgSerial & \REF{eucl 11b}	\\\hline
\QDef{246} & \QNeg{Empty} RgEucl IncTrans \QNeg{RgSerial} & \REF{incTrans 3}	\\\hline
\QDer{247} & RgUnique \QNeg{ASym} IncTrans \QNeg{RgSerial} & \REF{incTrans 9}	\\\hline
\QDef{248} & RgUnique ASym IncTrans RgSerial & \REF{incTrans 9}	\\\hline
\QDef{249} & RgUnique SemiOrd1 \QNeg{IncTrans} RgSerial & \REFF{semiOrd1 6}{4}	\\\hline
\QDef{250} & RgUnique \QNeg{ASym} SemiOrd2 \QNeg{RgSerial} & \REF{semiOrd2 7b}	\\\hline
\QDee{251} & ASym \QNeg{SemiOrd2} LfSerial RgSerial & \REF{exm 243,251}	\\\hline
\end{tabular}
\caption{Reported laws for level 4}
\LABEL{Reported laws for level 4}
\end{center}
\end{figure}

\begin{figure}
\begin{center}
\footnotesize
\begin{tabular}{|r|l|c|}
\hline
\QDee{252} & LfUnique RgUnique \QNeg{AntiTrans} \QNeg{AntiSym} \QNeg{QuasiTrans} & \REF{exm 252,255,257}	\\\hline
\QDef{253} & LfEucl Trans SemiOrd1 \QNeg{IncTrans} LfSerial & \REF{incTrans 11}	\\\hline
\QDef{254} & LfUnique RgUnique \QNeg{Trans} SemiOrd2 \QNeg{LfSerial} & \REF{semiOrd2 6}	\\\hline
\QDee{255} & LfUnique \QNeg{AntiTrans} \QNeg{AntiSym} \QNeg{QuasiTrans} \QNeg{LfSerial} & \REF{exm 252,255,257}	\\\hline
\QDef{256} & RgEucl Trans SemiOrd1 \QNeg{IncTrans} RgSerial & \REF{incTrans 11}	\\\hline
\QDee{257} & RgUnique \QNeg{AntiTrans} \QNeg{AntiSym} \QNeg{QuasiTrans} \QNeg{RgSerial} & \REF{exm 252,255,257}	\\\hline
\QDee{258} & Trans \QNeg{SemiOrd1} SemiOrd2 LfSerial RgSerial & \REF{exm 258}	\\\hline
\end{tabular}
\caption{Reported laws for level 5}
\LABEL{Reported laws for level 5}
\end{center}
\end{figure}

\begin{figure}
\begin{center}
\footnotesize
\begin{tabular}{|r|l|c|}
\hline
\QDee{259} & \QNeg{Trans} SemiOrd1 AntiSym \QNeg{IncTrans} \QNeg{Dense} LfSerial & \REF{exm 259}	\\\hline
\QDee{260} & \QNeg{Trans} SemiOrd1 AntiSym \QNeg{IncTrans} \QNeg{Dense} RgSerial & \REF{exm 259}	\\\hline
\QDee{261} & \QNeg{Trans} SemiOrd1 AntiSym \QNeg{IncTrans} LfSerial RgSerial & \REF{exm 261}	\\\hline
\QDee{262} & \QNeg{ASym} Trans \QNeg{SemiOrd1} SemiOrd2 \QNeg{LfSerial} \QNeg{RgSerial} & \REF{exm 262}	\\\hline
\QDee{263} & Trans \QNeg{AntiSym} \QNeg{IncTrans} SemiOrd2 \QNeg{LfSerial} \QNeg{RgSerial} & \REF{exm 263}	\\\hline
\QDee{264} & Trans \QNeg{AntiSym} \QNeg{IncTrans} SemiOrd2 LfSerial RgSerial & \REF{exm 264}	\\\hline
\QDef{265} & AntiTrans \QNeg{IncTrans} SemiOrd2 QuasiTrans LfSerial RgSerial & \REFF{incTrans 14}{7}	\\\hline
\QDee{266} & Trans \QNeg{SemiOrd1} \QNeg{AntiSym} \QNeg{Dense} LfSerial RgSerial & \REF{exm 266}	\\\hline
\QDee{267} & \QNeg{Trans} SemiOrd1 AntiSym \QNeg{Dense} LfSerial RgSerial & \REF{exm 267}	\\\hline
\QDee{268} & Trans \QNeg{AntiSym} \QNeg{SemiConnex} \QNeg{Dense} LfSerial RgSerial & \REF{exm 268}	\\\hline
\QDee{269} & SemiOrd1 SemiConnex \QNeg{QuasiTrans} \QNeg{Dense} \QNeg{LfSerial} \QNeg{RgSerial} & \REF{exm 269}	\\\hline
\QDee{270} & Irrefl SemiConnex \QNeg{QuasiTrans} Dense \QNeg{LfSerial} \QNeg{RgSerial} & \REF{exm 270}	\\\hline
\end{tabular}
\caption{Reported laws for level 6}
\LABEL{Reported laws for level 6}
\end{center}
\end{figure}

\begin{figure}
\begin{center}
\footnotesize
\begin{tabular}{|r|l|c|}
\hline
\QDee{271} & Trans \QNeg{AntiSym} \QNeg{SemiConnex} IncTrans \QNeg{Dense} & \\ & \QNeg{LfSerial} \QNeg{RgSerial} & \REF{exm 271}	\\\hline
\QDee{272} & SemiOrd1 \QNeg{IncTrans} SemiOrd2 \QNeg{QuasiTrans} \QNeg{Dense} & \\ & LfSerial \QNeg{RgSerial} & \REF{exm 272}	\\\hline
\QDee{273} & SemiOrd1 \QNeg{IncTrans} SemiOrd2 \QNeg{QuasiTrans} \QNeg{Dense} & \\ & \QNeg{LfSerial} RgSerial & \REF{exm 272}	\\\hline
\end{tabular}
\caption{Reported laws for level 7}
\LABEL{Reported laws for level 7}
\end{center}
\end{figure}

\begin{figure}
\begin{center}
\footnotesize
\begin{tabular}{|r|l|c|}
\hline
\QDee{274} & \QNeg{Trans} SemiOrd1 \QNeg{IncTrans} SemiOrd2 QuasiTrans \QNeg{Dense} & \\ & LfSerial RgSerial & \REF{exm 274}	\\\hline
\end{tabular}
\caption{Reported laws for level 8}
\LABEL{Reported laws for level 8}
\end{center}
\end{figure}

{

\renewcommand{\t}[1]{\begin{tabular}{@{}l@{}}#1\end{tabular}}
\newcommand{\2}[1]{\hyperlink{QL#1}{\textcolor{coQii}{\sf #1}}}
\newcommand{\3}[1]{\hyperlink{QL#1}{\textcolor{coQiii}{\sf #1}}}
\newcommand{\4}[1]{\hyperlink{QL#1}{\textcolor{coQiv}{\sf #1}}}
\newcommand{\5}[1]{\hyperlink{QL#1}{\textcolor{coQv}{\sf #1}}}
\newcommand{\6}[1]{\hyperlink{QL#1}{\textcolor{coQvi}{\sf #1}}}

%%%%%%%% start of table + - %%%%%%%%
%%%%%%%% rows (em un cr le re lu ru sy at as co tr) %%%%%%%%
%%%%%%%% cols (em un cr le re lu ru sy at as co tr) %%%%%%%%

\begin{figure}
\begin{center}
\scriptsize
\begin{tabular}{|l@{$\;$}c|*{12}{|@{$\;$}c@{$\;$}}}
\hline
\multicolumn{2}{|r||@{$\;$}}{$_A \setminus ^B$} 
	& em & un & cr & le & re & lu & ru & sy & at & as & co & tr \\
\hline
\hline
Empty & em%
	{%
	} & \t{ % +em -em
		\x
	} & \t{ % +em -un
	} & \t{ % +em -cr
		\2{002}\\	% +em -cr
	} & \t{ % +em -le
		\2{004}\\	% +em -le
	} & \t{ % +em -re
		\2{007}\\	% +em -re
	} & \t{ % +em -lu
		\2{010}\\	% +em -lu
	} & \t{ % +em -ru
		\2{013}\\	% +em -ru
	} & \t{ % +em -sy
		\2{016}\\	% +em -sy
	} & \t{ % +em -at
		\2{019}\\	% +em -at
	} & \t{ % +em -as
		\2{021}\\	% +em -as
	} & \t{ % +em -co
	} & \t{ % +em -tr
		\2{030}\\	% +em -tr
	}
\\\hline
Univ & un%
	{%
	} & \t{ % +un -em
	} & \t{ % +un -un
		\x
	} & \t{ % +un -cr
	} & \t{ % +un -le
		\2{005}\\	% +un -le
	} & \t{ % +un -re
		\2{008}\\	% +un -re
	} & \t{ % +un -lu
	} & \t{ % +un -ru
	} & \t{ % +un -sy
		\2{017}\\	% +un -sy
	} & \t{ % +un -at
	} & \t{ % +un -as
	} & \t{ % +un -co
		\2{024}\\	% +un -co
	} & \t{ % +un -tr
		\2{031}\\	% +un -tr
	}
\\\hline
CoRefl & cr%
	{%
	} & \t{ % +cr -em
		%\3{102},%	% -em +cr +at
		%\3{105}\\	% -em +cr +as
		\3{124},%	% -em +cr +ir
		\3{161}\\	% -em +cr +it
		\3{169}\\	% -em +cr +s2
	} & \t{ % +cr -un
	} & \t{ % +cr -cr
		\x
	} & \t{ % +cr -le
		\2{006}\\	% +cr -le
	} & \t{ % +cr -re
		\2{009}\\	% +cr -re
	} & \t{ % +cr -lu
		\2{012}\\	% +cr -lu
	} & \t{ % +cr -ru
		\2{015}\\	% +cr -ru
	} & \t{ % +cr -sy
		\2{018}\\	% +cr -sy
	} & \t{ % +cr -at
	} & \t{ % +cr -as
	} & \t{ % +cr -co
	} & \t{ % +cr -tr
		\2{032}\\	% +cr -tr
	}
\\\hline
LfEucl & le%
	{%
	} & \t{ % +le -em
		%\3{103},%	% -em +le +at
		%\3{106}\\	% -em +le +as
		\3{125},%	% -em +le +ir
		\4{237}\\	% -em +le +it -ls
	} & \t{ % +le -un
		%\3{111}\\	% -un +le +co
	} & \t{ % +le -cr
		\3{096}\\	% -cr +le +ru		(+)
	} & \t{ % +le -le
		\x
	} & \t{ % +le -re
		\3{098}\\	% +le -re +sy
		\3{131}\\	% +le -re +rf
		\3{138}\\	% +le -re +qr
		\3{201}\\	% +le -re +rs
	} & \t{ % +le -lu
		\3{146}\\	% +le -lu +an
	} & \t{ % +le -ru
	} & \t{ % +le -sy
		\3{097}\\	% +le +re -sy
	} & \t{ % +le -at
	} & \t{ % +le -as
	} & \t{ % +le -co
	} & \t{ % +le -tr
		\3{114},%	% +le +re -tr
		\3{115}\\	% +le +lu -tr
		\3{154},%	% +le -tr +sc
		\3{162}\\	% +le -tr +it
		\3{170}\\	% +le -tr +s2
	}
\\\hline
RgEucl & re%
	{%
	} & \t{ % +re -em
		%\3{104},%	% -em +re +at
		%\3{107}\\	% -em +re +as
		\3{126},%	% -em +re +ir
		\4{246}\\	% -em +re +it -rs
	} & \t{ % +re -un
		%\3{112}\\	% -un +re +co
	} & \t{ % +re -cr
		\3{095}\\	% -cr +re +lu		(+)
	} & \t{ % +re -le
		\3{099}\\	% -le +re +sy
		\3{132}\\	% -le +re +rf
		\3{139}\\	% -le +re +qr
		\3{189}\\	% -le +re +ls
	} & \t{ % +re -re
		\x
	} & \t{ % +re -lu
	} & \t{ % +re -ru
		\3{148}\\	% +re -ru +an
	} & \t{ % +re -sy
		\3{097}\\	% +le +re -sy
	} & \t{ % +re -at
	} & \t{ % +re -as
	} & \t{ % +re -co
	} & \t{ % +re -tr
		\3{114},%	% +le +re -tr
		\3{116}\\	% +re +ru -tr
		\3{155},%	% +re -tr +sc
		\3{163}\\	% +re -tr +it
		\3{171}\\	% +re -tr +s2
	}
\\\hline
LfUnique & lu%
	{%
	} & \t{ % +lu -em
		\4{220}\\	% -em +lu +ru +it
	} & \t{ % +lu -un
	} & \t{ % +lu -cr
		\3{095}\\	% -cr +re +lu		(+)
		%\3{133}\\	% -cr +lu +rf
		\3{140}\\	% -cr +lu +qr
	} & \t{ % +lu -le
		\3{183}\\	% -le +lu +de
		\4{217}\\	% -le +lu -at +s1
		%\4{221}\\	% -le +lu -at +it
		\4{224}\\	% -le +lu -at +s2
		\4{236}\\	% -le +lu +tr +ls
	} & \t{ % +lu -re
	} & \t{ % +lu -lu
		\x
	} & \t{ % +lu -ru
		\3{100}\\	% +lu -ru +sy
		\3{202}\\	% +lu -ru +rs
	} & \t{ % +lu -sy
	} & \t{ % +lu -at
		%\3{109},%	% +lu -at +as
		\3{127},%	% +lu -at +ir
		\4{217}\\	% -le +lu -at +s1
		%\4{221},%	% -le +lu -at +it
		\4{224}\\	% -le +lu -at +s2
	} & \t{ % +lu -as
		\4{226},%	% +lu +ru -as +s2
		%\4{228}\\	% +lu -as +it -qt
		\4{230}\\	% +lu -as +s2 -qt
		%\4{239}\\	% +lu -as +it -ls
		\4{242}\\	% +lu -as +s2 -ls
	} & \t{ % +lu -co
	} & \t{ % +lu -tr
		\3{115},%	% +le +lu -tr
		\3{120}\\	% +lu -tr +s1
		\5{254}\\	% +lu +ru -tr +s2 -ls
	}
\\\hline
RgUnique & ru%
	{%
	} & \t{ % +ru -em
		\4{220}\\	% -em +lu +ru +it
	} & \t{ % +ru -un
	} & \t{ % +ru -cr
		\3{096}\\	% -cr +le +ru		(+)
		%\3{134}\\	% -cr +ru +rf
		\3{141}\\	% -cr +ru +qr
	} & \t{ % +ru -le
	} & \t{ % +ru -re
		\3{184}\\	% -re +ru +de
		\4{218}\\	% -re +ru -at +s1
		%\4{222}\\	% -re +ru -at +it
		\4{225}\\	% -re +ru -at +s2
		\4{245}\\	% -re +ru +tr +rs
	} & \t{ % +ru -lu
		\3{101}\\	% -lu +ru +sy
		\3{190}\\	% -lu +ru +ls
	} & \t{ % +ru -ru
		\x
	} & \t{ % +ru -sy
		\4{244}\\	% +ru -sy +qt +ls
	} & \t{ % +ru -at
		%\3{110},%	% +ru -at +as
		\3{128},%	% +ru -at +ir
		\4{218}\\	% -re +ru -at +s1
		%\4{222},%	% -re +ru -at +it
		\4{225}\\	% -re +ru -at +s2
	} & \t{ % +ru -as
		\4{226},%	% +lu +ru -as +s2
		%\4{229}\\	% +ru -as +it -qt
		\4{231}\\	% +ru -as +s2 -qt
		%\4{247}\\	% +ru -as +it -rs
		\4{250}\\	% +ru -as +s2 -rs
	} & \t{ % +ru -co
	} & \t{ % +ru -tr
		\3{116},%	% +re +ru -tr
		\3{121}\\	% +ru -tr +s1
		\5{254}\\	% +lu +ru -tr +s2 -ls
	}
\\\hline
Sym & sy%
	{%
	} & \t{ % +sy -em
		\3{108},%	% -em +sy +as
		\4{238}\\	% -em +sy +it -ls
	} & \t{ % +sy -un
		\3{113}\\	% -un +sy +co
	} & \t{ % +sy -cr
		\3{149}\\	% -cr +sy +an
	} & \t{ % +sy -le
		\3{099}\\	% -le +re +sy
		\3{117}\\	% -le +sy +tr
		\4{219}\\	% -le +sy +s1 +qr
	} & \t{ % +sy -re
		\3{098}\\	% +le -re +sy
	} & \t{ % +sy -lu
		\3{101}\\	% -lu +ru +sy
	} & \t{ % +sy -ru
		\3{100}\\	% +lu -ru +sy
	} & \t{ % +sy -sy
		\x
	} & \t{ % +sy -at
		\4{232}\\	% +sy -at +it -de
	} & \t{ % +sy -as
	} & \t{ % +sy -co
	} & \t{ % +sy -tr
	}
\\\hline
AntiTrans & at%
	{%
	} & \t{ % +at -em
		%\3{102},%	% -em +cr +at
		%\3{103}\\	% -em +le +at
		%\3{104},%	% -em +re +at
		%\3{142}\\	% -em +at +qr
		\3{185}\\	% -em +at +de
	} & \t{ % +at -un
	} & \t{ % +at -cr
	} & \t{ % +at -le
	} & \t{ % +at -re
	} & \t{ % +at -lu
	} & \t{ % +at -ru
	} & \t{ % +at -sy
	} & \t{ % +at -at
		\x
	} & \t{ % +at -as
		%\3{118},%	% +at -as +tr
		%\3{119}\\	% +at -as +s1
		%\3{150}\\	% +at -as +an
	} & \t{ % +at -co
	} & \t{ % +at -tr
		%\3{122}\\	% +at -tr +s1
	}
\\\hline
ASym & as%
	{%
	} & \t{ % +as -em
		%\3{105},%	% -em +cr +as
		%\3{106}\\	% -em +le +as
		%\3{107},%	% -em +re +as
		\3{108}\\	% -em +sy +as
		%\3{143}\\	% -em +as +qr
	} & \t{ % +as -un
	} & \t{ % +as -cr
	} & \t{ % +as -le
	} & \t{ % +as -re
	} & \t{ % +as -lu
	} & \t{ % +as -ru
	} & \t{ % +as -sy
	} & \t{ % +as -at
		%\3{109},%	% +lu -at +as
		%\3{110}\\	% +ru -at +as
	} & \t{ % +as -as
		\x
	} & \t{ % +as -co
	} & \t{ % +as -tr
		%\3{123}%	% +as -tr +s1
		%\3{181}\\	% +as -tr +qt
	}
\\\hline
Connex & co%
	{%
	} & \t{ % +co -em
	} & \t{ % +co -un
		%\3{111}\\	% -un +le +co
		%\3{112}\\	% -un +re +co
		\3{113}\\	% -un +sy +co
	} & \t{ % +co -cr
	} & \t{ % +co -le
	} & \t{ % +co -re
	} & \t{ % +co -lu
	} & \t{ % +co -ru
	} & \t{ % +co -sy
	} & \t{ % +co -at
	} & \t{ % +co -as
	} & \t{ % +co -co
		\x
	} & \t{ % +co -tr
	}
\\\hline
Trans & tr%
	{%
	} & \t{ % +tr -em
	} & \t{ % +tr -un
	} & \t{ % +tr -cr
	} & \t{ % +tr -le
		\3{117}\\	% -le +sy +tr
		\4{236}\\	% -le +lu +tr +ls
	} & \t{ % +tr -re
		\4{245}\\	% -re +ru +tr +rs
	} & \t{ % +tr -lu
	} & \t{ % +tr -ru
	} & \t{ % +tr -sy
	} & \t{ % +tr -at
	} & \t{ % +tr -as
		%\3{118}\\	% +at -as +tr
		\3{129}\\	% -as +tr +ir
	} & \t{ % +tr -co
	} & \t{ % +tr -tr
		\x

%%%%%%%% end of table + - %%%%%%%%
%%%%%%%% rows (em un cr le re lu ru sy at as co tr) %%%%%%%%
%%%%%%%% cols (em un cr le re lu ru sy at as co tr) %%%%%%%%

%%%%%%%% start of table + - %%%%%%%%
%%%%%%%% rows (s1 ir rf qr an sc it s2 qt de ls rs) %%%%%%%%
%%%%%%%% cols (em un cr le re lu ru sy at as co tr) %%%%%%%%

	}
\\\hline
\hline
SemiOrd1 & s1%
	{%
	} & \t{ % +s1 -em
	} & \t{ % +s1 -un
	} & \t{ % +s1 -cr
	} & \t{ % +s1 -le
		\4{217}\\	% -le +lu -at +s1
		\4{219}\\	% -le +sy +s1 +qr
	} & \t{ % +s1 -re
		\4{218}\\	% -re +ru -at +s1
	} & \t{ % +s1 -lu
	} & \t{ % +s1 -ru
	} & \t{ % +s1 -sy
	} & \t{ % +s1 -at
		\4{217},%	% -le +lu -at +s1
		\4{218}\\	% -re +ru -at +s1
	} & \t{ % +s1 -as
		%\3{119},%	% +at -as +s1
		%\3{130}\\	% -as +s1 +ir
	} & \t{ % +s1 -co
		\3{136}\\	% -co +s1 +rf
	} & \t{ % +s1 -tr
		\3{120},%	% +lu -tr +s1
		\3{121}\\	% +ru -tr +s1
		%\3{122},%	% +at -tr +s1
		%\3{123}\\	% +as -tr +s1
	}
\\\hline
Irrefl & ir%
	{%
	} & \t{ % +ir -em
		\3{124},%	% -em +cr +ir
		\3{125}\\	% -em +le +ir
		\3{126},%	% -em +re +ir
		\3{144}\\	% -em +ir +qr
	} & \t{ % +ir -un
	} & \t{ % +ir -cr
	} & \t{ % +ir -le
	} & \t{ % +ir -re
	} & \t{ % +ir -lu
	} & \t{ % +ir -ru
	} & \t{ % +ir -sy
	} & \t{ % +ir -at
		\3{127},%	% +lu -at +ir
		\3{128}\\	% +ru -at +ir
	} & \t{ % +ir -as
		\3{129},%	% -as +tr +ir
		%\3{130}\\	% -as +s1 +ir
		\3{153}\\	% -as +ir +an
	} & \t{ % +ir -co
	} & \t{ % +ir -tr
	}
\\\hline
Refl & rf%
	{%
	} & \t{ % +rf -em
	} & \t{ % +rf -un
	} & \t{ % +rf -cr
		%\3{133}\\	% -cr +lu +rf
		%\3{134}\\	% -cr +ru +rf
	} & \t{ % +rf -le
		\3{132}\\	% -le +re +rf
	} & \t{ % +rf -re
		\3{131}\\	% +le -re +rf
	} & \t{ % +rf -lu
	} & \t{ % +rf -ru
	} & \t{ % +rf -sy
	} & \t{ % +rf -at
	} & \t{ % +rf -as
	} & \t{ % +rf -co
		\3{136}\\	% -co +s1 +rf
		%\3{159}\\	% -co +rf +sc
		\3{167}\\	% -co +rf +it
		\3{175}\\	% -co +rf +s2
	} & \t{ % +rf -tr
	}
\\\hline
QuasiRefl & qr%
	{%
	} & \t{ % +qr -em
		%\3{142},%	% -em +at +qr
		%\3{143},%	% -em +as +qr
		\3{144}\\	% -em +ir +qr
		\4{223}\\	% -em -co +qr +it
	} & \t{ % +qr -un
	} & \t{ % +qr -cr
		\3{140}\\	% -cr +lu +qr
		\3{141}\\	% -cr +ru +qr
	} & \t{ % +qr -le
		\3{139}\\	% -le +re +qr
		\4{219}\\	% -le +sy +s1 +qr
	} & \t{ % +qr -re
		\3{138}\\	% +le -re +qr
	} & \t{ % +qr -lu
	} & \t{ % +qr -ru
	} & \t{ % +qr -sy
	} & \t{ % +qr -at
	} & \t{ % +qr -as
	} & \t{ % +qr -co
		\3{160}\\	% -co +qr +sc
		\4{223}\\	% -em -co +qr +it
	} & \t{ % +qr -tr
	}
\\\hline
AntiSym & an%
	{%
	} & \t{ % +an -em
	} & \t{ % +an -un
	} & \t{ % +an -cr
		\3{149}\\	% -cr +sy +an
	} & \t{ % +an -le
	} & \t{ % +an -re
	} & \t{ % +an -lu
		\3{146}\\	% +le -lu +an
	} & \t{ % +an -ru
		\3{148}\\	% +re -ru +an
	} & \t{ % +an -sy
	} & \t{ % +an -at
	} & \t{ % +an -as
		%\3{150},%	% +at -as +an
		\3{153}\\	% -as +ir +an
	} & \t{ % +an -co
	} & \t{ % +an -tr
		\3{182}\\	% -tr +an +qt
	}
\\\hline
SemiConnex & sc%
	\t{\\ \\%
	} & \t{ % +sc -em
	} & \t{ % +sc -un
	} & \t{ % +sc -cr
	} & \t{ % +sc -le
	} & \t{ % +sc -re
	} & \t{ % +sc -lu
	} & \t{ % +sc -ru
	} & \t{ % +sc -sy
	} & \t{ % +sc -at
	} & \t{ % +sc -as
	} & \t{ % +sc -co
		%\3{159},%	% -co +rf +sc
		\3{160}\\	% -co +qr +sc
	} & \t{ % +sc -tr
		\3{154},%	% +le -tr +sc
		\3{155}\\	% +re -tr +sc
	}
\\\hline
IncTrans & it%
	{%
	} & \t{ % +it -em
		\3{161},%	% -em +cr +it
		\4{220}\\	% -em +lu +ru +it
		\4{223},%	% -em -co +qr +it
		\4{237}\\	% -em +le +it -ls
		\4{238},%	% -em +sy +it -ls
		\4{246}\\	% -em +re +it -rs
	} & \t{ % +it -un
	} & \t{ % +it -cr
	} & \t{ % +it -le
		%\4{221}\\	% -le +lu -at +it
	} & \t{ % +it -re
		%\4{222}\\	% -re +ru -at +it
	} & \t{ % +it -lu
	} & \t{ % +it -ru
	} & \t{ % +it -sy
	} & \t{ % +it -at
		%\4{221},%	% -le +lu -at +it
		%\4{222}\\	% -re +ru -at +it
		\4{232}\\	% +sy -at +it -de
	} & \t{ % +it -as
		%\4{228},%	% +lu -as +it -qt
		%\4{229}\\	% +ru -as +it -qt
		%\4{239},%	% +lu -as +it -ls
		%\4{247}\\	% +ru -as +it -rs
	} & \t{ % +it -co
		\3{167}\\	% -co +rf +it
		\4{223}\\	% -em -co +qr +it
	} & \t{ % +it -tr
		\3{162},%	% +le -tr +it
		\3{163}\\	% +re -tr +it
	}
\\\hline
SemiOrd2 & s2%
	{%
	} & \t{ % +s2 -em
		\3{169}\\	% -em +cr +s2
	} & \t{ % +s2 -un
	} & \t{ % +s2 -cr
	} & \t{ % +s2 -le
		\4{224}\\	% -le +lu -at +s2
	} & \t{ % +s2 -re
		\4{225}\\	% -re +ru -at +s2
	} & \t{ % +s2 -lu
	} & \t{ % +s2 -ru
	} & \t{ % +s2 -sy
	} & \t{ % +s2 -at
		\4{224},%	% -le +lu -at +s2
		\4{225}\\	% -re +ru -at +s2
	} & \t{ % +s2 -as
		\4{226},%	% +lu +ru -as +s2
		\4{230}\\	% +lu -as +s2 -qt
		\4{231},%	% +ru -as +s2 -qt
		\4{242}\\	% +lu -as +s2 -ls
		\4{250}\\	% +ru -as +s2 -rs
	} & \t{ % +s2 -co
		\3{175}\\	% -co +rf +s2
	} & \t{ % +s2 -tr
		\3{170},%	% +le -tr +s2
		\3{171}\\	% +re -tr +s2
		\5{254}\\	% +lu +ru -tr +s2 -ls
	}
\\\hline
QuasiTrans & qt%
	{%
	} & \t{ % +qt -em
	} & \t{ % +qt -un
	} & \t{ % +qt -cr
	} & \t{ % +qt -le
	} & \t{ % +qt -re
	} & \t{ % +qt -lu
	} & \t{ % +qt -ru
	} & \t{ % +qt -sy
		\4{244}\\	% +ru -sy +qt +ls
	} & \t{ % +qt -at
	} & \t{ % +qt -as
	} & \t{ % +qt -co
	} & \t{ % +qt -tr
		%\3{181},%	% +as -tr +qt
		\3{182}\\	% -tr +an +qt
	}
\\\hline
Dense & de%
	{%
	} & \t{ % +de -em
		\3{185}\\	% -em +at +de
	} & \t{ % +de -un
	} & \t{ % +de -cr
	} & \t{ % +de -le
		\3{183}\\	% -le +lu +de
	} & \t{ % +de -re
		\3{184}\\	% -re +ru +de
	} & \t{ % +de -lu
	} & \t{ % +de -ru
	} & \t{ % +de -sy
	} & \t{ % +de -at
	} & \t{ % +de -as
	} & \t{ % +de -co
	} & \t{ % +de -tr
	}
\\\hline
LfSerial & ls%
	\t{\\ \\ \\%
	} & \t{ % +ls -em
	} & \t{ % +ls -un
	} & \t{ % +ls -cr
	} & \t{ % +ls -le
		\3{189}\\	% -le +re +ls
		\4{236}\\	% -le +lu +tr +ls
	} & \t{ % +ls -re
	} & \t{ % +ls -lu
		\3{190}\\	% -lu +ru +ls
	} & \t{ % +ls -ru
	} & \t{ % +ls -sy
		\4{244}\\	% +ru -sy +qt +ls
	} & \t{ % +ls -at
	} & \t{ % +ls -as
	} & \t{ % +ls -co
	} & \t{ % +ls -tr
	}
\\\hline
RgSerial & rs%
	\t{\\ \\ \\%
	} & \t{ % +rs -em
	} & \t{ % +rs -un
	} & \t{ % +rs -cr
	} & \t{ % +rs -le
	} & \t{ % +rs -re
		\3{201}\\	% +le -re +rs
		\4{245}\\	% -re +ru +tr +rs
	} & \t{ % +rs -lu
	} & \t{ % +rs -ru
		\3{202}\\	% +lu -ru +rs
	} & \t{ % +rs -sy
	} & \t{ % +rs -at
	} & \t{ % +rs -as
	} & \t{ % +rs -co
	} & \t{ % +rs -tr

	}
\\\hline
\end{tabular}
\caption{Law index lf ($A \Ra B$)}
%\rule{\textwidth}{0.1mm}
\LABEL{Law index + - lf}
\end{center}
\end{figure}

%%%%%%%% end of table + - %%%%%%%%
%%%%%%%% rows (s1 ir rf qr an sc it s2 qt de ls rs) %%%%%%%%
%%%%%%%% cols (em un cr le re lu ru sy at as co tr) %%%%%%%%

%%%%%%% start of table + - %%%%%%%%
%%%%%%% rows (em un cr le re lu ru sy at as co tr) %%%%%%%%
%%%%%%% cols (s1 ir rf qr an sc it s2 qt de ls rs) %%%%%%%%

\begin{figure}
\begin{center}
\scriptsize
\begin{tabular}{lc|*{12}{|@{$\;$}c@{$\;$}}|}
\hline
\multicolumn{2}{r||@{$\;$}}{$_A \setminus ^B$}
	& s1 & ir & rf & qr & an & sc & it & s2 & qt & de & ls & rs \\
\hline
\hline
Empty & em%
	{%
	} & \t{ % +em -s1
		\2{033}\\	% +em -s1
	} & \t{ % +em -ir
		\2{036}\\	% +em -ir
	} & \t{ % +em -rf
	} & \t{ % +em -qr
		\2{047}\\	% +em -qr
	} & \t{ % +em -an
		\2{052}\\	% +em -an
	} & \t{ % +em -sc
	} & \t{ % +em -it
		\2{063}\\	% +em -it
	} & \t{ % +em -s2
		\2{067}\\	% +em -s2
	} & \t{ % +em -qt
		\2{072}\\	% +em -qt
	} & \t{ % +em -de
		\2{079}\\	% +em -de
	} & \t{ % +em -ls
	} & \t{ % +em -rs
	}
\\\hline
Univ & un%
	{%
	} & \t{ % +un -s1
		\2{034}\\	% +un -s1
	} & \t{ % +un -ir
	} & \t{ % +un -rf
		\2{042}\\	% +un -rf
	} & \t{ % +un -qr
		\2{048}\\	% +un -qr
	} & \t{ % +un -an
	} & \t{ % +un -sc
		\2{057}\\	% +un -sc
	} & \t{ % +un -it
		\2{064}\\	% +un -it
	} & \t{ % +un -s2
		\2{068}\\	% +un -s2
	} & \t{ % +un -qt
		\2{073}\\	% +un -qt
	} & \t{ % +un -de
		\2{080}\\	% +un -de
	} & \t{ % +un -ls
		\2{088}\\	% +un -ls
	} & \t{ % +un -rs
		\2{092}\\	% +un -rs
	}
\\\hline
CoRefl & cr%
	{%
	} & \t{ % +cr -s1
	} & \t{ % +cr -ir
	} & \t{ % +cr -rf
		\3{195}\\	% +cr -rf +ls
		\3{207}\\	% +cr -rf +rs
	} & \t{ % +cr -qr
		\2{049}\\	% +cr -qr
	} & \t{ % +cr -an
		\2{054}\\	% +cr -an
	} & \t{ % +cr -sc
	} & \t{ % +cr -it
	} & \t{ % +cr -s2
	} & \t{ % +cr -qt
		%\2{074}\\	% +cr -qt
	} & \t{ % +cr -de
		%\2{081}\\	% +cr -de
	} & \t{ % +cr -ls
	} & \t{ % +cr -rs
	}
\\\hline
LfEucl & le%
	\t{\\ \\ \\ \\%
	} & \t{ % +le -s1
		%\3{156},%	% +le -s1 +sc
		\3{164}\\	% +le -s1 +it
		\3{173}\\	% +le -s1 +s2
	} & \t{ % +le -ir
	} & \t{ % +le -rf
		\3{208}\\	% +le -rf +rs		(+)
	} & \t{ % +le -qr
		\3{137}\\	% +le +re -qr
	} & \t{ % +le -an
		\3{145}\\	% +le +lu -an
	} & \t{ % +le -sc
	} & \t{ % +le -it
		\3{177},%	% +le -it +s2		(+)
		\5{253}\\	% +le +tr +s1 -it +ls
	} & \t{ % +le -s2
	} & \t{ % +le -qt
		\2{075}\\	% +le -qt
	} & \t{ % +le -de
		\2{082}\\	% +le -de		(+)
	} & \t{ % +le -ls
		\3{198}\\	% +le +sc -ls
		\4{237}\\	% -em +le +it -ls
	} & \t{ % +le -rs
	}
\\\hline
RgEucl & re%
	\t{\\ \\ \\ \\%
	} & \t{ % +re -s1
		%\3{157},%	% +re -s1 +sc
		\3{165}\\	% +re -s1 +it
		\3{174}\\	% +re -s1 +s2
	} & \t{ % +re -ir
	} & \t{ % +re -rf
		\3{196}\\	% +re -rf +ls		(+)
	} & \t{ % +re -qr
		\3{137}\\	% +le +re -qr
	} & \t{ % +re -an
		\3{147}\\	% +re +ru -an
	} & \t{ % +re -sc
	} & \t{ % +re -it
		\3{178},%	% +re -it +s2		(+)
		\5{256}\\	% +re +tr +s1 -it +rs
	} & \t{ % +re -s2
	} & \t{ % +re -qt
		\2{076}\\	% +re -qt
	} & \t{ % +re -de
		\2{083}\\	% +re -de		(+)
	} & \t{ % +re -ls
	} & \t{ % +re -rs
		\3{210}\\	% +re +sc -rs
		\4{246}\\	% -em +re +it -rs
	}
\\\hline
LfUnique & lu%
	\t{\\ \\ \\ \\%
	} & \t{ % +lu -s1
	} & \t{ % +lu -ir
	} & \t{ % +lu -rf
	} & \t{ % +lu -qr
	} & \t{ % +lu -an
		\3{145}\\	% +le +lu -an
		\3{151}\\	% +lu +tr -an
	} & \t{ % +lu -sc
	} & \t{ % +lu -it
		\4{241}\\	% +lu +s1 -it +ls
	} & \t{ % +lu -s2
	} & \t{ % +lu -qt
		%\4{228}\\	% +lu -as +it -qt
		\4{230}\\	% +lu -as +s2 -qt
	} & \t{ % +lu -de
	} & \t{ % +lu -ls
		\3{213}\\	% +lu -ls +rs
		%\4{239}\\	% +lu -as +it -ls
		\4{242}\\	% +lu -as +s2 -ls
		\5{254}\\	% +lu +ru -tr +s2 -ls
	} & \t{ % +lu -rs
	}
\\\hline
RgUnique & ru%
	\t{\\ \\ \\ \\%
	} & \t{ % +ru -s1
	} & \t{ % +ru -ir
	} & \t{ % +ru -rf
	} & \t{ % +ru -qr
	} & \t{ % +ru -an
		\3{147}\\	% +re +ru -an
		\3{152}\\	% +ru +tr -an
	} & \t{ % +ru -sc
	} & \t{ % +ru -it
		\4{249}\\	% +ru +s1 -it +rs
	} & \t{ % +ru -s2
	} & \t{ % +ru -qt
		%\4{229}\\	% +ru -as +it -qt
		\4{231}\\	% +ru -as +s2 -qt
	} & \t{ % +ru -de
	} & \t{ % +ru -ls
		\5{254}\\	% +lu +ru -tr +s2 -ls
	} & \t{ % +ru -rs
		\3{214}\\	% +ru +ls -rs
		%\4{247}\\	% +ru -as +it -rs
		\4{250}\\	% +ru -as +s2 -rs
	}
\\\hline
Sym & sy%
	{%
	} & \t{ % +sy -s1
	} & \t{ % +sy -ir
	} & \t{ % +sy -rf
	} & \t{ % +sy -qr
	} & \t{ % +sy -an
	} & \t{ % +sy -sc
	} & \t{ % +sy -it
		\3{179}\\	% +sy -it +s2
	} & \t{ % +sy -s2
	} & \t{ % +sy -qt
		\2{077}\\	% +sy -qt
	} & \t{ % +sy -de
		\3{187}\\	% +sy +s1 -de
		\3{188}\\	% +sy +sc -de
		\4{232}\\	% +sy -at +it -de
	} & \t{ % +sy -ls
		\3{199}\\	% +sy +sc -ls
		\3{216}\\	% +sy -ls +rs
		\4{238}\\	% -em +sy +it -ls
	} & \t{ % +sy -rs
		\3{211}\\	% +sy +sc -rs
		\3{215}\\	% +sy +ls -rs
	}
\\\hline
AntiTrans & at%
	{%
	} & \t{ % +at -s1
	} & \t{ % +at -ir
		\2{038}\\	% +at -ir
	} & \t{ % +at -rf
	} & \t{ % +at -qr
	} & \t{ % +at -an
	} & \t{ % +at -sc
	} & \t{ % +at -it
		\6{265}\\	% +at -it +s2 +qt +ls +rs
	} & \t{ % +at -s2
		\3{172}\\	% +at +tr -s2
	} & \t{ % +at -qt
	} & \t{ % +at -de
	} & \t{ % +at -ls
	} & \t{ % +at -rs
	}
\\\hline
ASym & as%
	\t{%
	} & \t{ % +as -s1
	} & \t{ % +as -ir
		\2{039}\\	% +as -ir
	} & \t{ % +as -rf
	} & \t{ % +as -qr
	} & \t{ % +as -an
		\2{055}\\	% +as -an
	} & \t{ % +as -sc
	} & \t{ % +as -it
	} & \t{ % +as -s2
	} & \t{ % +as -qt
	} & \t{ % +as -de
	} & \t{ % +as -ls
	} & \t{ % +as -rs
	}
\\\hline
Connex & co%
	\t{%
	} & \t{ % +co -s1
		\2{035}\\	% +co -s1
	} & \t{ % +co -ir
	} & \t{ % +co -rf
		\2{045}\\	% +co -rf
	} & \t{ % +co -qr
		%\2{050}\\	% +co -qr
	} & \t{ % +co -an
	} & \t{ % +co -sc
		\2{062}\\	% +co -sc
	} & \t{ % +co -it
		%\2{065}\\	% +co -it
	} & \t{ % +co -s2
		%\2{069}\\	% +co -s2
	} & \t{ % +co -qt
	} & \t{ % +co -de
		%\2{084}\\	% +co -de
	} & \t{ % +co -ls
		%\2{089}\\	% +co -ls
	} & \t{ % +co -rs
		%\2{093}\\	% +co -rs
	}
\\\hline
Trans & tr%
	{%
	} & \t{ % +tr -s1
		%\3{158},%	% +tr -s1 +sc
		\3{166}\\	% +tr -s1 +it
		\4{233}\\	% +tr -s1 +s2 +de
	} & \t{ % +tr -ir
	} & \t{ % +tr -rf
	} & \t{ % +tr -qr
	} & \t{ % +tr -an
		\3{151}\\	% +lu +tr -an
		\3{152}\\	% +ru +tr -an
	} & \t{ % +tr -sc
	} & \t{ % +tr -it
		\4{234},%	% +tr -it +s2 +de
		\5{253}\\	% +le +tr +s1 -it +ls
		\5{256}\\	% +re +tr +s1 -it +rs
	} & \t{ % +tr -s2
		\3{172}\\	% +at +tr -s2
	} & \t{ % +tr -qt
		\2{078}\\	% +tr -qt
	} & \t{ % +tr -de
	} & \t{ % +tr -ls
	} & \t{ % +tr -rs

%%%%%%%% end of table + - %%%%%%%%
%%%%%%%% rows (em un cr le re lu ru sy at as co tr) %%%%%%%%
%%%%%%%% cols (s1 ir rf qr an sc it s2 qt de ls rs) %%%%%%%%

%%%%%%%% start of table + - %%%%%%%%
%%%%%%%% rows (s1 ir rf qr an sc it s2 qt de ls rs) %%%%%%%%
%%%%%%%% cols (s1 ir rf qr an sc it s2 qt de ls rs) %%%%%%%%

	}
\\\hline
\hline
SemiOrd1 & s1%
	{%
	} & \t{ % +s1 -s1
		\x
	} & \t{ % +s1 -ir
	} & \t{ % +s1 -rf
	} & \t{ % +s1 -qr
	} & \t{ % +s1 -an
	} & \t{ % +s1 -sc
	} & \t{ % +s1 -it
		\4{241},%	% +lu +s1 -it +ls
		\4{249}\\	% +ru +s1 -it +rs
		\5{253},%	% +le +tr +s1 -it +ls
		\5{256}\\	% +re +tr +s1 -it +rs
	} & \t{ % +s1 -s2
	} & \t{ % +s1 -qt
	} & \t{ % +s1 -de
		\3{187}\\	% +sy +s1 -de
	} & \t{ % +s1 -ls
	} & \t{ % +s1 -rs
	}
\\\hline
Irrefl & ir%
	\t{\\ \\%
	} & \t{ % +ir -s1
	} & \t{ % +ir -ir
		\x
	} & \t{ % +ir -rf
	} & \t{ % +ir -qr
	} & \t{ % +ir -an
	} & \t{ % +ir -sc
	} & \t{ % +ir -it
	} & \t{ % +ir -s2
	} & \t{ % +ir -qt
	} & \t{ % +ir -de
	} & \t{ % +ir -ls
	} & \t{ % +ir -rs
	}
\\\hline
Refl & rf%
	\t{\\ \\ \\%
	} & \t{ % +rf -s1
	} & \t{ % +rf -ir
	} & \t{ % +rf -rf
		\x
	} & \t{ % +rf -qr
		\2{051}\\	% +rf -qr
	} & \t{ % +rf -an
	} & \t{ % +rf -sc
	} & \t{ % +rf -it
	} & \t{ % +rf -s2
	} & \t{ % +rf -qt
	} & \t{ % +rf -de
		%\2{085}\\	% +rf -de
	} & \t{ % +rf -ls
		\2{090}\\	% +rf -ls
	} & \t{ % +rf -rs
		\2{094}\\	% +rf -rs
	}
\\\hline
QuasiRefl & qr%
	{%
	} & \t{ % +qr -s1
		\3{168},%	% -s1 +qr +it
		\3{176}\\	% -s1 +qr +s2
	} & \t{ % +qr -ir
	} & \t{ % +qr -rf
		\3{197}\\	% -rf +qr +ls
		\3{209}\\	% -rf +qr +rs
	} & \t{ % +qr -qr
		\x
	} & \t{ % +qr -an
	} & \t{ % +qr -sc
	} & \t{ % +qr -it
		\3{180}\\	% +qr -it +s2
	} & \t{ % +qr -s2
	} & \t{ % +qr -qt
	} & \t{ % +qr -de
		\2{086}\\	% +qr -de
	} & \t{ % +qr -ls
	} & \t{ % +qr -rs
	}
\\\hline
AntiSym & an%
	{%
	} & \t{ % +an -s1
	} & \t{ % +an -ir
	} & \t{ % +an -rf
	} & \t{ % +an -qr
	} & \t{ % +an -an
		\x
	} & \t{ % +an -sc
	} & \t{ % +an -it
	} & \t{ % +an -s2
	} & \t{ % +an -qt
	} & \t{ % +an -de
	} & \t{ % +an -ls
	} & \t{ % +an -rs
	}
\\\hline
SemiConnex & sc%
	{%
	} & \t{ % +sc -s1
		%\3{156},%	% +le -s1 +sc
		%\3{157}\\	% +re -s1 +sc
		%\3{158}\\	% +tr -s1 +sc
	} & \t{ % +sc -ir
	} & \t{ % +sc -rf
	} & \t{ % +sc -qr
	} & \t{ % +sc -an
	} & \t{ % +sc -sc
		\x
	} & \t{ % +sc -it
		\2{066}\\	% +sc -it
	} & \t{ % +sc -s2
		%\2{070}\\	% +sc -s2
	} & \t{ % +sc -qt
	} & \t{ % +sc -de
		\3{188}\\	% +sy +sc -de
	} & \t{ % +sc -ls
		\3{198}\\	% +le +sc -ls
		\3{199}\\	% +sy +sc -ls
	} & \t{ % +sc -rs
		\3{210}\\	% +re +sc -rs
		\3{211}\\	% +sy +sc -rs
	}
\\\hline
IncTrans & it%
	\t{\\ \\ \\%
	} & \t{ % +it -s1
		\3{164},%	% +le -s1 +it
		\3{165}\\	% +re -s1 +it
		\3{166},%	% +tr -s1 +it
		\3{168}\\	% -s1 +qr +it
	} & \t{ % +it -ir
	} & \t{ % +it -rf
	} & \t{ % +it -qr
	} & \t{ % +it -an
	} & \t{ % +it -sc
	} & \t{ % +it -it
		\x
	} & \t{ % +it -s2
		\2{071}\\	% +it -s2
	} & \t{ % +it -qt
		%\4{228}\\	% +lu -as +it -qt
		%\4{229}\\	% +ru -as +it -qt
	} & \t{ % +it -de
		\4{232}\\	% +sy -at +it -de
	} & \t{ % +it -ls
		\4{237}\\	% -em +le +it -ls
		\4{238}\\	% -em +sy +it -ls
		%\4{239}\\	% +lu -as +it -ls
	} & \t{ % +it -rs
		\4{246}\\	% -em +re +it -rs
		%\4{247}\\	% +ru -as +it -rs
	}
\\\hline
SemiOrd2 & s2%
	\t{\\ \\ \\%
	} & \t{ % +s2 -s1
		\3{173},%	% +le -s1 +s2
		\3{174}\\	% +re -s1 +s2
		\3{176},%	% -s1 +qr +s2
		\4{233}\\	% +tr -s1 +s2 +de
	} & \t{ % +s2 -ir
	} & \t{ % +s2 -rf
	} & \t{ % +s2 -qr
	} & \t{ % +s2 -an
	} & \t{ % +s2 -sc
	} & \t{ % +s2 -it
		\3{177},%	% +le -it +s2		(+)
		\3{178}\\	% +re -it +s2		(+)
		\3{179},%	% +sy -it +s2
		\3{180}\\	% +qr -it +s2
		\4{234},%	% +tr -it +s2 +de
		\6{265}\\	% +at -it +s2 +qt +ls +rs
	} & \t{ % +s2 -s2
		\x
	} & \t{ % +s2 -qt
		\4{230}\\	% +lu -as +s2 -qt
		\4{231}\\	% +ru -as +s2 -qt
	} & \t{ % +s2 -de
	} & \t{ % +s2 -ls
		\4{242}\\	% +lu -as +s2 -ls
		\5{254}\\	% +lu +ru -tr +s2 -ls
	} & \t{ % +s2 -rs
		\4{250}\\	% +ru -as +s2 -rs
	}
\\\hline
QuasiTrans & qt%
	{%
	} & \t{ % +qt -s1
	} & \t{ % +qt -ir
	} & \t{ % +qt -rf
	} & \t{ % +qt -qr
	} & \t{ % +qt -an
	} & \t{ % +qt -sc
	} & \t{ % +qt -it
		\6{265}\\	% +at -it +s2 +qt +ls +rs
	} & \t{ % +qt -s2
	} & \t{ % +qt -qt
		\x
	} & \t{ % +qt -de
	} & \t{ % +qt -ls
	} & \t{ % +qt -rs
	}
\\\hline
Dense & de%
	{%
	} & \t{ % +de -s1
		\4{233}\\	% +tr -s1 +s2 +de
	} & \t{ % +de -ir
	} & \t{ % +de -rf
	} & \t{ % +de -qr
	} & \t{ % +de -an
	} & \t{ % +de -sc
	} & \t{ % +de -it
		\4{234}\\	% +tr -it +s2 +de
	} & \t{ % +de -s2
	} & \t{ % +de -qt
	} & \t{ % +de -de
		\x
	} & \t{ % +de -ls
	} & \t{ % +de -rs
	}
\\\hline
LfSerial & ls%
	{%
	} & \t{ % +ls -s1
	} & \t{ % +ls -ir
	} & \t{ % +ls -rf
		\3{195}\\	% +cr -rf +ls
		\3{196}\\	% +re -rf +ls		(+)
		\3{197}\\	% -rf +qr +ls
	} & \t{ % +ls -qr
	} & \t{ % +ls -an
	} & \t{ % +ls -sc
	} & \t{ % +ls -it
		\4{241},%	% +lu +s1 -it +ls
		\5{253}\\	% +le +tr +s1 -it +ls
		\6{265}\\	% +at -it +s2 +qt +ls +rs
	} & \t{ % +ls -s2
	} & \t{ % +ls -qt
	} & \t{ % +ls -de
	} & \t{ % +ls -ls
		\x
	} & \t{ % +ls -rs
		\3{214}\\	% +ru +ls -rs
		\3{215}\\	% +sy +ls -rs
	}
\\\hline
RgSerial & rs%
	{%
	} & \t{ % +rs -s1
	} & \t{ % +rs -ir
	} & \t{ % +rs -rf
		\3{207}\\	% +cr -rf +rs
		\3{208}\\	% +le -rf +rs		(+)
		\3{209}\\	% -rf +qr +rs
	} & \t{ % +rs -qr
	} & \t{ % +rs -an
	} & \t{ % +rs -sc
	} & \t{ % +rs -it
		\4{249},%	% +ru +s1 -it +rs
		\5{256}\\	% +re +tr +s1 -it +rs
		\6{265}\\	% +at -it +s2 +qt +ls +rs
	} & \t{ % +rs -s2
	} & \t{ % +rs -qt
	} & \t{ % +rs -de
	} & \t{ % +rs -ls
		\3{213}\\	% +lu -ls +rs
		\3{216}\\	% +sy -ls +rs
	} & \t{ % +rs -rs
		\x

	}
\\\hline
\end{tabular}
\caption{Law index rg ($A \Ra B$)}
%\rule{\textwidth}{0.1mm}
\LABEL{Law index + - rg}
\end{center}
\end{figure}

%%%%%%%% end of table + - %%%%%%%%
%%%%%%%% rows (s1 ir rf qr an sc it s2 qt de ls rs) %%%%%%%%
%%%%%%%% cols (s1 ir rf qr an sc it s2 qt de ls rs) %%%%%%%%

%%%%%%%% start of table - - %%%%%%%%
%%%%%%%% rows (em un cr le re lu ru sy at as co tr) %%%%%%%%
%%%%%%%% cols (em un cr le re lu ru sy at as co tr) %%%%%%%%

\begin{figure}
\begin{center}
\scriptsize
\begin{tabular}{|l@{$\;$}c|*{12}{|@{$\;$}c@{$\;$}}}
\hline
\multicolumn{2}{|r||@{$\;$}}{$_A \setminus ^B$}
	& em & un & cr & le & re & lu & ru & sy & at & as & co & tr \\
\hline
\hline
Empty & em%
	{%
	} & \t{ % -em -em
		\x		% (inverting)
	} & \t{ % +em +un
		\2{001}\\	% +em +un
	} & \t{ % +em +cr
	} & \t{ % +em +le
	} & \t{ % +em +re
	} & \t{ % +em +lu
	} & \t{ % +em +ru
	} & \t{ % +em +sy
	} & \t{ % +em +at
	} & \t{ % +em +as
	} & \t{ % +em +co
		\2{023}\\	% +em +co
	} & \t{ % +em +tr
	}
\\\hline
Univ & un%
	{%
	} & \t{ % -un -em
	} & \t{ % -un -un
		\x		% (inverting)
	} & \t{ % +un +cr
		\2{003}\\	% +un +cr
	} & \t{ % +un +le
	} & \t{ % +un +re
	} & \t{ % +un +lu
		\2{011}\\	% +un +lu
	} & \t{ % +un +ru
		\2{014}\\	% +un +ru
	} & \t{ % +un +sy
	} & \t{ % +un +at
		\2{020}\\	% +un +at
	} & \t{ % +un +as
		\2{022}\\	% +un +as
	} & \t{ % +un +co
	} & \t{ % +un +tr
	}
\\\hline
CoRefl & cr%
	\t{\\ \\ \\%
	} & \t{ % -cr -em
	} & \t{ % -cr -un
	} & \t{ % -cr -cr
		\x		% (inverting)
	} & \t{ % +cr +le
	} & \t{ % +cr +re
	} & \t{ % +cr +lu
	} & \t{ % +cr +ru
	} & \t{ % +cr +sy
	} & \t{ % +cr +at
		%\3{102}\\	% -em +cr +at
	} & \t{ % +cr +as
		%\3{105}\\	% -em +cr +as
	} & \t{ % +cr +co
		%\2{025}\\	% +cr +co
	} & \t{ % +cr +tr
	}
\\\hline
LfEucl & le%
	{%
	} & \t{ % -le -em
	} & \t{ % -le -un
	} & \t{ % -le -cr
	} & \t{ % -le -le
		\x		% (inverting)
	} & \t{ % +le +re
		\3{097}\\	% +le +re -sy
		\3{114}\\	% +le +re -tr
		\3{137}\\	% +le +re -qr
	} & \t{ % +le +lu
		\3{115}\\	% +le +lu -tr
		\3{145}\\	% +le +lu -an
	} & \t{ % +le +ru
		\3{096}\\	% -cr +le +ru		(+)
	} & \t{ % +le +sy
		\3{098}\\	% +le -re +sy
	} & \t{ % +le +at
		%\3{103}\\	% -em +le +at
	} & \t{ % +le +as
		%\3{106}\\	% -em +le +as
	} & \t{ % +le +co
		%\3{111}\\	% -un +le +co
	} & \t{ % +le +tr
		\5{253}\\	% +le +tr +s1 -it +ls
	}
\\\hline
RgEucl & re%
	\t{\\ \\%
	} & \t{ % -re -em
	} & \t{ % -re -un
	} & \t{ % -re -cr
	} & \t{ % -re -le
	} & \t{ % -re -re
		\x		% (inverting)
	} & \t{ % +re +lu
		\3{095}\\	% -cr +re +lu		(+)
	} & \t{ % +re +ru
		\3{116}\\	% +re +ru -tr
		\3{147}\\	% +re +ru -an
	} & \t{ % +re +sy
		\3{099}\\	% -le +re +sy
	} & \t{ % +re +at
		%\3{104}\\	% -em +re +at
	} & \t{ % +re +as
		%\3{107}\\	% -em +re +as
	} & \t{ % +re +co
		%\3{112}\\	% -un +re +co
	} & \t{ % +re +tr
		\5{256}\\	% +re +tr +s1 -it +rs
	}
\\\hline
LfUnique & lu%
	\t{\\ \\ \\ \\%
	} & \t{ % -lu -em
	} & \t{ % -lu -un
	} & \t{ % -lu -cr
	} & \t{ % -lu -le
	} & \t{ % -lu -re
	} & \t{ % -lu -lu
		\x		% (inverting)
	} & \t{ % +lu +ru
		\4{220}\\	% -em +lu +ru +it
		\4{226}\\	% +lu +ru -as +s2
		\5{254}\\	% +lu +ru -tr +s2 -ls
	} & \t{ % +lu +sy
		\3{100}\\	% +lu -ru +sy
		\4{235}\\	% +lu +sy +at +ls
	} & \t{ % +lu +at
		\4{235}\\	% +lu +sy +at +ls
	} & \t{ % +lu +as
		%\3{109}\\	% +lu -at +as
		\4{240}\\	% +lu +as +it +ls
	} & \t{ % +lu +co
		%\2{026}\\	% +lu +co
	} & \t{ % +lu +tr
		\3{151},%	% +lu +tr -an
		\4{236}\\	% -le +lu +tr +ls
	}
\\\hline
RgUnique & ru%
	\t{\\ \\ \\ \\%
	} & \t{ % -ru -em
	} & \t{ % -ru -un
	} & \t{ % -ru -cr
	} & \t{ % -ru -le
	} & \t{ % -ru -re
	} & \t{ % -ru -lu
	} & \t{ % -ru -ru
		\x		% (inverting)
	} & \t{ % +ru +sy
		\3{101}\\	% -lu +ru +sy
	} & \t{ % +ru +at
	} & \t{ % +ru +as
		%\3{110}\\	% +ru -at +as
		\4{248}\\	% +ru +as +it +rs
	} & \t{ % +ru +co
		%\2{027}\\	% +ru +co
	} & \t{ % +ru +tr
		\3{152},%	% +ru +tr -an
		\4{245}\\	% -re +ru +tr +rs
	}
\\\hline
Sym & sy%
	\t{\\ \\ \\%
	} & \t{ % -sy -em
	} & \t{ % -sy -un
	} & \t{ % -sy -cr
	} & \t{ % -sy -le
	} & \t{ % -sy -re
	} & \t{ % -sy -lu
	} & \t{ % -sy -ru
	} & \t{ % -sy -sy
		\x		% (inverting)
	} & \t{ % +sy +at
		\4{235}\\	% +lu +sy +at +ls
	} & \t{ % +sy +as
		\3{108}\\	% -em +sy +as
	} & \t{ % +sy +co
		\3{113}\\	% -un +sy +co
	} & \t{ % +sy +tr
		\3{117}\\	% -le +sy +tr
	}
\\\hline
AntiTrans & at%
	{%
	} & \t{ % -at -em
	} & \t{ % -at -un
	} & \t{ % -at -cr
	} & \t{ % -at -le
		\4{217}\\	% -le +lu -at +s1
		%\4{221}\\	% -le +lu -at +it
		\4{224}\\	% -le +lu -at +s2
	} & \t{ % -at -re
		\4{218}\\	% -re +ru -at +s1
		%\4{222}\\	% -re +ru -at +it
		\4{225}\\	% -re +ru -at +s2
	} & \t{ % -at -lu
	} & \t{ % -at -ru
	} & \t{ % -at -sy
	} & \t{ % -at -at
		\x		% (inverting)
	} & \t{ % +at +as
	} & \t{ % +at +co
		%\2{028}\\	% +at +co
	} & \t{ % +at +tr
		%\3{118},%	% +at -as +tr
		\3{172},%	% +at +tr -s2
		\3{191}\\	% +at +tr +ls
		\3{203}\\	% +at +tr +rs
	}
\\\hline
ASym & as%
	{%
	} & \t{ % -as -em
	} & \t{ % -as -un
	} & \t{ % -as -cr
	} & \t{ % -as -le
	} & \t{ % -as -re
	} & \t{ % -as -lu
	} & \t{ % -as -ru
	} & \t{ % -as -sy
	} & \t{ % -as -at
	} & \t{ % -as -as
		\x		% (inverting)
	} & \t{ % +as +co
		%\2{029}\\	% +as +co
	} & \t{ % +as +tr
		\3{192},%	% +as +tr +ls
		\3{204}\\	% +as +tr +rs
	}
\\\hline
Connex & co%
	{%
	} & \t{ % -co -em
		\4{223}\\	% -em -co +qr +it
	} & \t{ % -co -un
	} & \t{ % -co -cr
	} & \t{ % -co -le
	} & \t{ % -co -re
	} & \t{ % -co -lu
	} & \t{ % -co -ru
	} & \t{ % -co -sy
	} & \t{ % -co -at
	} & \t{ % -co -as
	} & \t{ % -co -co
		\x		% (inverting)
	} & \t{ % +co +tr
	}
\\\hline
Trans & tr%
	\t{\\ \\%
	} & \t{ % -tr -em
	} & \t{ % -tr -un
	} & \t{ % -tr -cr
	} & \t{ % -tr -le
	} & \t{ % -tr -re
	} & \t{ % -tr -lu
	} & \t{ % -tr -ru
	} & \t{ % -tr -sy
	} & \t{ % -tr -at
	} & \t{ % -tr -as
	} & \t{ % -tr -co
	} & \t{ % -tr -tr
		\x		% (inverting)

%%%%%%%% end of table - - %%%%%%%%
%%%%%%%% rows (em un cr le re lu ru sy at as co tr) %%%%%%%%
%%%%%%%% cols (em un cr le re lu ru sy at as co tr) %%%%%%%%

%%%%%%%% start of table - - %%%%%%%%
%%%%%%%% rows (s1 ir rf qr an sc it s2 qt de ls rs) %%%%%%%%
%%%%%%%% cols (em un cr le re lu ru sy at as co tr) %%%%%%%%

	}
\\\hline
\hline
SemiOrd1 & s1%
	\t{\\ \\%
	} & \t{ % -s1 -em
	} & \t{ % -s1 -un
	} & \t{ % -s1 -cr
	} & \t{ % -s1 -le
	} & \t{ % -s1 -re
	} & \t{ % -s1 -lu
	} & \t{ % -s1 -ru
	} & \t{ % -s1 -sy
	} & \t{ % -s1 -at
	} & \t{ % -s1 -as
	} & \t{ % -s1 -co
	} & \t{ % -s1 -tr
	}
\\\hline
Irrefl & ir%
	{%
	} & \t{ % -ir -em
	} & \t{ % -ir -un
	} & \t{ % -ir -cr
	} & \t{ % -ir -le
	} & \t{ % -ir -re
	} & \t{ % -ir -lu
	} & \t{ % -ir -ru
	} & \t{ % -ir -sy
	} & \t{ % -ir -at
	} & \t{ % -ir -as
	} & \t{ % -ir -co
	} & \t{ % -ir -tr
	}
\\\hline
Refl & rf%
	{%
	} & \t{ % -rf -em
	} & \t{ % -rf -un
	} & \t{ % -rf -cr
	} & \t{ % -rf -le
	} & \t{ % -rf -re
	} & \t{ % -rf -lu
	} & \t{ % -rf -ru
	} & \t{ % -rf -sy
	} & \t{ % -rf -at
	} & \t{ % -rf -as
	} & \t{ % -rf -co
	} & \t{ % -rf -tr
	}
\\\hline
QuasiRefl & qr%
	{%
	} & \t{ % -qr -em
	} & \t{ % -qr -un
	} & \t{ % -qr -cr
	} & \t{ % -qr -le
	} & \t{ % -qr -re
	} & \t{ % -qr -lu
	} & \t{ % -qr -ru
	} & \t{ % -qr -sy
	} & \t{ % -qr -at
	} & \t{ % -qr -as
	} & \t{ % -qr -co
	} & \t{ % -qr -tr
	}
\\\hline
AntiSym & an%
	{%
	} & \t{ % -an -em
	} & \t{ % -an -un
	} & \t{ % -an -cr
	} & \t{ % -an -le
	} & \t{ % -an -re
	} & \t{ % -an -lu
	} & \t{ % -an -ru
	} & \t{ % -an -sy
	} & \t{ % -an -at
	} & \t{ % -an -as
	} & \t{ % -an -co
	} & \t{ % -an -tr
	}
\\\hline
SemiConnex & sc%
	{%
	} & \t{ % -sc -em
	} & \t{ % -sc -un
	} & \t{ % -sc -cr
	} & \t{ % -sc -le
	} & \t{ % -sc -re
	} & \t{ % -sc -lu
	} & \t{ % -sc -ru
	} & \t{ % -sc -sy
	} & \t{ % -sc -at
	} & \t{ % -sc -as
	} & \t{ % -sc -co
	} & \t{ % -sc -tr
	}
\\\hline
IncTrans & it%
	{%
	} & \t{ % -it -em
	} & \t{ % -it -un
	} & \t{ % -it -cr
	} & \t{ % -it -le
	} & \t{ % -it -re
	} & \t{ % -it -lu
	} & \t{ % -it -ru
	} & \t{ % -it -sy
	} & \t{ % -it -at
	} & \t{ % -it -as
	} & \t{ % -it -co
	} & \t{ % -it -tr
	}
\\\hline
SemiOrd2 & s2%
	\t{\\ \\%
	} & \t{ % -s2 -em
	} & \t{ % -s2 -un
	} & \t{ % -s2 -cr
	} & \t{ % -s2 -le
	} & \t{ % -s2 -re
	} & \t{ % -s2 -lu
	} & \t{ % -s2 -ru
	} & \t{ % -s2 -sy
	} & \t{ % -s2 -at
	} & \t{ % -s2 -as
	} & \t{ % -s2 -co
	} & \t{ % -s2 -tr
	}
\\\hline
QuasiTrans & qt%
	{%
	} & \t{ % -qt -em
	} & \t{ % -qt -un
	} & \t{ % -qt -cr
	} & \t{ % -qt -le
	} & \t{ % -qt -re
	} & \t{ % -qt -lu
	} & \t{ % -qt -ru
	} & \t{ % -qt -sy
	} & \t{ % -qt -at
	} & \t{ % -qt -as
		%\4{228}\\	% +lu -as +it -qt
		%\4{229}\\	% +ru -as +it -qt
		\4{230}\\	% +lu -as +s2 -qt
		\4{231}\\	% +ru -as +s2 -qt
	} & \t{ % -qt -co
	} & \t{ % -qt -tr
	}
\\\hline
Dense & de%
	{%
	} & \t{ % -de -em
	} & \t{ % -de -un
	} & \t{ % -de -cr
	} & \t{ % -de -le
	} & \t{ % -de -re
	} & \t{ % -de -lu
	} & \t{ % -de -ru
	} & \t{ % -de -sy
	} & \t{ % -de -at
		\4{232}\\	% +sy -at +it -de
	} & \t{ % -de -as
	} & \t{ % -de -co
	} & \t{ % -de -tr
	}
\\\hline
LfSerial & ls%
	{%
	} & \t{ % -ls -em
		\4{237}\\	% -em +le +it -ls
		\4{238}\\	% -em +sy +it -ls
	} & \t{ % -ls -un
	} & \t{ % -ls -cr
	} & \t{ % -ls -le
	} & \t{ % -ls -re
	} & \t{ % -ls -lu
	} & \t{ % -ls -ru
	} & \t{ % -ls -sy
	} & \t{ % -ls -at
	} & \t{ % -ls -as
		%\4{239}\\	% +lu -as +it -ls
		\4{242}\\	% +lu -as +s2 -ls
	} & \t{ % -ls -co
	} & \t{ % -ls -tr
		\5{254}\\	% +lu +ru -tr +s2 -ls
	}
\\\hline
RgSerial & rs%
	{%
	} & \t{ % -rs -em
		\4{246}\\	% -em +re +it -rs
	} & \t{ % -rs -un
	} & \t{ % -rs -cr
	} & \t{ % -rs -le
	} & \t{ % -rs -re
	} & \t{ % -rs -lu
	} & \t{ % -rs -ru
	} & \t{ % -rs -sy
	} & \t{ % -rs -at
	} & \t{ % -rs -as
		%\4{247}\\	% +ru -as +it -rs
		\4{250}\\	% +ru -as +s2 -rs
	} & \t{ % -rs -co
	} & \t{ % -rs -tr

	}
\\\hline
\end{tabular}
\caption{Law index lf (bot lf: $A \lor B$ required, top rg: $A \land B$ impossible)}
%\rule{\textwidth}{0.1mm}
\LABEL{Law index - - / + + lf}
\end{center}
\end{figure}

%%%%%%%% end of table - - %%%%%%%%
%%%%%%%% rows (s1 ir rf qr an sc it s2 qt de ls rs) %%%%%%%%
%%%%%%%% cols (em un cr le re lu ru sy at as co tr) %%%%%%%%

%%%%%%%% start of table + + %%%%%%%%
%%%%%%%% rows (em un cr le re lu ru sy at as co tr) %%%%%%%%
%%%%%%%% cols (s1 ir rf qr an sc it s2 qt de ls rs) %%%%%%%%

\begin{figure}
\begin{center}
\scriptsize
\begin{tabular}{lc|*{12}{|@{$\;$}c@{$\;$}}|}
\hline
\multicolumn{2}{r||@{$\;$}}{$_A \setminus ^B$}
	& s1 & ir & rf & qr & an & sc & it & s2 & qt & de & ls & rs \\
\hline
\hline
Empty & em%
	{%
	} & \t{ % +em +s1
	} & \t{ % +em +ir
	} & \t{ % +em +rf
		\2{041}\\	% +em +rf
	} & \t{ % +em +qr
	} & \t{ % +em +an
	} & \t{ % +em +sc
		\2{056}\\	% +em +sc
	} & \t{ % +em +it
	} & \t{ % +em +s2
	} & \t{ % +em +qt
	} & \t{ % +em +de
	} & \t{ % +em +ls
		\2{087}\\	% +em +ls
	} & \t{ % +em +rs
		\2{091}\\	% +em +rs
	}
\\\hline
Univ & un%
	{%
	} & \t{ % +un +s1
	} & \t{ % +un +ir
		\2{037}\\	% +un +ir
	} & \t{ % +un +rf
	} & \t{ % +un +qr
	} & \t{ % +un +an
		\2{053}\\	% +un +an
	} & \t{ % +un +sc
	} & \t{ % +un +it
	} & \t{ % +un +s2
	} & \t{ % +un +qt
	} & \t{ % +un +de
	} & \t{ % +un +ls
	} & \t{ % +un +rs
	}
\\\hline
CoRefl & cr%
	{%
	} & \t{ % +cr +s1
		\3{135}\\	% +cr +s1 +rf
		\3{193}\\	% +cr +s1 +ls
		\3{205}\\	% +cr +s1 +rs
	} & \t{ % +cr +ir
		\3{124}\\	% -em +cr +ir
	} & \t{ % +cr +rf
		\3{135}\\	% +cr +s1 +rf
	} & \t{ % +cr +qr
	} & \t{ % +cr +an
	} & \t{ % +cr +sc
		\2{058}\\	% +cr +sc
	} & \t{ % +cr +it
		\3{161}\\	% -em +cr +it
	} & \t{ % +cr +s2
		\3{169}\\	% -em +cr +s2
	} & \t{ % +cr +qt
	} & \t{ % +cr +de
	} & \t{ % +cr +ls
		\3{193},%	% +cr +s1 +ls
		\3{195}\\	% +cr -rf +ls
	} & \t{ % +cr +rs
		\3{205},%	% +cr +s1 +rs
		\3{207}\\	% +cr -rf +rs
	}
\\\hline
LfEucl & le%
	\t{\\ \\ \\%
	} & \t{ % +le +s1
		\5{253}\\	% +le +tr +s1 -it +ls
	} & \t{ % +le +ir
		\3{125}\\	% -em +le +ir
	} & \t{ % +le +rf
		\3{131}\\	% +le -re +rf
	} & \t{ % +le +qr
		\3{138}\\	% +le -re +qr
	} & \t{ % +le +an
		\3{146}\\	% +le -lu +an
	} & \t{ % +le +sc
		\3{154}\\	% +le -tr +sc
		%\3{156}\\	% +le -s1 +sc
		\3{198}\\	% +le +sc -ls
	} & \t{ % +le +it
		\3{162},%	% +le -tr +it
		\3{164}\\	% +le -s1 +it
		\4{237}\\	% -em +le +it -ls
	} & \t{ % +le +s2
		\3{170},%	% +le -tr +s2
		\3{173}\\	% +le -s1 +s2
		\3{177}\\	% +le -it +s2		(+)
	} & \t{ % +le +qt
	} & \t{ % +le +de
	} & \t{ % +le +ls
		\5{253}\\	% +le +tr +s1 -it +ls
	} & \t{ % +le +rs
		\3{201}%	% +le -re +rs
		\3{208}\\	% +le -rf +rs		(+)
	}
\\\hline
RgEucl & re%
	\t{\\ \\%
	} & \t{ % +re +s1
		\5{256}\\	% +re +tr +s1 -it +rs
	} & \t{ % +re +ir
		\3{126}\\	% -em +re +ir
	} & \t{ % +re +rf
		\3{132}\\	% -le +re +rf
	} & \t{ % +re +qr
		\3{139}\\	% -le +re +qr
	} & \t{ % +re +an
		\3{148}\\	% +re -ru +an
	} & \t{ % +re +sc
		\3{155}\\	% +re -tr +sc
		%\3{157}\\	% +re -s1 +sc
		\3{210}\\	% +re +sc -rs
	} & \t{ % +re +it
		\3{163},%	% +re -tr +it
		\3{165}\\	% +re -s1 +it
		\4{246}\\	% -em +re +it -rs
	} & \t{ % +re +s2
		\3{171},%	% +re -tr +s2
		\3{174}\\	% +re -s1 +s2
		\3{178}\\	% +re -it +s2		(+)
	} & \t{ % +re +qt
	} & \t{ % +re +de
	} & \t{ % +re +ls
		\3{189},%	% -le +re +ls
		\3{196}\\	% +re -rf +ls		(+)
	} & \t{ % +re +rs
		\5{256}\\	% +re +tr +s1 -it +rs
	}
\\\hline
LfUnique & lu%
	{%
	} & \t{ % +lu +s1
		\3{120}\\	% +lu -tr +s1
		\3{206}\\	% +lu +s1 +rs
		\4{217}\\	% -le +lu -at +s1
		\4{241}\\	% +lu +s1 -it +ls
	} & \t{ % +lu +ir
		\3{127}\\	% +lu -at +ir
	} & \t{ % +lu +rf
		%\3{133}\\	% -cr +lu +rf
	} & \t{ % +lu +qr
		\3{140}\\	% -cr +lu +qr
	} & \t{ % +lu +an
	} & \t{ % +lu +sc
		\2{059}\\	% +lu +sc
	} & \t{ % +lu +it
		\3{212},%	% +lu +it +rs
		\4{220}\\	% -em +lu +ru +it
		%\4{221},%	% -le +lu -at +it
		%\4{228}\\	% +lu -as +it -qt
		%\4{239},%	% +lu -as +it -ls
		\4{240}\\	% +lu +as +it +ls
	} & \t{ % +lu +s2
		\4{224},%	% -le +lu -at +s2
		\4{226}\\	% +lu +ru -as +s2
		\4{230},%	% +lu -as +s2 -qt
		\4{242}\\	% +lu -as +s2 -ls
		\5{254}\\	% +lu +ru -tr +s2 -ls
	} & \t{ % +lu +qt
	} & \t{ % +lu +de
		\3{183}\\	% -le +lu +de
	} & \t{ % +lu +ls
		\4{235},%	% +lu +sy +at +ls
		\4{236}\\	% -le +lu +tr +ls
		\4{240},%	% +lu +as +it +ls
		\4{241}\\	% +lu +s1 -it +ls
	} & \t{ % +lu +rs
		\3{202},%	% +lu -ru +rs
		\3{206}\\	% +lu +s1 +rs
		\3{212},%	% +lu +it +rs
		\3{213}\\	% +lu -ls +rs
	}
\\\hline
RgUnique & ru%
	{%
	} & \t{ % +ru +s1
		\3{121}\\	% +ru -tr +s1
		\3{194}\\	% +ru +s1 +ls
		\4{218}\\	% -re +ru -at +s1
		\4{249}\\	% +ru +s1 -it +rs
	} & \t{ % +ru +ir
		\3{128}\\	% +ru -at +ir
	} & \t{ % +ru +rf
		%\3{134}\\	% -cr +ru +rf
	} & \t{ % +ru +qr
		\3{141}\\	% -cr +ru +qr
	} & \t{ % +ru +an
	} & \t{ % +ru +sc
		\2{060}\\	% +ru +sc
	} & \t{ % +ru +it
		\3{200},%	% +ru +it +ls
		\4{220}\\	% -em +lu +ru +it
		%\4{222},%	% -re +ru -at +it
		%\4{229}\\	% +ru -as +it -qt
		%\4{247},%	% +ru -as +it -rs
		\4{248}\\	% +ru +as +it +rs
	} & \t{ % +ru +s2
		\4{225},%	% -re +ru -at +s2
		\4{226}\\	% +lu +ru -as +s2
		\4{231},%	% +ru -as +s2 -qt
		\4{250}\\	% +ru -as +s2 -rs
		\5{254}\\	% +lu +ru -tr +s2 -ls
	} & \t{ % +ru +qt
		\4{244}\\	% +ru -sy +qt +ls
	} & \t{ % +ru +de
		\3{184}\\	% -re +ru +de
	} & \t{ % +ru +ls
		\3{190},%	% -lu +ru +ls
		\3{194}\\	% +ru +s1 +ls
		\3{200},%	% +ru +it +ls
		\3{214}\\	% +ru +ls -rs
		\4{244}\\	% +ru -sy +qt +ls
	} & \t{ % +ru +rs
		\4{245},%	% -re +ru +tr +rs
		\4{248}\\	% +ru +as +it +rs
		\4{249}\\	% +ru +s1 -it +rs
	}
\\\hline
Sym & sy%
	{%
	} & \t{ % +sy +s1
		\3{187}\\	% +sy +s1 -de
		\4{219}\\	% -le +sy +s1 +qr
	} & \t{ % +sy +ir
	} & \t{ % +sy +rf
	} & \t{ % +sy +qr
		\4{219}\\	% -le +sy +s1 +qr
	} & \t{ % +sy +an
		\3{149}\\	% -cr +sy +an
	} & \t{ % +sy +sc
		\3{188}\\	% +sy +sc -de
		\3{199}\\	% +sy +sc -ls
		\3{211}\\	% +sy +sc -rs
	} & \t{ % +sy +it
		\4{232},%	% +sy -at +it -de
		\4{238}\\	% -em +sy +it -ls
	} & \t{ % +sy +s2
		\3{179}\\	% +sy -it +s2
	} & \t{ % +sy +qt
	} & \t{ % +sy +de
	} & \t{ % +sy +ls
		\3{215},%	% +sy +ls -rs
		\4{235}\\	% +lu +sy +at +ls
	} & \t{ % +sy +rs
		\3{216}\\	% +sy -ls +rs
	}
\\\hline
AntiTrans & at%
	\t{\\ \\%
	} & \t{ % +at +s1
		%\3{119}\\	% +at -as +s1
		%\3{122}\\	% +at -tr +s1
	} & \t{ % +at +ir
	} & \t{ % +at +rf
		%\2{043}\\	% +at +rf
	} & \t{ % +at +qr
		%\3{142}\\	% -em +at +qr
	} & \t{ % +at +an
		%\3{150}\\	% +at -as +an
	} & \t{ % +at +sc
		\2{061}\\	% +at +sc
	} & \t{ % +at +it
	} & \t{ % +at +s2
		\6{265}\\	% +at -it +s2 +qt +ls +rs
	} & \t{ % +at +qt
		\6{265}\\	% +at -it +s2 +qt +ls +rs
	} & \t{ % +at +de
		\3{185}\\	% -em +at +de
	} & \t{ % +at +ls
		\3{191},%	% +at +tr +ls
		\4{235}\\	% +lu +sy +at +ls
		\6{265}\\	% +at -it +s2 +qt +ls +rs
	} & \t{ % +at +rs
		\3{203},%	% +at +tr +rs
		\6{265}\\	% +at -it +s2 +qt +ls +rs
	}
\\\hline
ASym & as%
	{%
	} & \t{ % +as +s1
		%\3{123}\\	% +as -tr +s1
	} & \t{ % +as +ir
	} & \t{ % +as +rf
		%\2{044}\\	% +as +rf
	} & \t{ % +as +qr
		%\3{143}\\	% -em +as +qr
	} & \t{ % +as +an
	} & \t{ % +as +sc
	} & \t{ % +as +it
		\4{240},%	% +lu +as +it +ls
		\4{248}\\	% +ru +as +it +rs
	} & \t{ % +as +s2
	} & \t{ % +as +qt
		%\3{181}\\	% +as -tr +qt
	} & \t{ % +as +de
	} & \t{ % +as +ls
		\3{192},%	% +as +tr +ls
		\4{240}\\	% +lu +as +it +ls
	} & \t{ % +as +rs
		\3{204},%	% +as +tr +rs
		\4{248}\\	% +ru +as +it +rs
	}
\\\hline
Connex & co%
	{%
	} & \t{ % +co +s1
	} & \t{ % +co +ir
		%\2{040}\\	% +co +ir
	} & \t{ % +co +rf
	} & \t{ % +co +qr
	} & \t{ % +co +an
	} & \t{ % +co +sc
	} & \t{ % +co +it
	} & \t{ % +co +s2
	} & \t{ % +co +qt
	} & \t{ % +co +de
	} & \t{ % +co +ls
	} & \t{ % +co +rs
	}
\\\hline
Trans & tr%
	{%
	} & \t{ % +tr +s1
		\5{253}\\	% +le +tr +s1 -it +ls
		\5{256}\\	% +re +tr +s1 -it +rs
	} & \t{ % +tr +ir
		\3{129}\\	% -as +tr +ir
	} & \t{ % +tr +rf
	} & \t{ % +tr +qr
	} & \t{ % +tr +an
	} & \t{ % +tr +sc
		%\3{158}\\	% +tr -s1 +sc
	} & \t{ % +tr +it
		\3{166}\\	% +tr -s1 +it
	} & \t{ % +tr +s2
		\4{233},%	% +tr -s1 +s2 +de
		\4{234}\\	% +tr -it +s2 +de
	} & \t{ % +tr +qt
	} & \t{ % +tr +de
		\4{233}\\	% +tr -s1 +s2 +de
		\4{234}\\	% +tr -it +s2 +de
	} & \t{ % +tr +ls
		\3{191},%	% +at +tr +ls
		\3{192}\\	% +as +tr +ls
		\4{236},%	% -le +lu +tr +ls
		\5{253}\\	% +le +tr +s1 -it +ls
	} & \t{ % +tr +rs
		\3{203},%	% +at +tr +rs
		\3{204}\\	% +as +tr +rs
		\4{245},%	% -re +ru +tr +rs
		\5{256}\\	% +re +tr +s1 -it +rs

%%%%%%%% end of table + + %%%%%%%%
%%%%%%%% rows (em un cr le re lu ru sy at as co tr) %%%%%%%%
%%%%%%%% cols (s1 ir rf qr an sc it s2 qt de ls rs) %%%%%%%%

%%%%%%%% start of table - - %%%%%%%%
%%%%%%%% rows (s1 ir rf qr an sc it s2 qt de ls rs) %%%%%%%%
%%%%%%%% cols (s1 ir rf qr an sc it s2 qt de ls rs) %%%%%%%%

	}
\\\hline
\hline
SemiOrd1 & s1%
	{%
	} & \t{ % -s1 -s1
		\x		% (inverting)
	} & \t{ % +s1 +ir
		%\3{130}\\	% -as +s1 +ir
	} & \t{ % +s1 +rf
		\3{135}\\	% +cr +s1 +rf
		\3{136}\\	% -co +s1 +rf
	} & \t{ % +s1 +qr
		\4{219}\\	% -le +sy +s1 +qr
	} & \t{ % +s1 +an
	} & \t{ % +s1 +sc
	} & \t{ % +s1 +it
	} & \t{ % +s1 +s2
	} & \t{ % +s1 +qt
	} & \t{ % +s1 +de
	} & \t{ % +s1 +ls
		\3{193},%	% +cr +s1 +ls
		\3{194}\\	% +ru +s1 +ls
		\4{241},%	% +lu +s1 -it +ls
		\5{253}\\	% +le +tr +s1 -it +ls
	} & \t{ % +s1 +rs
		\3{205},%	% +cr +s1 +rs
		\3{206}\\	% +lu +s1 +rs
		\4{249},%	% +ru +s1 -it +rs
		\5{256}\\	% +re +tr +s1 -it +rs
	}
\\\hline
Irrefl & ir%
	{%
	} & \t{ % -ir -s1
	} & \t{ % -ir -ir
		\x		% (inverting)
	} & \t{ % +ir +rf
		\2{046}\\	% +ir +rf
	} & \t{ % +ir +qr
		\3{144}\\	% -em +ir +qr
	} & \t{ % +ir +an
		\3{153}\\	% -as +ir +an
	} & \t{ % +ir +sc
	} & \t{ % +ir +it
	} & \t{ % +ir +s2
	} & \t{ % +ir +qt
	} & \t{ % +ir +de
	} & \t{ % +ir +ls
	} & \t{ % +ir +rs
	}
\\\hline
Refl & rf%
	{%
	} & \t{ % -rf -s1
	} & \t{ % -rf -ir
	} & \t{ % -rf -rf
		\x		% (inverting)
	} & \t{ % +rf +qr
	} & \t{ % +rf +an
	} & \t{ % +rf +sc
		%\3{159}\\	% -co +rf +sc
	} & \t{ % +rf +it
		\3{167}\\	% -co +rf +it
	} & \t{ % +rf +s2
		\3{175}\\	% -co +rf +s2
	} & \t{ % +rf +qt
	} & \t{ % +rf +de
	} & \t{ % +rf +ls
	} & \t{ % +rf +rs
	}
\\\hline
QuasiRefl & qr%
	{%
	} & \t{ % -qr -s1
	} & \t{ % -qr -ir
	} & \t{ % -qr -rf
	} & \t{ % -qr -qr
		\x		% (inverting)
	} & \t{ % +qr +an
	} & \t{ % +qr +sc
		\3{160}\\	% -co +qr +sc
	} & \t{ % +qr +it
		\3{168},%	% -s1 +qr +it
		\4{223}\\	% -em -co +qr +it
	} & \t{ % +qr +s2
		\3{176},%	% -s1 +qr +s2
		\3{180}\\	% +qr -it +s2
	} & \t{ % +qr +qt
	} & \t{ % +qr +de
	} & \t{ % +qr +ls
		\3{197}\\	% -rf +qr +ls
	} & \t{ % +qr +rs
		\3{209}\\	% -rf +qr +rs
	}
\\\hline
AntiSym & an%
	{%
	} & \t{ % -an -s1
	} & \t{ % -an -ir
	} & \t{ % -an -rf
	} & \t{ % -an -qr
	} & \t{ % -an -an
		\x		% (inverting)
	} & \t{ % +an +sc
	} & \t{ % +an +it
	} & \t{ % +an +s2
	} & \t{ % +an +qt
		\3{182}\\	% -tr +an +qt
	} & \t{ % +an +de
	} & \t{ % +an +ls
	} & \t{ % +an +rs
	}
\\\hline
SemiConnex & sc%
	{%
	} & \t{ % -sc -s1
	} & \t{ % -sc -ir
	} & \t{ % -sc -rf
	} & \t{ % -sc -qr
	} & \t{ % -sc -an
	} & \t{ % -sc -sc
		\x		% (inverting)
	} & \t{ % +sc +it
	} & \t{ % +sc +s2
	} & \t{ % +sc +qt
	} & \t{ % +sc +de
	} & \t{ % +sc +ls
	} & \t{ % +sc +rs
	}
\\\hline
IncTrans & it%
	{%
	} & \t{ % -it -s1
	} & \t{ % -it -ir
	} & \t{ % -it -rf
	} & \t{ % -it -qr
	} & \t{ % -it -an
	} & \t{ % -it -sc
	} & \t{ % -it -it
		\x		% (inverting)
	} & \t{ % +it +s2
	} & \t{ % +it +qt
	} & \t{ % +it +de
	} & \t{ % +it +ls
		\3{200},%	% +ru +it +ls
		\4{240}\\	% +lu +as +it +ls
	} & \t{ % +it +rs
		\3{212},%	% +lu +it +rs
		\4{248}\\	% +ru +as +it +rs
	}
\\\hline
SemiOrd2 & s2%
	{%
	} & \t{ % -s2 -s1
	} & \t{ % -s2 -ir
	} & \t{ % -s2 -rf
	} & \t{ % -s2 -qr
	} & \t{ % -s2 -an
	} & \t{ % -s2 -sc
	} & \t{ % -s2 -it
	} & \t{ % -s2 -s2
		\x		% (inverting)
	} & \t{ % +s2 +qt
		\6{265}\\	% +at -it +s2 +qt +ls +rs
	} & \t{ % +s2 +de
		\4{233}\\	% +tr -s1 +s2 +de
		\4{234}\\	% +tr -it +s2 +de
	} & \t{ % +s2 +ls
		\6{265}\\	% +at -it +s2 +qt +ls +rs
	} & \t{ % +s2 +rs
		\6{265}\\	% +at -it +s2 +qt +ls +rs
	}
\\\hline
QuasiTrans & qt%
	\t{\\ \\%
	} & \t{ % -qt -s1
	} & \t{ % -qt -ir
	} & \t{ % -qt -rf
	} & \t{ % -qt -qr
	} & \t{ % -qt -an
	} & \t{ % -qt -sc
	} & \t{ % -qt -it
	} & \t{ % -qt -s2
	} & \t{ % -qt -qt
		\x		% (inverting)
	} & \t{ % +qt +de
	} & \t{ % +qt +ls
		\4{244},%	% +ru -sy +qt +ls
		\6{265}\\	% +at -it +s2 +qt +ls +rs
	} & \t{ % +qt +rs
		\6{265}\\	% +at -it +s2 +qt +ls +rs
	}
\\\hline
Dense & de%
	{%
	} & \t{ % -de -s1
	} & \t{ % -de -ir
	} & \t{ % -de -rf
	} & \t{ % -de -qr
	} & \t{ % -de -an
	} & \t{ % -de -sc
	} & \t{ % -de -it
	} & \t{ % -de -s2
	} & \t{ % -de -qt
	} & \t{ % -de -de
		\x		% (inverting)
	} & \t{ % +de +ls
	} & \t{ % +de +rs
	}
\\\hline
LfSerial & ls%
	\t{\\ \\%
	} & \t{ % -ls -s1
	} & \t{ % -ls -ir
	} & \t{ % -ls -rf
	} & \t{ % -ls -qr
	} & \t{ % -ls -an
	} & \t{ % -ls -sc
	} & \t{ % -ls -it
	} & \t{ % -ls -s2
	} & \t{ % -ls -qt
	} & \t{ % -ls -de
	} & \t{ % -ls -ls
		\x		% (inverting)
	} & \t{ % +ls +rs
		\6{265}\\	% +at -it +s2 +qt +ls +rs
	}
\\\hline
RgSerial & rs%
	\t{\\%
	} & \t{ % -rs -s1
	} & \t{ % -rs -ir
	} & \t{ % -rs -rf
	} & \t{ % -rs -qr
	} & \t{ % -rs -an
	} & \t{ % -rs -sc
	} & \t{ % -rs -it
	} & \t{ % -rs -s2
	} & \t{ % -rs -qt
	} & \t{ % -rs -de
	} & \t{ % -rs -ls
	} & \t{ % -rs -rs
		\x		% (inverting)

	}
\\\hline
\end{tabular}
\caption{Law index rg (bot lf: $A \lor B$ required, top rg: $A \land B$ impossible)}
%\rule{\textwidth}{0.1mm}
\LABEL{Law index - - / + + rg}
\end{center}
\end{figure}

%%%%%%%% end of table - - %%%%%%%%
%%%%%%%% rows (s1 ir rf qr an sc it s2 qt de ls rs) %%%%%%%%
%%%%%%%% cols (s1 ir rf qr an sc it s2 qt de ls rs) %%%%%%%%

}

\begin{figure}
\begin{center}
\includegraphics[width=\textwidth]{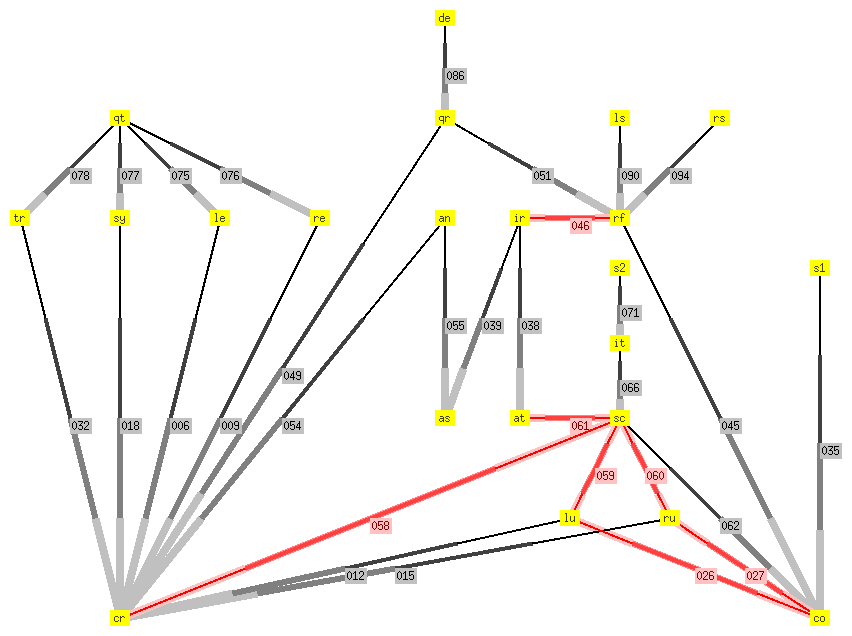}
\caption{Implications (black) and incompatibilities (red)
	between properties}
\LABEL{Implications and incompatibilities between properties}
\end{center}
\end{figure}

\clearpage
\section{Formal proofs of property laws}
\LABEL{Formal proofs of property laws}

Most of the law suggestions from Fig.~\REF{Law index + - lf}
to~\REF{Law index - - / + + rg}
could be proven to hold for all relations, on finite and on infinite
sets.
%~
Some suggestions turned out to hold only for relations on a
sufficiently large\footnote{
	For the following laws, we need a universe of at least
	$2$ elements:
		\QRef{025},
		\QRef{058},
		\QRef{135},
		\QRef{160},
		\QRef{161},
		\QRef{169},
		\QRef{193},
		\QRef{194},
		\QRef{198},
		\QRef{199},
		\QRef{205},
		\QRef{206},
		\QRef{210},
		\QRef{211};
	$3$ elements:
		Lem.~\REF{eucl 7},
		\QRef{188},
		\QRef{226};
	$4$ elements:
		\QRef{028},
		\QRef{059},
		\QRef{060},
		\QRef{061};
	$5$ elements:
		\QRef{200},
		\QRef{212},
		\QRef{220},
		\QRef{240},
		\QRef{248},
		\QRef{254};
	finite cardinality:
		Lem.~\REF{ser 5},
		\QRef{190},
		\QRef{202},
		\QRef{213},
		\QRef{214};
	finite and odd cardinality:
		\QRef{235}.
}
set $X$.
%~
Suggestion \QRef{235} turned out to hold only for a finite set
$X$ of odd cardinality (cf.\ Lem.~\REF{uniq 3}).

We considered all these suggestions to be laws, when appropriate
cardinality restrictions are added; their proofs are given in this
section.\footnote{
	We also gave proofs for well-known laws, and even for trivial ones.
}
%~
All remaining suggestions were considered non-laws; they are discussed
in section ``Examples'' (\REF{Examples}).

We loosely grouped the proven laws by some ``main property'', usually
the most unfamiliar property; for example,
Lem.~\REF{quasiTrans 4} relates
symmetry, transitivity, and quasi-transitivity, it is shown in
the ``Quasi-transitivity'' section (\REF{Quasi-transitivity}).
%~
Sometimes, we listed a result multiply, accepting some redundancy as a
trade-off for local completeness.
%~
The grouping is still far from being objective, and it is
doubtful that the latter is possible at all.

Due to the grouping we have some forward references in our proofs.
%~
For example, the proof of Lem.~\REFF{corefl 2}{2} uses Lem.~\REF{ser 1}.
%~
In order to establish the absence of cycles, we computed the
proof depth of each lemma to be one more than the maximal proof
depth of all its used lemmas.
%~
If a lemma would refer to itself directly or indirectly in its proof,
no finite proof depth could be assigned to it.
%~
We indicate the proof depth by a small superscript, e.g.\
``Lemma~\REF{corefl 2}. \ldots\ \REF{corefl 2 2}.\dpd{4}''
indicates that Lem.~\REFF{corefl 2}{2} has proof depth $4$.

\subsection{Co-reflexivity}
\LABEL{Co-reflexivity}

\LEMMA{}{%
\LABEL{corefl 1}\dpd{1}%
The union of a co-reflexive relation and a transitive relation
is always transitive.
}
\PROOF{
Let $C$ be co-reflexive and $T$ be transitive.
%~
Let $R = C \cup T$.
%~
Assume $x R y \land y R z$.
%~
We distinguish four cases:
\begin{enumerate}
\item
	If $x T y \land y T z$,
	then $x T z$ by transitivity of $T$,
	and hence $x R z$.
\item
	If $x T y \land y C z$,
	then $y = z$ by co-reflexivity of $C$,
	hence $x T z$ by substitutivity,
	hence $x R z$.
\item
	Similarly,
	$x C y \land y T z \Ra x = y T z \Ra x R z$.
\item
	If $x C y \land y C z$,
	then $x = y C z$ implies $x R z$.
\qed
\end{enumerate}
}

\LEMMA{Identity relation}{%
\LABEL{corefl 2}%
Given a set $X$,
the identity relation $I = \set{ \tpl{x,x} \mid x \in X }$
is uniquely characterized by any of the following properties
(\QRef{195}, \QRef{207}):
\begin{enumerate}
\item\LABEL{corefl 2 1}\dpd{1}%
	It is the only relation on $X$ that is both co-reflexive and
	reflexive.
\item\LABEL{corefl 2 2}\dpd{4}%
	It is the only relation on $X$ that is both co-reflexive and
	left serial.
\item\LABEL{corefl 2 3}\dpd{4}%
	It is the only relation on $X$ that is both co-reflexive and
	right serial.
\end{enumerate}
It has the following properties:
\begin{enumerate}
\setcounter{enumi}{3}
\item\LABEL{corefl 2 4}\dpd{1}%
	It doesn't satisfy semi-order property 1,
	if $X$ has at least $2$ elements
	(\QRef{135}, \QRef{193}, \QRef{205}).
\end{enumerate}
}
\PROOF{
\begin{enumerate}
\item%corefl 2 1
	The conjunction of
	Def.~\REFF{def}{Refl} and~\REFF{def}{CoRefl} is
	$\forall x,y \in X. \;\; xRy \lra x=y$
	which is the defining condition of $I$.
\item%corefl 2 2
	$I$ is left serial by Lem.~\REF{ser 1}\dpr{3}.
	%~
	If some relation $R$ is both co-reflexive and left serial,
	then $\forall x \in X \; \exists x' \in X. \;\; x' R x$
	holds; hence $\forall x \in X . \;\; x R x$ holds;
	therefore $R = I$ by case~\REF{corefl 2 1}\dpr{1}.
\item%corefl 2 3
	Dual to case~\REF{corefl 2 2}.
\item%corefl 2 4
	Let $x \neq y$, then $xIx$, $yIy$, and $x,y$ are incomparable
	w.r.t.\ $I$.
	%~
	If $I$ would satisfy semi-order property 1,
	then $xIy$ would hold, contradicting our assumption.
\qed
\end{enumerate}
}

\LEMMA{Sufficient for co-reflexivity implying emptiness}{%
\LABEL{corefl 3}%
On a set $X$ of at least $2$ elements,
a co-reflexive relation $R$
is empty if one of the following sufficient conditions holds:
\begin{enumerate}
\item\LABEL{corefl 3 1}\dpd{1}
	$R$ satisfies semi-order property 2 (\QRef{169});
\item\LABEL{corefl 3 2}\dpd{1}
	$R$ is incomparability-transitive (\QRef{161}).
\end{enumerate}
}
\PROOF{
Assume for contradiction $xRy$, then $x=y$.
%~
Let $w \neq x$.
\begin{enumerate}
\item%corefl 3 1
	Applying semi-order property 2 to $xRx \land xRx$ and $w$
	yields the contradiction $wRx \lor xRw$, i.e. $w=x$.
\item%corefl 3 2
	Applying incomparability-transitivity to $xRx$ and $w$
	yields the contradiction that $x,x$ must be incomparable.
\qed
\end{enumerate}
}

\LEMMA{Sufficient for co-reflexivity}{%
\LABEL{corefl 6}%
A relation $R$ is co-reflexive if one of the following sufficient
conditions holds:
\begin{enumerate}
\item\LABEL{corefl 6 1}\dpd{1}%
	$R$ is right quasi-reflexive and left unique (\QRef{140});
\item\LABEL{corefl 6 2}\dpd{1}%
	$R$ is left quasi-reflexive and right unique (\QRef{141});
\item\LABEL{corefl 6 3}\dpd{2}%
	$R$ is right Euclidean and left unique (\QRef{095});
\item\LABEL{corefl 6 4}\dpd{2}%
	$R$ is left Euclidean and right unique (\QRef{096});
\item\LABEL{corefl 6 4a}\dpd{2}%
	$R$ is reflexive and left unique (\QRef{133});
\item\LABEL{corefl 6 4b}\dpd{2}%
	$R$ is reflexive and right unique (\QRef{134}); or
\item\LABEL{corefl 6 5}\dpd{1}%
	$R$ is symmetric and anti-symmetric (\QRef{149}).
\end{enumerate}
}
\PROOF{
\begin{enumerate}
\item%corefl 6 1
	If $xRy$, then $yRy$ by right quasi-reflexivity,
	hence $x=y$ by left uniqueness.
\item%corefl 6 2
	Dual to~\REF{corefl 6 1}.
\item%corefl 6 3
	Follows from~\REF{corefl 6 1}\dpr{1},
	since right Euclidean relation is right quasi-reflexive
	by Lem.~\REF{eucl 11}\dpr{1}.
\item%corefl 6 4
	Dual to~\REF{corefl 6 3}.
\item%corefl 6 4a
	Follows from~\REF{corefl 6 1}\dpr{1}
	and Lem.~\REF{quasiRefl 2}\dpr{1}.
\item%corefl 6 4b
	Dual to~\REF{corefl 6 4a}.
\item%corefl 6 5
	Assume for contradiction $xTy$ holds for some $x \neq y$.
	%~
	Then $yRx$ by symmetry, while $\lnot yRx$ by anti-symmetry.
\qed
\end{enumerate}
}

\LEMMA{Necessary for co-reflexivity}{%
\LABEL{corefl 5}%
If a relation $R$ on a set $X$ is co-reflexive, then it satisfies the
following necessary conditions:
\begin{enumerate}
\item\LABEL{corefl 5 1}\dpd{1}%
	$R$ is left Euclidean (\QRef{006});
\item\LABEL{corefl 5 2}\dpd{1}%
	$R$ is right Euclidean (\QRef{009});
\item\LABEL{corefl 5 3}\dpd{1}%
	$R$ is left unique (\QRef{012});
\item\LABEL{corefl 5 4}\dpd{1}%
	$R$ is right unique (\QRef{015});
\item\LABEL{corefl 5 5}\dpd{1}%
	$R$ is symmetric (\QRef{018});
\item\LABEL{corefl 5 6}\dpd{1}%
	$R$ is anti-symmetric (\QRef{054});
\item\LABEL{corefl 5 7}\dpd{3}%
	$R$ is transitive (\QRef{032}),
	hence quasi-transitive (\QRef{074});
\item\LABEL{corefl 5 9}\dpd{3}%
	$R$ is not semi-connex (\QRef{058}),
	hence not connex (\QRef{025}),
	provided $X$ has at least $2$ elements.
\item\LABEL{corefl 5 10}\dpd{1}%
	$R$ is quasi-reflexive (\QRef{049});
\item\LABEL{corefl 5 8}\dpd{2}%
	$R$ is dense (\QRef{081}).
\end{enumerate}
}
\PROOF{
\begin{enumerate}
\item%corefl 5 1
	If $yRx$ and $zRx$, then $y=x=z$, hence $yRz$.
\item%corefl 5 2
	Dual to case~\REF{corefl 5 1}.
\item%corefl 5 3
	If $x_1 R y \land x_2 R y$, then $x_1 = y = x_2$.
\item%corefl 5 4
	Dual to case~\REF{corefl 5 3}.
\item%corefl 5 5
	If $xRy$, then $x=y$, hence $yRx$.
\item%corefl 5 6
	If $xRy$ [and $yRx$], then $x=y$.
\item%corefl 5 7
	If $xRy \land yRz$, the $x=y=z$, hence $xRz$.
	%~
	Quasi-transitivity follows by Lem.~\REF{quasiTrans 4}\dpr{2}.
\item%corefl 5 9
	If $x \neq y$, then neither $xRy$ nor $yRx$ is possible.
	%~
	By Lem.~\REF{conn 1}\dpr{2}, $R$ can't be connex either.
\item%corefl 5 10
	If $xRy$, then $x=y$, hence $xRx$ and $yRy$.
\item%corefl 5 8
	Follows from~\REF{corefl 5 10}\dpr{1}, since a quasi-reflexive
	relation is always dense by Lem.~\REFF{den 1}{3}\dpr{1}.
\qed
\end{enumerate}
}

%\clearpage
\subsection{Reflexivity}
\LABEL{Reflexivity}

\LEMMA{Necessary for reflexivity}{%
\LABEL{quasiRefl 2}\dpd{1}%
A reflexive relation is always quasi-reflexive (\QRef{051}).
}
\PROOF{
If $xRy$ holds, and even if not, then $xRx \land yRy$ holds by reflexivity.
\qed
}

\LEMMA{Incompatibilities of reflexivity}{%
\LABEL{irrefl 2}\dpd{3}%
An irreflexive relation $R$ on a non-empty $X$ cannot be reflexive
(\QRef{046}).
%~
As a consequence, an anti-transitive relation on $X$ can't be
reflexive (\QRef{043}), and neither can an asymmetric relation on $X$
(\QRef{044}).
%~
Moreover, an irreflexive (\QRef{040}), anti-transitive (\QRef{028}), or
asymmetric (\QRef{029}) relation cannot be connex.
}
\PROOF{
Let $x \in X$, then $\lnot xRx$ by irreflexivity, and $xRx$ by
reflexivity.
%~
By Lem.~\REF{antitrans 1}\dpr{1}, each anti-transitive relation is
irreflexive.
%~
By Lem.~\REFF{asym 1b}{1}\dpr{1}, each asymmetric relation is irreflexive.
%~
By Lem.~\REF{conn 1}\dpr{2}, each connex relation is reflexive.
\qed
}

%\clearpage
\subsection{Irreflexivity}
\LABEL{Irreflexivity}

\LEMMA{Sufficient for irreflexivity implying emptiness}{%
\LABEL{irrefl 1}%
An irreflexive relation $R$ needs to be empty
if one of the following sufficient conditions is satisfied:
\begin{enumerate}
\item\LABEL{irrefl 1 1}\dpd{1}%
	$R$ is co-reflexive (\QRef{124});
\item\LABEL{irrefl 1 2}\dpd{1}%
	$R$ is left quasi-reflexive (\QRef{144});
\item\LABEL{irrefl 1 3}\dpd{1}%
	$R$ is right quasi-reflexive (\QRef{144});
\item\LABEL{irrefl 1 4}\dpd{2}%
	$R$ is left Euclidean (\QRef{125});
\item\LABEL{irrefl 1 5}\dpd{2}%
	$R$ is right Euclidean (\QRef{126}).
\end{enumerate}
As a consequence, an anti-transitive relation
(\QRef{102}, \QRef{142}, \QRef{103}, \QRef{104})
as well as an asymmetric relation
(\QRef{105}, \QRef{143}, \QRef{106}, \QRef{107})
needs to be empty under
the same sufficient conditions.
}
\PROOF{
For case~\REF{irrefl 1 1}
to~\REF{irrefl 1 3},
assume for contradiction that $R$ is irreflexive and $aRb$ holds,
i.e.\ $a \neq b$.
\begin{enumerate}
\item%irrefl 1 1
	Then $a=b$ by co-reflexivity, contradicting our assumption.
\item%irrefl 1 2
	Then $aRa$, contradicting irreflexivity.
\item%irrefl 1 3
	Then $bRb$, contradicting irreflexivity.
\item%irrefl 1 4
	Follows from~\REF{irrefl 1 2}\dpr{1}
	using Lem.~\REF{eucl 11}\dpr{1}.
\item%irrefl 1 5
	Follows from~\REF{irrefl 1 3}\dpr{1}
	using Lem.~\REF{eucl 11}\dpr{1}.
\end{enumerate}
Each anti-transitive relation is irreflexive by
Lem.~\REF{antitrans 1}\dpr{1}.
%~
Each asymmetric relation is irreflexive by
Lem.~\REFF{asym 1b}{1}\dpr{1}.
\qed
}

%\clearpage
\subsection{Asymmetry}
\LABEL{Asymmetry}

\LEMMA{Sufficient for asymmetry}{%
\LABEL{asym 1a}%
A relation $R$
is asymmetric if one of the following sufficient conditions holds:
\begin{enumerate}
\item\LABEL{asym 1a 1}\dpd{1}%
	$R$ is irreflexive and anti-symmetric (\QRef{153});
\item\LABEL{asym 1a 2}\dpd{1}%
	$R$ is irreflexive and transitive (\QRef{129});
\item\LABEL{asym 1a 4}\dpd{2}%
	$R$ is irreflexive and satisfies semi-order property 1
	(\QRef{130}).
\item\LABEL{asym 1a 3}\dpd{2}%
	$R$ is anti-transitive and anti-symmetric (\QRef{150});
\item\LABEL{asym 1a 6}\dpd{2}%
	$R$ is anti-transitive and transitive (\QRef{118});
	or
\item\LABEL{asym 1a 5}\dpd{3}%
	$R$ is anti-transitive and satisfies semi-order property 1
	(\QRef{119});
\end{enumerate}
}
\PROOF{
\begin{enumerate}
\item%asym 1a 1
	If $xRy$,
	then $x \neq y$ by irreflexivity,
	hence, $\lnot yRx$ by anti-symmetry.
\item%asym 1a 2
	Let $xRy$ hold.
	%~
	If $yRx$, then $xRx$ by transitivity,
	which contradicts irreflexivity.
\item%asym 1a 4
	Follows from~\REF{asym 1a 2}\dpr{1}
	since $R$ is transitive
	by Lem.~\REFF{semiOrd1 12}{3}\dpr{1}.
\item%asym 1a 3
	Follows from~\REF{asym 1a 1}\dpr{1},
	since $R$ is irreflexive by Lem.~\REF{antitrans 1}\dpr{1}.
\item%asym 1a 6
	Follows from~\REF{asym 1a 2}\dpr{1},
	by the same argument.
	%~
	Anti-transitivity as well as transitivity is vacuous in this case;
	cf.\ Lem.~\REF{antitrans 4}.
\item%asym 1a 5
	Follows from~\REF{asym 1a 4}\dpr{2},
	by the same argument.
\qed
\end{enumerate}
}

\LEMMA{Necessary for asymmetry}{%
\LABEL{asym 1b}%
Let $R$ be asymmetric.
%~
Then $R$ is necessarily
\begin{enumerate}
\item\LABEL{asym 1b 1}\dpd{1}%
	irreflexive (\QRef{039}); and
\item\LABEL{asym 1b 2}\dpd{1}%
	anti-symmetric (\QRef{055}).
\end{enumerate}
}
\PROOF{
\begin{enumerate}
\item%asym 1b 1
	$xRx$ would imply the contradiction $\lnot xRx$.
\item%asym 1b 2
	$xRy \lor x \neq y$ implies $\lnot yRx$,
	since its left disjunct does.
\qed
\end{enumerate}
}

\LEMMA{Incompatibilities of asymmetry}{%
\LABEL{asym 3}\dpd{2}%
On a finite set $X$, an asymmetric and transitive relation
can neither be
left (\QRef{192}) nor right (\QRef{204}) serial.
%~
On the infinite set $\Z$
of integer numbers, the usual order $<$ satisfies
all four properties simultaneously.
}
\PROOF{
By Lem.~\REFF{asym 1b}{1}\dpr{1},
such a relation is an irreflexive partial order.
%~
On a finite set, it must have a smallest and a largest element;
thus it can't be serial.
\qed
}

%\clearpage
\subsection{Symmetry}
\LABEL{Symmetry}

\LEMMA{}{%
\LABEL{sym 1}%
A symmetric relation $R$ on a set $X$ is
\begin{enumerate}
\item\LABEL{sym 1 1}\dpd{1}%
	left quasi-reflexive iff it is right quasi-reflexive;
\item\LABEL{sym 1 2}\dpd{1}%
	left Euclidean iff it is right Euclidean (\QRef{098}, \QRef{099});
\item\LABEL{sym 1 3}\dpd{1}%
	left serial iff it is right serial (\QRef{215}, \QRef{216});
\item\LABEL{sym 1 4}\dpd{1}%
	left unique iff it is right unique (\QRef{100}, \QRef{101}).
\end{enumerate}
}
\PROOF{
Let $R$ be symmetric.
\begin{enumerate}
\item%sym 1 1
	If $R$ is left quasi-reflexive
	and $xRy$ holds,
	then $yRx$ by symmetry,
	hence $yRy$.
\item%sym 1 2
	If $R$ is left Euclidean,
	and $xRy$ and $xRz$ holds,
	then $yRx$ and $zRx$ by symmetry,
	hence $yRz$ by left Euclideanness.
\item%sym 1 3
	If $R$ is left serial
	and $x \in X$,
	then find some $y$ with $yRx$ by left seriality
	hence $xRy$ by symmetry.
\item%sym 1 4
	If $R$ is left unique,
	and $x R y_1$ and $x R y_2$ holds,
	then $y_1 R x$ and $y_2 R x$ by symmetry,
	hence $y_1 = y_2$ by left uniqueness.
\end{enumerate}
The converse directions are shown similarly.
\qed
}

\LEMMA{}{%
\LABEL{sym 2}\dpd{1}%
A symmetric and asymmetric relation is always
empty (\QRef{108}).
}
\PROOF{
Assume for contradiction $xRy$ holds.
%~
Then $yRx$ by symmetry, and $\lnot yRx$ by asymmetry.
\qed
}

%\clearpage
\subsection{Quasi-transitivity}
\LABEL{Quasi-transitivity}

\LEMMA{}{%
\LABEL{quasiTrans 1}%
\begin{enumerate}
\item\LABEL{quasiTrans 1 1}\dpd{1}%
	$R$ is a quasi-transitive relation
	iff
	$R = I \disjUnion P$
	for some symmetric relation $I$
	and some transitive relation $P$,
	where ``$\disjUnion$'' denotes the disjoint union..
\item\LABEL{quasiTrans 1 2}\dpd{1}%
	$I$ and $P$ are not uniquely determined by a given $R$.
\item\LABEL{quasiTrans 1 3}\dpd{1}%
	The definitions
	$x I y :\Lra x R y \land y R x$
	and
	$x P y :\Lra x R y \land \lnot y R x$
	lead to the minimal $P$.
\end{enumerate}
}
\PROOF{
\begin{enumerate}
\item%quasiTrans 1 1
	``$\Ra$'':
	%~
	Let $R$ be quasi-transitive.
	%~
	Following \CITEPAGE{p.381}{Sen.1969},
	define $x I y :\Lra x R y \land y R x$
	and $x P y :\Lra x R y \land \lnot y R x$.
	%~
	Then
	\begin{itemize}
	\item $I$ and $P$ are disjoint:
		\\
		$\begin{array}{cll}
		& x I y \land x P y	\\
		\Ra & y R x \land \lnot y R x
			& \mbox{ using the definitions of } I
			\mbox{ and } P	\\
		\Ra & \false	\\
		\end{array}$
	\item Their union is $R$:
		\\
		$\begin{array}{cll}
		& x I y \lor x P y	\\
		\Lra & (x R y \land y R x) \lor (x R y \land \lnot y R x)
			& \mbox{ by definition of } I \mbox{ and } P	\\
		\Lra & x R y \land (y R x \lor \lnot y R x)
			& \mbox{ by distributivity}	\\
		\Lra & x R y	\\
		\end{array}$
	\item $I$ is symmetric:
		\\
		$\begin{array}{cll}
		& x I y	\\
		\Ra & x R y \land y R x	\\
		\Ra & y I x	\\
		\end{array}$
	\item $P$ is transitive:
		\\
		$\begin{array}{cll}
		& x P y \land y P z	\\
		\Ra & x R y \land \lnot y R x \land y R z \land \lnot z R y
			& \mbox{ by definition of } P	\\
		\Ra & x R z \land \lnot z R x
			& \mbox{ by quasi-transitivity of } R	\\
		\Ra & x P z
			& \mbox{ by definition of } P	\\
		\end{array}$
	\end{itemize}
	``$\La$'':
	%~
	Let $R = I \disjUnion P$ for some symmetric relation $I$ and some
	transitive relation $R$.
	%~
	Assume $x R y$ and $y R z$ hold,
	but neither $y R x$ nor $z R y$ does.
	%~
	We observe the following facts:
	\begin{enumerate}
	\item\LABEL{quasiTrans 1 <= 1}%
		$x I y$ is false,
		since else $x I y \Ra y I x \Ra y R x$,
		contradicting our assumptions.
	\item\LABEL{quasiTrans 1 <= 2}%
		$x P y$ holds,
		since $x R y \Ra x I y \lor x P y \Ra x P y$
		by~\REF{quasiTrans 1 <= 1}.
	\item\LABEL{quasiTrans 1 <= 3}%
		$y P z$ follows by an argument similar
		to~\REF{quasiTrans 1 <= 1} and~\REF{quasiTrans 1 <= 2}.
	\item\LABEL{quasiTrans 1 <= 4}%
		Hence $x P z$ holds, by transitivity of $P$.
	\item\LABEL{quasiTrans 1 <= 5}%
		Hence $x R z$.
	\item\LABEL{quasiTrans 1 <= 6}%
		Since $I$ and $P$ are disjoint,
		we obtain $\lnot x I z$
		from~\REF{quasiTrans 1 <= 4};
		hence $\lnot z I x$ by symmetry of $I$.
	\item\LABEL{quasiTrans 1 <= 7}%
		Finally, we have $\lnot z R x$, since else
		$z P x$
		by~\REF{quasiTrans 1 <= 6},
		which in turn would imply
		$z P y$ by~\REF{quasiTrans 1 <= 2}
		and the transitivity of $P$,
		which would imply $z R y$,
		contradicting our assumptions.
	\end{enumerate}
	From~\REF{quasiTrans 1 <= 5} and~\REF{quasiTrans 1 <= 7},
	we conclude the quasi-transitivity of $R$.
\item%quasiTrans 1 2
	For example, if $R$ is an equivalence relation,
	$I$ may be chosen as the empty relation, or as $R$ itself,
	and $P$ as its complement.
\item%quasiTrans 1 3
	Given $R$, whenever $xRy \land \lnot yRx$ holds,
	the pair $\tpl{x,y}$ can't belong to
	the symmetric part, but must belong to the transitive part.
\qed
\end{enumerate}
}

\LEMMA{}{%
\LABEL{quasiTrans 4}\dpd{2}%
Each symmetric relation is quasi-transitive (\QRef{077});
each transitive relation is quasi-transitive (\QRef{078}).
}
\PROOF{
Follows from Lem.~\REF{quasiTrans 1}\dpr{1} and the
transitivity (Exm.~\REFF{empty}{7}\dpr{1})
and symmetry (\REFF{empty}{4}\dpr{1})
of the empty relation.
\qed
}

\LEMMA{}{%
\LABEL{quasiTrans 5}\dpd{2}%
A quasi-transitive relation is transitive if it is anti-symmetric
(\QRef{182}), hence in particular if it is asymmetric (\QRef{181}).
}
\PROOF{
Let $R$ be anti-symmetric and quasi-transitive.
%~
We use the definitions of $I$ and $P$ from
Lem.~\REFF{quasiTrans 1}{3}\dpr{1}.
%~
We have $x I y \Ra x R y \land y R x \Ra x = y$
by anti-symmetry, hence $I$ is co-reflexive.
%~
By Lem.~\REF{corefl 1}\dpr{1},
$R = I \cup P$ is transitive.
\qed
}

\LEMMA{}{%
\LABEL{quasiTrans 6}\dpd{2}%
If $P$ is a semi-order
(Def.~\REF{Kinds of binary relations}.\REF{semi-order})
then $R$ defined by
$xRy :\Lra xPy \lor (\lnot xPy \land \lnot yPx)$
is quasi-transitive.
%~
This may be meant by Sen's remark that semi-orders are a special case
of quasi-transitivity (\CITEPAGE{p.314}{Sen.1971}).
%~
Note that $P$ itself is transitive by Lem.~\REFF{semiOrd1 12}{5},
hence trivially quasi-transitive by Lem.~\REF{quasiTrans 4}.
}
\PROOF{
Let $P$ be a semi-order.
%~
Define $xIy :\Lra \lnot xPy \land \lnot yPx$,
then $I$ is symmetric, and disjoint from $P$.
%~
Since $P$ is asymmetric
by Def.~\REF{Kinds of binary relations}.\REF{semi-order},
it is irreflexive by Lem.~\REFF{asym 1b}{1}\dpr{1}
and hence transitive by~\REFF{semiOrd1 12}{3}\dpr{1}.
%~
Hence $R = I \disjUnion P$ is quasi-transitive
by Lem.~\REF{quasiTrans 1}\dpr{1}.
\qed
}

\LEMMA{Sufficient for quasi-transitivity implying symmetry}{%
\LABEL{quasiTrans 8}%
A quasi-transitive relation $R$ is symmetric if one of the following
sufficient conditions holds:
\begin{enumerate}
\item\LABEL{quasiTrans 8 1}\dpd{1}%
	$R$ is right unique and left serial (\QRef{244}); or
\item\LABEL{quasiTrans 8 2}\dpd{1}%
	$R$ is left unique and right serial.
\end{enumerate}
}
\PROOF{
We show~\REF{quasiTrans 8 1}; the proof of~\REF{quasiTrans 8 2} is similar.
%~
Let $yRz$ hold; assume for contradiction $\lnot zRy$.
%~
Obtain $xRy$ by left seriality.
%~
We distinguish two cases:
\begin{itemize}
\item $\lnot yRx$ holds.
	%~
	Then $xRz$ by quasi-transitivity,
	hence $y=z$ by right uniqueness,
	hence $zRy$, contradicting our assumption.
\item $yRx$ holds.
	%~
	Then $x=z$ by right uniqueness,
	hence $zRy$, contradicting our assumption.
\qed
\end{itemize}
}

%\clearpage
\subsection{Anti-transitivity}
\LABEL{Anti-transitivity}

\LEMMA{}{%
\LABEL{antitrans 1}\dpd{1}%
An anti-transitive relation is always irreflexive (\QRef{038}).
}
\PROOF{
Assume $xRx$ holds.
%~
Then $xRx \land xRx$ implies $\lnot xRx$ by anti-transitivity,
which is a contradiction.
\qed
}

\LEMMA{}{%
\LABEL{antitrans 3}\dpd{2}%
An irreflexive and left unique relation is always
anti-transitive; and so is  an irreflexive and right unique
relation (\QRef{127}, \QRef{128}).
%~
In particular, each asymmetric and left or right unique
relation is anti-transitive (\QRef{109}, \QRef{110}).
}
\PROOF{
Let $R$ be irreflexive and left unique,
assume for contradiction $xRy$, and $yRz$, but $xRz$.
%~
Then $x \neq y$ due to irreflexivity,
hence $yRz \land xRz$ contradicts left uniqueness.
%~
The proof for right uniqueness is similar.
%~
Each asymmetric relation is irreflexive
by Lem.~\REFF{asym 1b}{1}\dpr{1}.
\qed
}

\LEMMA{Necessary for transitivity and anti-transitivity}{%
\LABEL{antitrans 4}\dpd{3}%
On a nonempty set $X$, a relation that is both transitive
and anti-transitive
can for trivial reasons neither be left (\QRef{191}) nor right
(\QRef{203}) serial,
is must be asymmetric (\QRef{118}) and satisfy semi-order
property 2 (\QRef{172}).
}
\PROOF{
Let $R$ be transitive and anti-transitive,
then $xRy \land yRz$ cannot be satisfied for any $x,y,z$.
%~
Hence, $R$ vacuously satisfies semi-order property 2.
%~
If $R$ is left serial and $z \in X$, we have $yRz$ for some
$y$, and $xRy$ for some $x$, contradicting the above.
%~
Similarly, $R$ can't be right serial.
%~
Asymmetry has been shown in Lem.~\REFF{asym 1a}{6}\dpr{2}.
\qed
}

%\clearpage
\subsection{Incomparability-transitivity}
\LABEL{Incomparability-transitivity}

\LEMMA{}{%
\LABEL{incTrans 1}\dpd{3}%
Each semi-connex relation is incomparability-transitive (\QRef{066})
and hence satisfies semi-order property 2 (\QRef{070}).
%~
In particular, this applies to each connex relation
(\QRef{065}, \QRef{069}).
}
\PROOF{
\begin{itemize}
\item
	If $R$ is semi-connex and
	$x,y$ and
	$y,z$ are incomparable, then
	$x=y$ and $y=z$.
	%~
	Due to the latter, $x,z$ are incomparable.
\item
	Each incomparability-transitive relation
	satisfies semi-order property 2
	by Lem.~\REF{incTrans 13}\dpr{1}.
\item
	Each connex relation is semi-connex by
	Lem.~\REF{conn 1}\dpr{2}.
\qed
\end{itemize}
}

\LEMMA{}{%
\LABEL{incTrans 11}\dpd{1}%
If a relation is left Euclidean, left serial, and transitive, and
satisfies semi-order property 1, then it is also
incomparability-transitive (\QRef{253}).
%~
Dually, right Euclideanness, right seriality, transitivity, and
semi-order property 1 imply incomparability-transitivity (\QRef{256}).
}
\PROOF{
To show the first claim, assume for contradiction $R$ is not
incomparability-transitive.
%~
Let $aRb$ hold, and $c$ be incomparable both to $a$ and to $b$.
%~
By seriality, obtain $c' R c$.
%~
By semi-order property 1, $c' R b$ must hold.
%~
Hence, by Euclideanness, $a R c'$ holds.
%~
But transitivity then implies $aRc$, contradicting incomparability.

The proof for the dual claim is similar.
\qed
}

\LEMMA{}{%
\LABEL{incTrans 2}\dpd{1}%
Let $R$ be an incomparability-transitive relation on $X$.
%~
Whenever $xRx$ holds for some $x \in X$, then $x$ is comparable
to every $y \in X$.
%~
In particular, a reflexive relation $R$ can only be vacuously
incomparability-transitive,
that is, when $R$ is also connex (\QRef{167}).
}
\PROOF{
Let $xRx$ hold, let $y$ be arbitrary.
%~
If $x$ and $y$ were incomparable, then so were $y$ and $x$
due to symmetry,
hence also $x$ and $x$ by incomparability-transitivity,
contradicting $xRx$.
\qed
}

\LEMMA{}{%
\LABEL{incTrans 3}\dpd{4}%
If a left Euclidean is also incomparability-transitive,
then it is also transitive (\QRef{162})
and trivially satisfies semi-order property 1 (\QRef{164}),
moreover it is left serial or empty (\QRef{237}).
%~
The dual applies to a right Euclidean relation
(\QRef{163}, \QRef{165}, \QRef{246}).
}
\PROOF{
\begin{enumerate}
\item Transitivity:
	%~
	Let $R$ be left Euclidean; let $xRy$ and $yRz$ hold.
	%~
	We have $xRx$ by Lem.~\REF{eucl 3}\dpr{3},
	hence $x$ and $z$ are comparable
	by Lem.~\REF{incTrans 2}\dpr{1}.
	%~
	If $xRz$, we have transitivity immediately.
	%~
	Else, we have $zRx$, hence $x,y,z \in \dom(R)$, implying
	transitivity by Lem.~\REF{eucl 3}\dpr{3}.
\item Semi-order property 1:
	%~
	The antecedent of that property (Def.~\REFF{def}{SemiOrd1})
	cannot hold, since $yRz$ implies $yRy$
	by Lem.~\REF{eucl 3}\dpr{3},
	hence $y$ and $x$ are comparable
	by Lem.~\REF{incTrans 2}\dpr{1}.
\item Left seriality:
	%~
	Let $R$ be non-empty; let $aRb$ hold.
	%~
	An arbitrary $y$ must be comparable to $a$ or to $b$.
	%~
	If $aRy$ or $bRy$ holds, we are done immediately.
	%~
	If $yRa$ or $yRb$ holds, we have $yRy$
	by Lem.~\REF{eucl 3}\dpr{3}.
\end{enumerate}
The proof for a right Euclidean $R$ is similar.
\qed
}

\LEMMA{}{%
\LABEL{incTrans 6}\dpd{2}%
Each nonempty, quasi-reflexive and incomparability-transitive
relation is reflexive, and hence connex (\QRef{223});
i.e.\ its incomparability-transitivity is vacuous.
}
\PROOF{
Let $aRb$, hence also $aRa$ and $bRb$ hold.
%~
Let $x \in X$ be arbitrary.
%~
By Lem.~\REF{incTrans 2}\dpr{1}, $a$ is comparable to $x$.
%~
By quasi-reflexivity, $xRx$ holds.
%~
Hence $R$ is reflexive.
%~
Again by Lem.~\REF{incTrans 2}\dpr{1} we obtain that $R$ is connex.
\qed
}

\LEMMA{}{%
\LABEL{incTrans 7}\dpd{1}%
A symmetric and incomparability-transitive relation is
anti-transitive or dense (\QRef{232}).
}
\PROOF{
Let $R$ be symmetric, incomparability-transitive,
and not anti-transitive; let $aRb$, $bRc$, but $aRc$ hold.
%~
An arbitrary $x$ can be incomparable to at most one of $a,b,c$
(as a side remark: therefore $R$ needn't be semi-connex).
%~
If $xRy$ holds, then $x$ and $y$ must both be comparable to at
least one of $a,b,c$, we assume w.l.o.g\ $a$.
%~
Due to the symmetry of $R$, we have $xRa$ and $aRy$;
therefore $R$ is dense.
\qed
}

\LEMMA{}{%
\LABEL{incTrans 8}\dpd{1}%
A non-empty symmetric and incomparability-transitive relation is
always left (\QRef{238}) and right serial.
}
\PROOF{
Let $aRb$ hold.
%~
By symmetry, we have $bRa$.
%~
An arbitrary $y$ cannot be incomparable to both $a$ and $b$,
hence w.l.o.g.\ $aRy$, using symmetry.
\qed
}

\LEMMA{Necessary for uniqueness and incomparability-transitivity}{%
\LABEL{incTrans 9}%
Let a relation $R$ on a set
$X$ be left unique and incomparability-transitive.
%~
Then $R$ is necessarily
\begin{enumerate}
\item\LABEL{incTrans 9 2}\dpd{3}%
	left Euclidean or anti-transitive (\QRef{221});
\item\LABEL{incTrans 9 3}\dpd{2}%
	asymmetric or vacuously quasi-transitive (\QRef{228}); and
\item\LABEL{incTrans 9 4}\dpd{2}%
	asymmetric or left serial (\QRef{239}).
\end{enumerate}
Moreover, if $X$ has at least $5$ elements, then $R$ is necessarily
\begin{enumerate}
\setcounter{enumi}{3}
\item\LABEL{incTrans 9 1}\dpd{1}%
	empty or not right unique (\QRef{220}); and
\item\LABEL{incTrans 9 5}\dpd{2}%
	not both asymmetric and left serial (\QRef{240}).
\end{enumerate}
On the
%four-element
set $\set{a,b,c,d}$,
the relation
$\set{ \tpl{a,b}, \tpl{b,c}, \tpl{c,d}, \tpl{d,a} }$
is a counter-example for~\REF{incTrans 9 1} and~\REF{incTrans 9 5}.

Dually, let a relation $R$ on a set
$X$ be right unique and incomparability-transitive.
%~
Then $R$ is necessarily
\begin{enumerate}
\item right Euclidean or anti-transitive (\QRef{222});
\item asymmetric or vacuously quasi-transitive (\QRef{229}); and
\item asymmetric or right serial (\QRef{247}).
\end{enumerate}
If $X$ has at least $5$ elements, then $R$ is necessarily
\begin{enumerate}
\setcounter{enumi}{3}
\item empty or not left unique
	(coincides with dual, \QRef{220}); and
\item not both asymmetric and right serial (\QRef{248}).
\end{enumerate}
}
\PROOF{
\begin{enumerate}
\item%incTrans 9 2
	Shown in Lem.~\REFF{eucl 11b}{3}\dpr{2}.
\item%incTrans 9 3
	Follows from Lem.~\REFF{semiOrd2 7b}{3}\dpr{1},
	since each
	incomparability-transitive relation satisfies semi-order
	property 2 by Lem.~\REF{incTrans 13}\dpr{1}.
\item%incTrans 9 4
	Follows from Lem.~\REFF{semiOrd2 7b}{4}\dpr{1},
	since each
	incomparability-transitive relation satisfies semi-order
	property 2 by Lem.~\REF{incTrans 13}\dpr{1}.
\item%incTrans 9 1
	Let $a,b,x,y,z$ be five distinct elements of $X$,
	let $aRb$ hold.
	%~
	Consider the directed graph corresponding to $R$, with its
	vertices being the elements of $X$, and its edges being the
	pairs related by $R$.
	%~
	Due to the uniqueness properties, no two edges can go out from,
	or come in to, the same vertex.
	%~
	Hence, between $a$ and $x,y,z$, we can have at most one vertex
	(an incoming one).
	%~
	Similarly, between $b$ and $x,y,z$, we can have at most one
	vertex (an outgoing one).
	%~
	Hence, two of $x,y,z$ are unrelated to $a$, and two are
	unrelated to $b$.
	%~
	Hence, at least one of $x,y,z$ is unrelated to both $a$ and $b$.
	%~
	But this contradicts incomparability-transitivity.
\item%incTrans 9 5
	Define $\:,1n{x_\i}$ to be a cycle of length $n$
	if $x_i R x_{i+1}$ holds for $i=\:,1{n-1}\i$, and $x_n R x_1$
	holds.
	%~
	We first show the existence of a cycle of length $3$ or $4$.
	%~
	Using seriality, obtain $x_1, x_2, x_3, x_4$ such that
	$x_i R x_{i+1}$ $i=1,2,3$.
	%~
	We can't have $x_1 R x_3$, since this would imply $x_1 = x_2$,
	contradicting Lem.~\REFF{asym 1b}{1}\dpr{1}.
	%~
	If $x_3 R x_1$, we have a cycle of length $3$, and are done.
	%~
	Else, $x_1$ and $x_3$ are incomparable;
	hence due to incomparability-transitivity
	$x_1$ and $x_4$ can't be incomparable, too.
	%~
	We can't have $x_1 R x_4$, since then $x_1 = x_3$ by uniqueness,
	and we would have $x_1 R x_2$ and $x_2 R x_1$, contradicting
	asymmetry.
	%~
	Therefore, $x_4 R x_1$, and we have a cycle of length $4$.

	Now let a cycle $\:,1n{x_\i}$ of length $n=3$ or $n=4$ be given.
	%~
	Let $y$ be an element distinct from all cycle members.
	%~
	We can't have $y R x_i$ for any $i$ by uniqueness.
	%~
	We can't have $x_i R y$ for more than one $i$, again by
	uniqueness.
	%~
	Hence, $y$ must be incomparable to $n-1$ cycle members.
	%~
	However, for both $n=3$ and $n=4$ this implies that $y$ is
	incomparable to two adjacent cycle members, w.l.o.g.\ to $x_1$
	and $x_2$, contradicting incomparability-transitivity.
\end{enumerate}
The proof of the dual claims is similar.
\qed
}

\LEMMA{}{%
\LABEL{incTrans 12}\dpd{2}%
On set set $X$ of at least $5$ elements,
a left unique and right serial relation $R$ cannot be
incomparability-transitive (\QRef{212}),
and neither can a right unique and left serial relation (\QRef{200}).
%~
On the $4$-element set $X = \set{a,b,c,d}$, the relation
$R = \set{ \tpl{a,b}, \tpl{b,c}, \tpl{c,d}, \tpl{d,a} }$ is a
counter-example for the first claim.
}
\PROOF{
Assume for contradiction $X$ has at least $5$ elements,
and $R$ is a left
unique, right serial, and incomparability-transitive relation on $X$.

First, $xRx$ cannot hold for any $x$.
%~
Else, we had by Lem.~\REF{incTrans 2}\dpr{1} that
$x$ is comparable to every $y \in X \setminus \set{x}$.
%~
Since $yRx$ would imply the contradiction $x=y$ by uniqueness,
we even had $xRy$ for every $y \in X \setminus \set{x}$.
%~
By seriality, every such $y$ has an $R$-successor; by uniqueness, at
most one such successor can be $x$.
%~
Hence we can find $y_1, y_2 \in X \setminus \set{x}$ with $y_1 R y_2$.
%~
But this contradicts $x R y_2$ and uniqueness.

Second, by seriality, we find a chain
$x_1 R x_2 \land x_2 R x_3 \land \ldots$.
%~
Let $m$ be maximal such that $\:,1n{x_\i}$ are pairwise distinct;
by our first observation, we have $m \geq 2$, and $x_i \neq x_{i+1}$
for all chain members.
%~
For $1 \leq i \leq m$ and $2 \leq j \leq m$, we can't have $x_i R x_j$
when $i \neq j-1$, since else $x_i = x_{j-1}$ by uniqueness,
contradicting distinctness.
%~
Therefore, if $m \geq 5$,
we had $x_2$ incomparable to both $x_4$ and $x_5$,
contradicting $x_r R x_5$.
%~
In particular, $m$ can't be infinite.

We thus have $x_{m+1} = x_k$ for some $k \in \set{ \:,1n\i }$,
that is, $x_m R x_k$, which by uniqueness enforces $k=1$.
%~
That is, starting from an arbitrary member $x_1$, we always find a
cycle $x_1 R x_2 \land \ldots \land x_{m-1} R x_m \land x_m R x_1$
with $2 \leq m \leq 4$.

Since we have $\geq 5$ elements, we find another cycle
$y_1 R y_2 \land \ldots \land y_{n-1} R y_n \land y_n R y_1$ of some
length $n$.
%~
Then each $x_i$ is incomparable to each $y_j$, since $x_i R y_j$ would
imply $x_i = y_{j-1}$ or $x_i = y_m$, i.e.\ both
cycles would be identical; by symmetry, $y_j R x_i$ would imply the
same contradiction.
%~
But $y_1$ incomparable both to $x_1$ and to $x_2$
contradicts $x_1 R x_2$.
\qed
}

\LEMMA{}{%
\LABEL{incTrans 13}\dpd{1}%
If $R$ is incomparability-transitive, then it always
satisfies semi-order property 2 (\QRef{071}).
}
\PROOF{
Let $xRy \land yRz$ hold.
%~
If both $x,w$ and $w,y$ were incomparable, then so would be $x,y$,
contradicting $xRy$.
\qed
}

\LEMMA{}{%
\LABEL{incTrans 14}%
Let $R$ satisfy semi-order property 2.
%~
Then $R$ is incomparability-transitive if one of the following
sufficient conditions holds:
\begin{enumerate}
\item\LABEL{incTrans 14 1}\dpd{1}%
	$R$ is left quasi-reflexive (\QRef{180});
\item\LABEL{incTrans 14 2}\dpd{1}%
	$R$ is right quasi-reflexive (\QRef{180});
\item\LABEL{incTrans 14 3}\dpd{1}%
	$R$ is symmetric (\QRef{179});
\item\LABEL{incTrans 14 4}\dpd{1}%
	$R$ is transitive and dense (\QRef{234});
\item\LABEL{incTrans 14 5}\dpd{2}%
	$R$ is left Euclidean (\QRef{177});
\item\LABEL{incTrans 14 6}\dpd{2}%
	$R$ is right Euclidean (\QRef{178});
\item\LABEL{incTrans 14 7}\dpd{1}%
	$R$ is anti-transitive, quasi-transitive,
	and left and right serial (\QRef{265}).
\end{enumerate}
By Lem~\REF{incTrans 13}\dpr{1},
if any of the conditions~\REF{incTrans 14 1}
to~\REF{incTrans 14 7} holds, then $R$ satisfies
semi-order property 2 iff $R$ is incomparability-transitive.
%~
The latter doesn't hold without some extra conditions:
on the set $X=\set{a,b,c}$, the relation $R=\set{\tpl{a,c}}$
satisfies semi-order property 2, but isn't
incomparability-transitive.
}
\PROOF{
For cases~\REF{incTrans 14 1} to~\REF{incTrans 14 4},
assume for contradiction $aRb$ holds
and $c$ is incomparable both to $a$ and to $b$.
%~
In each of these cases,
we construct a chain $xRy \land yRz$ such that
$c$ is incomparable to all of $x,y,z$,
thus contradicting semi-order property 2.
\begin{enumerate}
\item%incTrans 14 1
	If $R$ is left quasi-reflexive, we have $aRa$.
	%~
	Choose $x,y,z$ to be $a,a,b$.
\item%incTrans 14 2
	If $R$ is right quasi-reflexive, we have $bRb$.
	%~
	Choose $x,y,z$ to be $a,b,b$.
\item%incTrans 14 3
	If $R$ is symmetric, we have $bRa$.
	%~
	Choose $x,y,z$ to be $a,b,a$.
\item%incTrans 14 4
	If $R$ is dense, we have $a R a' \land a' R b$.
	%~
	Choose $x,y,z$ to be $a,a',b$, we find that $c$ must be
	comparable to $a'$.
	%~
	However, $a' R c$ implies $aRc$, while $c R a'$ implies $cRb$,
	both by transitivity, and both contradicting our
	incomparability assumptions.
\item%incTrans 14 5
	Follows from~\REF{incTrans 14 1}\dpr{1},
	since each left Euclidean
	relation is left quasi-reflexive by Lem.~\REF{eucl 11}\dpr{1}.
\item%incTrans 14 6
	Follows similarly from~\REF{incTrans 14 2}\dpr{1}
	and Lem.~\REF{eucl 11}\dpr{1}.
\item%incTrans 14 7
	Assume for contradiction $R$ satisfies semi-order property 2
	and all properties from~\REF{incTrans 14 7},
	but isn't incomparability-transitive.
	\begin{enumerate}
	\item First, from the conjunction of anti-transitivity and
		quasi-transitivity we can draw some strong conclusions:
		%~
		Whenever $xRy \land yRz$ holds,
		then we must have $yRx \lor zRy$, and
		$x$ and $z$ must be incomparable.
		%~
		If neither $yRx$ nor $zRy$ held, then $xRz$ and
		its negation would follow by quasi-transitivity and
		anti-transitivity, respectively.
		%~
		$xRz$ is forbidden by anti-transitivity.
		%~
		If $zRx$ held, then $yRx$
		would imply $yRz \land zRx \land yRx$,
		while $zRy$ would imply $zRx \land xRy \land zRy$;
		both cases contradicting anti-transitivity.
	\item Second, since $R$ isn't incomparability-transitive,
		we have $aRb$ and
		$a$ as well as $b$ is incomparable to some $c$.
		%~
		By semi-order property 2, we can't have $bRa$.
		%~
		By left and right seriality applied to $a$ and $b$,
		we find $a' R a$
		and $b R b'$, respectively.
		%~
		From the first observation,
		we can conclude that $a R a'$ as well as
		$b' R b$ must hold, too,
		while $a'$ and $b$ must be incomparable, and
		so must be $a$ and $b'$.
	\item Third, by semi-order property 2,
		$b'$ can't be incomparable to $a'$, since it is to $a$,
		and $a' R a \land a R a'$ holds.
		%~
		Similarly, $c$ can't be incomparable to $a'$, and
		neither to $b'$,
		since it is to $a$, and to $b$, respectively.
		%~
		Moreover, we cannot have $a' R c \land c R a'$,
		since $b$ is
		incomparable to both $a'$ and $c$;
		similarly, we can't have $b' R c \land c R b'$.
		%~
		And we can't have $a' R c \land c R b'$, since this
		would imply
		incomparability of $a'$ and $b'$ by our first observation;
		for the same reason, we can't have $b' R c \land c R a'$.
	\item Altogether, two possibilities remain:
		\begin{enumerate}
		\item $a' R c \land b' R c$.

			Then $a' R b'$ would imply
			$a' R b' \land b' R c \land a' R c$,
			contradicting anti-transitivity;
			and $b' R a'$ would yield a symmetric
			contradiction.
		\item $c R a' \land c R b'$.

			Then $a' R b'$ would imply
			$c R a' \land a' R b' \land c R b'$,
			again contradicting anti-transitivity;
			similar for $b' R a'$.
		\end{enumerate}
	\end{enumerate}
\end{enumerate}
For the converse direction,
let $R$ be incomparability-transitive and
let $xRy \land yRz$ hold.
%~
If both $x,w$ and $w,y$ were incomparable, then so would be $x,y$,
contradicting $xRy$.
\qed
}

%\clearpage
\subsection{Euclideanness}
\LABEL{Euclideanness}

\LEMMA{}{%
\LABEL{eucl 1}\dpd{1}%
For symmetric relations, transitivity, right Euclideanness,
and left Euclideanness all coincide
(\QRef{098}, \QRef{099}, \QRef{117}).
%~
In particular, each equivalence relation is both left and right Euclidean.
}
\PROOF{
Let $R$ be symmetric.
%~
Transitivity implies right Euclideanness:
%~
Given $xRy$ and $xRz$,
we have $yRx \land xRz$ by symmetry,
hence $yRz$ by transitivity.
%~ 
The proof that right implies left Euclideanness and the proof that the
latter implies transitivity are similar.
\qed
}

\LEMMA{}{%
\LABEL{eucl 2}\dpd{2}%
A right Euclidean and left quasi-reflexive relation
is always symmetric, and hence transitive and left Euclidean
(\QRef{139}).
%~
Dually, a left Euclidean and right quasi-reflexive relation
is always symmetric and hence transitive and right Euclidean
(\QRef{138}).
%~
As a consequence,
a right and left Euclidean relation is symmetric (\QRef{097})
and hence transitive (\QRef{114}).
%~
A reflexive and right \emph{or} left Euclidean relation
is an equivalence, and both left (\QRef{132})
\emph{and} right (\QRef{131})
Euclidean.
%~
On the two-element set $X = \set{a, b}$,
the relation $R := \set{ \tpl{a,a} }$ is left and right quasi-reflexive,
left and right Euclidean, symmetric, transitive, but not reflexive,
hence no equivalence.
}
\PROOF{
\begin{itemize}
\item
	Let $R$ be right Euclidean and left quasi-reflexive.
	%~
	Then $R$ is also symmetric, since $xRy$ implies $xRx$
	by quasi-reflexivity,
	and both together imply $yRx$ by right Euclideanness.
	%~
	Hence, by Lem.~\REF{eucl 1}\dpr{1},
	$R$ is also transitive and left Euclidean.
	%~
	The proof for a left Euclidean $R$ is similar.
\item
	If $R$ is left Euclidean, then it is left quasi-reflexive
	by Lem.~\REF{eucl 11}\dpr{1}.
	%~
	Hence if $R$ is also right Euclidean, then
	it is symmetric and transitive as shown above.
\item
	If $R$ is reflexive and right Euclidean,
	then it is quasi-reflexive by Lem.~\REF{quasiRefl 2}\dpr{1},
	and hence symmetric, transitive and left Euclidean
	as shown above.
	%~
	The proof for a left Euclidean relation is similar.
\qed
\end{itemize}
}

\LEMMA{}{%
\LABEL{eucl 3}\dpd{3}%
The range of a right Euclidean relation is always a subset
of its domain.
%~
The restriction of a right Euclidean relation to its range
is always an equivalence.
%~
Similarly, the domain of a left Euclidean relation is a subset
of its range, and the restriction of a left Euclidean relation
to its domain is an equivalence.
%~
In particular, each left serial and right Euclidean relation
is an equivalence (\QRef{189}),
and so is each right serial and left Euclidean relation
(\QRef{201}).
}
\PROOF{
If $y$ is in the range of $R$, then $xRy \land xRy$ implies
$yRy$, for some suitable $x$.
%~
This also proves that $y$ is in the domain of $R$.
%~
By Lem.~\REF{eucl 2}\dpr{2}, $R$ is therefore an equivalence.

If $R$ is left serial, then every element is in the range of $R$.

The proofs for the dual claims are similar.
\qed
}

\LEMMA{}{%
\LABEL{eucl 4}\dpd{4}%
A relation $R$ is both left and right Euclidean, iff
the domain and the range set of $R$ agree, and $R$ is an
equivalence relation on that set (\QRef{097}, \QRef{114}).
}
\PROOF{
``$\Ra$'':
%~
follows by Lem.~\REF{eucl 3}\dpr{3}.

``$\La$'':
%~
Assume $aRb$ and $aRc$, then $a,b,c$ are
members of the domain and range of $R$, hence $bRc$ by symmetry
and transitivity.
%~
Left Euclideanness of $R$ follows similarly.
\qed
}

\LEMMA{}{%
\LABEL{eucl 5}\dpd{4}%
A right Euclidean relation is always
vacuously quasi-transitive (\QRef{076}),
and so is a left Euclidean relation (\QRef{075}).
}
\PROOF{
Let $R$ be right Euclidean.
%~
Let $xRy \land \lnot yRx \land yRz \land \lnot zRy$ hold.
%~
Observe that $y, z \in \ran(R)$.
%~
By Lem.~\REF{eucl 3}\dpr{3}, $R$ is symmetric on $\ran(R)$,
hence $yRz$ implies $zRy$, which is a contradiction.
%~
Hence, $R$ is vacuously quasi-transitive, since the assumptions
about $x,y,z$ can never be met.

A similar argument applies to left Euclidean relations,
exploiting that $x,y \in \dom(R)$.
\qed
}

\LEMMA{}{%
\LABEL{eucl 6}\dpd{4}%
A semi-connex right Euclidean relation is always transitive
(\QRef{155}),
and so is a semi-connex left Euclidean relation (\QRef{154}).
%~
On the set $X = \set{ a,b }$,
the relation $R = \set{ \tpl{a,a}, \tpl{a,b} }$
is semi-connex and left Euclidean, but not symmetric.
}
\PROOF{
Let $R$ be semi-connex and right Euclidean.
%~
Let $xRy \land yRz$ hold.
%~
Observe again that $y, z \in \ran(R)$.
%~
Since $R$ is semi-connex,
the following case distinction is exhaustive:
\begin{itemize}
\item
	$xRz$ holds.

	Then we are done immediately.
\item
	$zRx$ holds.

	Then also $x \in \ran(R)$;
	hence $xRz$, since $R$ is symmetric on its range
	by Lem.~\REF{eucl 3}\dpr{3}.
\item
	$x=z$.

	Then also $x \in \ran(R)$;
	hence $xRz$, since $R$ is reflexive on its range
	by Lem.~\REF{eucl 3}\dpr{3}.
\end{itemize}
Again, a similar argument applies to semi-connex and Euclidean
relations, using $x,y \in \dom(R)$.
\qed
}

\LEMMA{}{%
\LABEL{eucl 7}\dpd{4}%
If $X$ has at least 3 elements,
a semi-connex right Euclidean relation on $X$
is never anti-symmetric,
and neither is a semi-connex left Euclidean relation on $X$.
}
\PROOF{
Let $R$ be semi-connex and right Euclidean.
%~
By Lem.~\REF{conn 3}\dpr{1}, at most one element of $X$
is not in the range of $R$.
%~
Hence, by assumption,
two distinct elements $x,y \in \ran(R)$ exist.
%~
Since $R$ is semi-connex and $x \neq y$, we have $xRy$ or $yRx$.
%~
By Lem.~\REF{eucl 3}\dpr{3},
we obtain both $xRy$ and $yRx$.
%~
This contradicts the anti-symmetry requirement.
\qed
}

\definecolor{fEqcA}	{rgb}{0.99,0.00,0.00}
\definecolor{fEqcB}	{rgb}{0.00,0.99,0.00}
\definecolor{fEqcC}	{rgb}{0.00,0.00,0.99}
\definecolor{fEqcD}	{rgb}{0.99,0.99,0.00}
\definecolor{fEqcE}	{rgb}{0.00,0.99,0.99}

\definecolor{bEqcA}	{rgb}{0.99,0.50,0.50}
\definecolor{bEqcB}	{rgb}{0.50,0.99,0.50}
\definecolor{bEqcC}	{rgb}{0.50,0.50,0.99}
\definecolor{bEqcD}	{rgb}{0.99,0.99,0.50}
\definecolor{bEqcE}	{rgb}{0.50,0.99,0.99}

\definecolor{bDomXr}	{rgb}{0.85,0.85,0.85}
\definecolor{bDomXi}	{rgb}{0.97,0.97,0.97}

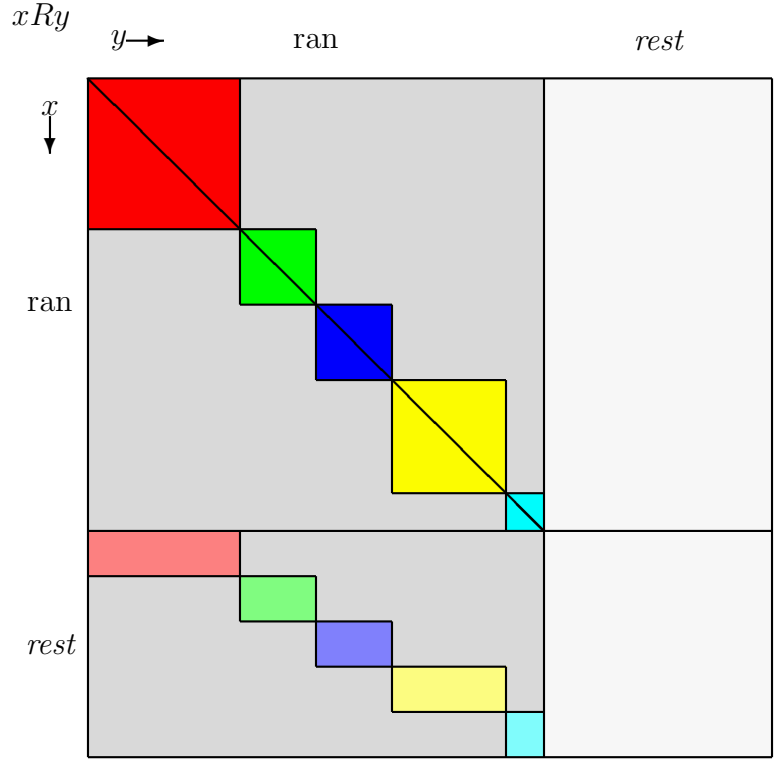
\begin{figure}
\begin{center}
\begin{picture}(100,100)
	%\put(0,0){\makebox(0,0){$+$}}
	%\put(100,100){\makebox(0,0){$+$}}
% dom,ran split
\textcolor{bDomXr}{\put(10,90){\makebox(0,0)[tl]{\rule{60mm}{90mm}}}}%
\textcolor{bDomXi}{\put(70,90){\makebox(0,0)[tl]{\rule{30mm}{90mm}}}}%
% equaivalence class backgrounds in X_r
\textcolor{fEqcA}{\put(10,90){\makebox(0,0)[tl]{\rule{20mm}{20mm}}}}%
\textcolor{fEqcB}{\put(30,70){\makebox(0,0)[tl]{\rule{10mm}{10mm}}}}%
\textcolor{fEqcC}{\put(40,60){\makebox(0,0)[tl]{\rule{10mm}{10mm}}}}%
\textcolor{fEqcD}{\put(50,50){\makebox(0,0)[tl]{\rule{15mm}{15mm}}}}%
\textcolor{fEqcE}{\put(65,35){\makebox(0,0)[tl]{\rule{05mm}{05mm}}}}%
% equaivalence class echo backgrounds in X_i
\textcolor{bEqcA}{\put(10,30){\makebox(0,0)[tl]{\rule{20mm}{06mm}}}}%
\textcolor{bEqcB}{\put(30,24){\makebox(0,0)[tl]{\rule{10mm}{06mm}}}}%
\textcolor{bEqcC}{\put(40,18){\makebox(0,0)[tl]{\rule{10mm}{06mm}}}}%
\textcolor{bEqcD}{\put(50,12){\makebox(0,0)[tl]{\rule{15mm}{06mm}}}}%
\textcolor{bEqcE}{\put(65,06){\makebox(0,0)[tl]{\rule{05mm}{06mm}}}}%
\thicklines%
% outer frame
\put(10,90){\line(1,0){90}}%
\put(10,0){\line(1,0){90}}%
\put(10,0){\line(0,1){90}}%
\put(100,0){\line(0,1){90}}%
% separation X_r,X_i
\put(70,0){\line(0,1){90}}%
\put(10,30){\line(1,0){90}}%
% reflexive diagonal in X_r
\put(10,90){\line(1,-1){60}}%
% equivalence class 1,2 in X_r
\put(10,70){\line(1,0){30}}%
\put(30,90){\line(0,-1){30}}%
% equivalence class 2,3 in X_r
\put(30,60){\line(1,0){20}}%
\put(40,70){\line(0,-1){20}}%
% equivalence class 3,4 in X_r
\put(40,50){\line(1,0){25}}%
\put(50,60){\line(0,-1){25}}%
% equivalence class 4,5 in X_r
\put(50,35){\line(1,0){20}}%
\put(65,50){\line(0,-1){20}}%
% equivalence class 1,2 echo in X_i
\put(10,24){\line(1,0){30}}
\put(30,30){\line(0,-1){12}}
% equivalence class 2,3 echo in X_i
\put(30,18){\line(1,0){20}}
\put(40,24){\line(0,-1){12}}
% equivalence class 3,4 echo in X_i
\put(40,12){\line(1,0){25}}
\put(50,18){\line(0,-1){12}}
% equivalence class 4,5 echo in X_i
\put(50,06){\line(1,0){20}}
\put(65,12){\line(0,-1){12}}
% captions X_r, X_i
\put(5,60){\makebox(0,0){$\ran$}}
\put(5,15){\makebox(0,0){\it rest}}
\put(40,95){\makebox(0,0){$\ran$}}
\put(85,95){\makebox(0,0){\it rest}}
% captions row,col = x,y
\put(0,100){\makebox(0,0)[tl]{$xRy$}}
\put(5,87){\makebox(0,0)[t]{$x$}}
\put(13,95){\makebox(0,0)[l]{$y$}}
\put(5,85){\vector(0,-1){5}}
\put(15,95){\vector(1,0){5}}
\end{picture}
\caption{Right Euclidean relation}
\LABEL{Right Euclidean relation}
\end{center}
\end{figure}

\begin{figure}
\begin{center}
\begin{picture}(100,100)
	%\put(0,0){\makebox(0,0){$+$}}
	%\put(100,100){\makebox(0,0){$+$}}
% dom,ran split
\textcolor{bDomXr}{\put(10,90){\makebox(0,0)[tl]{\rule{90mm}{60mm}}}}%
\textcolor{bDomXi}{\put(10,30){\makebox(0,0)[tl]{\rule{90mm}{30mm}}}}%
% equaivalence class backgrounds in X_r
\textcolor{fEqcA}{\put(10,90){\makebox(0,0)[tl]{\rule{20mm}{20mm}}}}%
\textcolor{fEqcB}{\put(30,70){\makebox(0,0)[tl]{\rule{10mm}{10mm}}}}%
\textcolor{fEqcC}{\put(40,60){\makebox(0,0)[tl]{\rule{10mm}{10mm}}}}%
\textcolor{fEqcD}{\put(50,50){\makebox(0,0)[tl]{\rule{15mm}{15mm}}}}%
\textcolor{fEqcE}{\put(65,35){\makebox(0,0)[tl]{\rule{05mm}{05mm}}}}%
% equaivalence class echo backgrounds in X_i
\textcolor{bEqcA}{\put(70,90){\makebox(0,0)[tl]{\rule{06mm}{20mm}}}}%
\textcolor{bEqcB}{\put(76,70){\makebox(0,0)[tl]{\rule{06mm}{10mm}}}}%
\textcolor{bEqcC}{\put(82,60){\makebox(0,0)[tl]{\rule{06mm}{10mm}}}}%
\textcolor{bEqcD}{\put(88,50){\makebox(0,0)[tl]{\rule{06mm}{15mm}}}}%
\textcolor{bEqcE}{\put(94,35){\makebox(0,0)[tl]{\rule{06mm}{05mm}}}}%
\thicklines%
% outer frame
\put(10,90){\line(1,0){90}}%
\put(10,0){\line(1,0){90}}%
\put(10,0){\line(0,1){90}}%
\put(100,0){\line(0,1){90}}%
% separation X_r,X_i
\put(70,0){\line(0,1){90}}%
\put(10,30){\line(1,0){90}}%
% reflexive diagonal in X_r
\put(10,90){\line(1,-1){60}}%
% equivalence class 1,2 in X_r
\put(10,70){\line(1,0){30}}%
\put(30,90){\line(0,-1){30}}%
% equivalence class 2,3 in X_r
\put(30,60){\line(1,0){20}}%
\put(40,70){\line(0,-1){20}}%
% equivalence class 3,4 in X_r
\put(40,50){\line(1,0){25}}%
\put(50,60){\line(0,-1){25}}%
% equivalence class 4,5 in X_r
\put(50,35){\line(1,0){20}}%
\put(65,50){\line(0,-1){20}}%
% equivalence class 1,2 echo in X_i
\put(70,70){\line(1,0){12}}
\put(76,90){\line(0,-1){30}}
% equivalence class 2,3 echo in X_i
\put(76,60){\line(1,0){12}}
\put(82,70){\line(0,-1){20}}
% equivalence class 3,4 echo in X_i
\put(82,50){\line(1,0){12}}
\put(88,60){\line(0,-1){25}}
% equivalence class 4,5 echo in X_i
\put(88,35){\line(1,0){12}}
\put(94,50){\line(0,-1){20}}
% captions X_r, X_i
\put(5,60){\makebox(0,0){$\dom$}}
\put(5,15){\makebox(0,0){\it rest}}
\put(40,95){\makebox(0,0){$\dom$}}
\put(85,95){\makebox(0,0){\it rest}}
% captions row,col = x,y
\put(0,100){\makebox(0,0)[tl]{$xRy$}}
\put(5,87){\makebox(0,0)[t]{$x$}}
\put(13,95){\makebox(0,0)[l]{$y$}}
\put(5,85){\vector(0,-1){5}}
\put(15,95){\vector(1,0){5}}
\end{picture}
\caption{Left Euclidean relation}
\LABEL{Left Euclidean relation}
\end{center}
\end{figure}

\LEMMA{}{%
\LABEL{eucl 8}\dpd{4}%
A relation $R$ on a set $X$ is right Euclidean iff
$R' := \restrict{R}{\ran(R)}$ is an equivalence
and
$\forall x \in X \.\setminus \ran(R) \; \exists y \in \ran(R). \;
xR \subseteq \eqc{y}{R'}$,
cf.\ Fig.~\REF{Right Euclidean relation}.
%~
Similarly, $R$ on $X$ is left Euclidean iff
$R' := \restrict{R}{\dom(R)}$ is an equivalence and
$\forall y \in X \.\setminus \dom(R) \; \exists x \in \dom(R). \;
Ry \subseteq \eqc{x}{R'}$,
cf.\ Fig.~\REF{Left Euclidean relation}.
}
\PROOF{
``$\Ra$'':
%~
By Lem.~\REF{eucl 3}\dpr{3}, $\restrict{R}{\ran(R)}$ is an equivalence.
%~
Let $x \in X \.\setminus \ran(R)$.
%~
If $x R y_1$ and $x R y_2$, then $y_1, y_2 \in \ran(R)$,
and $y_1 R y_2$ by right Euclideanness of $R$,
that is, $y_1, y_2$ belong to the same equivalence class
w.r.t.\ $R'$.

``$\La$'':
%~
Let $x, y, z \in X$ such that $xRy \land xRz$, we show $yRz$.
%~
Observe $y, z \in \ran(R)$.
%~
We distinguish two cases:
\begin{itemize}
\item If $x \in \ran(R)$,

	then $x R' y \land x R' z$,
	hence $y R' z$ by symmetry and transitivity of $R'$,
	hence $y R z$.
\item If $x \not\in \ran(R)$,

	then let $w \in \ran(R)$ with $xR \subseteq \eqc{w}{R'}$.
	%~
	We have $y, z \in \eqc{w}{R'}$ by assumption,
	hence $y R' z$,
	hence $y R z$.
\qed
\end{itemize}
}

Based on Lem.~\REF{eucl 8},
Fig.~\REF{Right Euclidean relation}
shows a schematized Right Euclidean relation.
%~
Deeply-colored squares indicate equivalence classes of
$\restrict{R}{\ran(R)}$, assuming $X$'s elements are arranged in such
a way that equivalent ones are adjacent.
%~
Pale-colored rectangles indicate possible relationships of elements in
$X \setminus \ran(R)$, again assuming them to be arranged in
convenient order.
%~
In these rectangles, relationships may, or may not, hold.
%~
A light grey color indicates that the element corresponding to the
line is unrelated to that corresponding to the column; in particular,
the lighter grey right rectangle indicates that
no element at all can be related to some in the set
${\it rest} := X \setminus \ran(R)$.
%~
The diagonal line indicates that $xRx$ holds iff $x \in \ran(R)$.

Figure~\REF{Left Euclidean relation}
shows a similar schema for a left Euclidean relation,

\LEMMA{}{%
\LABEL{eucl 9}\dpd{4}%
A left Euclidean and left unique relation is always
transitive, and so is a right Euclidean and right unique
relation (\QRef{115}, \QRef{116}).
%~
More particularly,
in both cases no chains
$xRy \land yRz$ with $x \neq y \land y \neq z$ can exist.
}
\PROOF{
Let $R$ be left Euclidean and left unique.
%~
Let $xRy$ and $yRz$ hold.
%~
By Lem.~\REF{eucl 3}\dpr{3},
$y \in \dom(R)$ implies $yRy$,
hence $x=y$,
hence $xRz$.
%~
The proof for right relations is similar.
\qed
}

\LEMMA{}{%
\LABEL{eucl 10}\dpd{1}%
A left Euclidean relation is left unique iff it is
anti-symmetric (\QRef{145}, \QRef{146}).
%~
Dually, a right Euclidean relation is right unique iff it is
anti-symmetric (\QRef{147}, \QRef{148}).
}
\PROOF{
Let $R$ be left Euclidean.

``$\Ra$'':
%~
If $xRy$ holds, then $xRx$ by Euclideanness.
%~
If also $yRx$ holds, we therefore have $x=y$ by uniqueness.

``$\La$'':
%~
If $x_1 R y$ and $x_2 R y$,
then both $x_1 R x_2$ and $x_2 R x_1$ follows by Euclideanness,
hence $x_1 = x_2$ by anti-symmetry.

The proof for a right Euclidean $R$ is similar.
\qed
}

\LEMMA{}{%
\LABEL{eucl 11}\dpd{1}%
Each left Euclidean relation is left quasi-reflexive.
%~
For left unique relations, the converse also holds.
%~
Dually, each right Euclidean relation is right quasi-reflexive,
and each right unique and right quasi-reflexive relation is right
Euclidean.
}
\PROOF{
\begin{itemize}
\item
	Let $R$ be left Euclidean.
	%~
	Then $xRy \land xRy$ implies $xRx$.
\item
	Let $R$ be left unique and left quasi-reflexive.
	%~
	If $yRx$ and $zRx$, then $y=z$ by left uniqueness,
	and $yRy$ by left quasi-reflexivity, hence $yRz$.
\item
	The proof for a right relations is similar.
\qed
\end{itemize}
}

\LEMMA{Sufficient for uniqueness implying Euclideanness}{%
\LABEL{eucl 11b}\dpd{2}%
A left unique relation $R$ is left Euclidean if one of the
following sufficient conditions holds:
\begin{enumerate}
\item\LABEL{eucl 11b 1}%
	$R$ is dense (\QRef{183});
\item\LABEL{eucl 11b 2}%
	$R$ is transitive and left serial (\QRef{236});
\item\LABEL{eucl 11b 4}%
	$R$ satisfies semi-order property 1 and
	is not anti-transitive (\QRef{217});
\item\LABEL{eucl 11b 5}%
	$R$ satisfies semi-order property 2 and
	is not anti-transitive (\QRef{224});
	or
\item\LABEL{eucl 11b 3}%
	$R$ is incomparability-transitive and not anti-transitive
	(\QRef{221}).
\end{enumerate}
Dually, a right unique relation is right Euclidean when
it additionally meets one of the following restrictions:
\begin{enumerate}
\item dense (\QRef{184});
\item transitive and right serial (\QRef{245});
\item semi-order property 1 and not anti-transitive (\QRef{218});
\item semi-order property 2 and not anti-transitive (\QRef{225}); or
\item incomparability-transitive and not anti-transitive
	(\QRef{222}).
\end{enumerate}
}
\PROOF{
Let $R$ be left unique.
%~
By Lem.~\REF{eucl 11}\dpr{1}, is is sufficient to show that $R$ is left
quasi-reflexive.
\begin{enumerate}
\item%eucl 11b 1
	Let $R$ additionally be dense.

	If $xRy$,
	then $xRw$ and $wRy$ for some $w$ by density,
	hence $x=w$ by left uniqueness,
	i.e.\ $xRx$.
\item%eucl 11b 2
	Let $R$ additionally be transitive and left serial.
	Let $xRy$,
	let $x' R x$ by seriality,
	then $x' R y$ by transitivity,
	hence $x'=x$ by uniqueness,
	i.e.\ $xRx$.
\item%eucl 11b 4
	Let $aRb$, $bRc$, and $aRc$ be a counter-example to
	anti-transitivity.
	%~
	Then $a=b$ by left uniqueness, and hence $bRb$.
	%~
	Let $xRy$ hold.
	%~
	Considering $b$ and $y$, three cases are possible:
	\begin{enumerate}
	\item $yRb$ holds.
		%~
		Then $y=a=b$ by uniqueness,
		hence $xRb$, hence $x=b$ by uniqueness,
		i.e.\ $xRx$.
	\item $bRy$ holds.
		%~
		Then $a=b=x$ by uniqueness, i.e.\ $xRx$.
	\item\LABEL{eucl 11b 4 c}
		$b$ and $y$ are incomparable.
		%~
		Then $xRc$ by semi-order property 1,
		hence $x=b$ by uniqueness,
		that is $xRx$.
	\end{enumerate}
	In each case, we have established left quasi-reflexivity.
\item%eucl 11b 5
	The proof proceeds as in~\REF{eucl 11b 4},
	except that in case~\REF{eucl 11b 4 c},
	applying semi-order property 2
	to $y$ and $bRb \land bRb$ obtains a contradiction.
\item%eucl 11b 3
	Follows from~\REF{eucl 11b 5}, since
	an incomparability-transitive relation satisfies semi-order 2
	by Lem.~\REF{incTrans 13}\dpr{1}.
\end{enumerate}
The proofs for a right unique $R$ are similar.
\qed
}

%\clearpage
\subsection{Density}
\LABEL{Density}

\LEMMA{Sufficient for density}{
\LABEL{den 1}%
A relation $R$ on a set $X$
is dense if one of the following sufficient conditions holds:
\begin{enumerate}
\item\LABEL{den 1 1}\dpd{1}%
	$R$ is reflexive (\QRef{085});
\item\LABEL{den 1 2}\dpd{4}%
	$R$ is co-reflexive (\QRef{081});
\item\LABEL{den 1 3}\dpd{1}%
	$R$ is left quasi-reflexive (\QRef{086});
\item\LABEL{den 1 4}\dpd{1}%
	$R$ is right quasi-reflexive (\QRef{086});
\item\LABEL{den 1 5}\dpd{2}%
	$R$ is left Euclidean (\QRef{082});
\item\LABEL{den 1 6}\dpd{2}%
	$R$ is right Euclidean (\QRef{083}).
\item\LABEL{den 1 7}\dpd{1}%
	$R$ is symmetric and satisfies semi-order property 1 (\QRef{187});
\item\LABEL{den 1 8}\dpd{3}%
	$R$ is connex (\QRef{084}).
\end{enumerate}
If $X$ has at least $3$ elements, then $R$ is also dense if
\begin{enumerate}
\setcounter{enumi}{8}
\item\LABEL{den 1 9}\dpd{1}%
	$R$ is symmetric and semi-connex (\QRef{188}).
\end{enumerate}
Even the conjunction of conditions~\REF{den 1 1}
to~\REF{den 1 9}
isn't a necessary condition for density.
}
\PROOF{
Let $x, z \in X$ be given such that $xRz$ holds.
\begin{enumerate}
\item%den 1 1
	If $R$ is reflexive, then $xRx \land xRz$ holds.
\item%den 1 2
	Shown in Lem.~\REF{corefl 5 8}\dpr{3}.
\item%den 1 3
	If $R$ is left quasi-reflexive, then $xRz$ implies $xRx$,
	hence $xRx \land xRz$ holds.
\item%den 1 4
	If $R$ is right quasi-reflexive, then $xRz$ implies $zRz$,
	hence $xRz \land zRz$ holds.
\item%den 1 5
	Follows from~\REF{den 1 3}\dpr{1},
	since $R$ is left quasi-reflexive
	by Lem.~\REF{eucl 11}\dpr{1}.
\item%den 1 6
	Follows from~\REF{den 1 4}\dpr{1},
	since $R$ is right quasi-reflexive
	by Lem.~\REF{eucl 11}\dpr{1}.
\item%den 1 7
	If $R$ is symmetric and satisfies semi-order property 1,
	then $xRz$ implies $zRx$, and both imply $xRx \lor zRz$ by
	semi-order property 1.
	%~
	Density follows similar to case~\REF{den 1 3} and~\REF{den 1 4}.
\item%den 1 8
	If $R$ is connex, then it is reflexive
	by Lem.~\REF{conn 1}\dpr{2},
	and hence dense by case~\REF{den 1 1}\dpr{1}.
\item%den 1 9
	Let $R$ be symmetric and semi-connex,
	let $x, z \in X$ be given such that $xRz$ holds.
	%~
	Let $y \in X$ be distinct from both $x$ and $z$,
	then $xRy \lor yRx$, and $yRz \lor zRy$ holds,
	since $R$ is semi-connex.
	%~
	By $R$'s symmetry, $xRy$ and $yRz$ holds, hence we are done.

	Note that
	on the two-element set $X = \set{a,b}$,
	the relation $R = \set{ \tpl{a,b}, \tpl{b,a} }$
	is symmetric and semi-connex, but not dense.
\end{enumerate}
On the set $X = \set{a,b,c}$,
the relation $R = \set{ \tpl{a,a}, \tpl{a,b}, \tpl{b,c}, \tpl{c,c} }$
is dense, but neither reflexive, nor co-reflexive,
nor left or right quasi-reflexive, nor left or right Euclidean,
nor symmetric, nor satisfying semi-order property 1, nor semi-connex.
\qed
}

\LEMMA{}{%
\LABEL{den 3}\dpd{1}%
A non-empty dense relation cannot be anti-transitive (\QRef{185}).
}
\PROOF{
Assume for contradiction that
$R$ is non-empty, dense, and anti-transitive.
%~
Due to the first property, $xRz$ holds for some $x, z$;
hence due to the second one, $xRy \land xRz$ holds for some
$y$;
hence due to the third one, $\lnot xRz$ holds, which is a
contradiction.
\qed
}

%\clearpage
\subsection{Connex and semi-connex relations}
\LABEL{Connex and semi-connex relations}

\LEMMA{}{%
\LABEL{conn 1}\dpd{2}%
A relation is connex iff it is semi-connex and reflexive
(\QRef{045}, \QRef{062}, \QRef{159}).
%~
If $X$ has at least 2 elements,
a relation $R$ on $X$ is connex iff it is semi-connex and
left and right quasi-reflexive (\QRef{050}, \QRef{062}, \QRef{160}).
%~
On a singleton set $X$, the empty relation is semi-connex and
quasi-reflexive, but not connex.
}
\PROOF{
\begin{itemize}
\item
	If $R$ is connex,
	the semi-connex property and the reflexivity follow trivially.
	%~
	The latter implies quasi-reflexivity by Lem.~\REF{quasiRefl 2}.
\item
	Conversely, let $R$ be semi-connex and reflexive.
	%~
	For $x \neq y$, the semi-connex property implies $xRy \lor yRx$.
	%~
	For $x=y$, reflexivity implies $xRy$.
\item
	Finally, let $R$ be semi-connex and quasi-reflexive.
	%~
	For $x \neq y$, the semi-connex property again
	implies $xRy \lor yRx$.
	%~
	For $x=y$, choose an arbitrary $z \neq x$,
	then $xRz$ or $zRx$ by the semi-connex property.
	%~
	Both cases imply $xRx$, i.e.\ $xRy$, by quasi-reflexivity.
\qed
\end{itemize}
}

\LEMMA{}{%
\LABEL{conn 2}\dpd{3}%
If a set $X$ has at least 4 elements,
then a semi-connex relation $R$ on can neither be
anti-transitive (\QRef{061}),
nor left (\QRef{059}) nor right (\QRef{060}) unique.
%~
The same applies in particular to a connex relation on $X$
(\QRef{028}, \QRef{026}, \QRef{027}).
}
\PROOF{
First, assume $R$ is both semi-connex and anti-transitive.
%~
Consider the directed graph corresponding to $R$, with its
vertices being the elements of $X$, and its edges being the
pairs related by $R$.

Consider three arbitrary distinct vertices.
%~
By the semi-connex property, each pair of them must be
connected by an edge.
%~
By anti-transitivity, (*) none of them may be the source
of more than one edge.
%~
Hence, the three edges must be oriented in such a way that
they for a directed cycle.

Let $w,x,y,z$ be four distinct elements,
%~
W.l.o.g.\ assume the subgraph for $x,y,z$ is oriented a
directed cycle corresponding to $xRy \land yRz \land zRx$.
%~
The subgraph for $w,x,y$ must be oriented as a directed
cycle, too; therefore $wRx \land xRy \land yRw$ must hold.
%~
But then, the subgraph for $w,x,z$ is not oriented as a cycle,
since $wRx \land zRx$.
%~
This contradicts the cycle-property shown above.

If $R$ is semi-connex and right unique, the latter property
implies (*) that no vertex may be the source of two edges,
and the proof is similar.

If $R$ is semi-connex and left unique, no vertex may be the
target of two edges, leading again to a similar proof.

Each connex relation is semi-connex
by Lem.~\REF{conn 1}\dpr{2}.
\qed
}

\LEMMA{}{%
\LABEL{conn 3}\dpd{1}%
If $R$ is a semi-connex relation on $X$,
then the set $X \setminus \ran(R)$ has at most one element;
the same applies to $X \setminus \dom(R)$.
}
\PROOF{
Let $x, y \in X \setminus \ran(R)$.
%~
Since $R$ is semi-connex, $xRy$ or $yRx$ or $x=y$ must hold.
%~
The first two possibilities are ruled out by assumption, so
the third one must apply, i.e.\ $x$ and $y$ can't be distinct.
%~
A similar argument applies to $\dom(R)$.
\qed
}

\LEMMA{Sufficient for connex implying universality}{%
\LABEL{conn 4}%
A connex relation $R$ on a set $X$
needs to be universal if one of the following
sufficient conditions is satisfied:
\begin{enumerate}
\item\LABEL{conn 4 1}\dpd{1}%
	$R$ is symmetric (\QRef{113});
\item\LABEL{conn 4 2}\dpd{3}%
	$R$ is left Euclidean (\QRef{111});
\item\LABEL{conn 4 3}\dpd{3}%
	$R$ is right Euclidean (\QRef{112}).
\end{enumerate}
}
\PROOF{
\begin{enumerate}
\item%conn 4 1
	Let $x,y \in X$.
	%~
	Then $xRy \lor yRx$ by the connex property.
	%~
	Hence $xRy \land yRx$ by symmetry.
\item%conn 4 2
	$R$ is reflexive by Lem.~\REF{conn 1}\dpr{2},
	hence symmetric by Lem.~\REF{eucl 2}\dpr{2}.
	%~
	So we are done using case~\REF{conn 4 1}\dpr{1}.
\item%conn 4 3
	The proof is dual to case~\REF{conn 4 2}.
\qed
\end{enumerate}
}

%\clearpage
\subsection{Seriality}
\LABEL{Seriality}

\LEMMA{}{%
\LABEL{ser 1}\dpd{3}%
A reflexive relation is always both right (\QRef{094}) and left
(\QRef{090}) serial.
%~
In particular, a connex relation has both properties (\QRef{093},
\QRef{089}).
}
\PROOF{
Given $x$, choose $x$ both as an $R$-successor and an $R$-predecessor.
%~
By Lem.~\REF{conn 1}\dpr{2},
in particular each connex relation is reflexive.
\qed
}

\LEMMA{}{%
\LABEL{ser 3}%
A right serial relation $R$ is reflexive if one of the following
sufficient conditions is met:
\begin{enumerate}
\item\LABEL{ser 3 1}\dpd{1}%
	$R$ is co-reflexive (\QRef{207});
\item\LABEL{ser 3 3}\dpd{1}%
	$R$ is left quasi-reflexive (\QRef{209}); or
\item\LABEL{ser 3 2}\dpd{2}%
	$R$ is left Euclidean (\QRef{208}).
\end{enumerate}
Dually, a left serial relation $R$ is reflexive if one of the
following sufficient conditions is met:
\begin{enumerate}
\item
	$R$ is co-reflexive (\QRef{195});
\item
	$R$ is right quasi-reflexive (\QRef{197}); or
\item
	$R$ is right Euclidean (\QRef{196}).
\end{enumerate}
}
\PROOF{
Let $R$ be a right serial relation on the set $X$ and let $x \in X$.
%~
Then $xRy$ holds for some $y$.
%~
We have to show $xRx$.
\begin{enumerate}
\item%ser 3 1
	If $R$ is also co-reflexive, then $y=x$.
\item%ser 3 3
	If $R$ is also left quasi-reflexive,
	then $xRy$ implies $xRx$.
\item%ser 3 2
	Follows from~\REF{ser 3 3}\dpr{1},
	since $R$ is left quasi-reflexive by Lem.~\REF{eucl 11}\dpr{1}.
\end{enumerate}
If $R$ is left serial, the proofs are similar.
\qed
}

\LEMMA{Sufficient for semi-connex implying seriality}{%
\LABEL{ser 4}%
A semi-connex relation $R$ on a set $X$ with at least $2$ elements
is right serial if one of the following sufficient conditions is met:
\begin{enumerate}
\item\LABEL{ser 4 1}\dpd{1}%
	$R$ is right Euclidean (\QRef{210}); or
\item\LABEL{ser 4 2}\dpd{1}%
	$R$ is symmetric (\QRef{211}).
\end{enumerate}
Dually, a semi-connex relation
$R$ on $X$ is left serial if one of the following sufficient
conditions is met:
\begin{enumerate}
\item
	$R$ is left Euclidean (\QRef{198}); or
\item
	$R$ is symmetric (\QRef{199}).
\end{enumerate}
On a singleton set $X$, the empty relation is trivially
symmetric, semi-connex, and left and right Euclidean,
but neither right nor left serial.
}
\PROOF{
Let $R$ on $X$ be semi-connex.
%~
Given $x \in X$,
let $y \neq x$ be some other member of $X$;
then $xRy \lor yRx$ holds.
%~
In the former case, we are done immediately;
so we only need to consider the case $yRx$:
\begin{enumerate}
\item%ser 4 1
	If $R$ is also right Euclidean,
	$yRx \land yRx$ implies $xRx$.
\item%ser 4 2
	If $R$ is also symmetric,
	$yRx$ implies $xRy$.
\end{enumerate}
Left seriality properties follow similarly.
\qed
}

\LEMMA{}{%
\LABEL{ser 5}\dpd{1}%
On a nonempty, but finite domain $X$,
an irreflexive and transitive relation
cannot be right serial; neither can it be left serial.
}
\PROOF{
In this proof,
we write ``$<$'' instead of ``$R$'' for readability.
%~
Use induction on $n$ to show the existence of arbitrarily long
chains
$\:<1n{x_\i}$
for all $n \in \N$.
%~
Since $<$ is transitive, we have $x_i < x_j$
for all $1 \leq i < j \leq n$.
%~
Since $<$ is irreflexive, this implies $x_i \neq x_j$ for all
these $i, j$.
%~
For each $n$ larger than the finite cardinality of $X$, this
is a contradiction.
%~
The proof for left seriality is similar.
\qed
}

%\clearpage
\subsection{Uniqueness}
\LABEL{Uniqueness}

\LEMMA{}{%
\LABEL{uniq 1}\dpd{1}%
A left unique and transitive relation is always anti-symmetric
(\QRef{151}),
and so is a right unique and transitive relation (\QRef{152}).
}
\PROOF{
If both $xRy$ and $yRx$, then $xRx$ by transitivity,
hence $x=y$ by left or right uniqueness.
\qed
}

\LEMMA{}{%
\LABEL{uniq 2}\dpd{1}%
On a finite set $X$, a relation is
both right unique and left serial iff
it is both left unique and right serial
(\QRef{190}, \QRef{202}, \QRef{213}, \QRef{214}).
%~
On the set of natural numbers, the relation $y=x/\!\!/2 \land x>0$
is right unique (i.e.\ a partial function) and left serial
(i.e.\ surjective), but neither left unique (injective) nor
right serial (total), where ``$/\!\!/$'' denotes truncating
integer division (see Fig.~\REF{Counter-example in Lem.uniq 2});
the converse relation is a counter-example
for the opposite direction.
}
\PROOF{
``$\Ra$'':
%~
Let $R$ be right unique and left serial,
i.e.\ a partial function and surjective.
%~
Then $\ran(R) = X$.
%~
Since for each $y \in \ran(R)$ we have an $x \in \dom(R)$
such that $xRy$,
and since no two $y$ share an $x$, we have that $\dom(R)$
has no less elements than $\ran(R)$.
%~
Since $X$ is finite, this implies $\dom(R) = X$, i.e.\ $R$
is right serial, i.e.\ $R$ is a surjective total function.
%~
From set theory, we know that $R$ then also must be injective,
i.e.\ left unique.

``$\La$'':
%~
Apply the ``$\Ra$'' proof to the converse relation $R^{-1}$.
\qed
}

\LEMMA{}{%
\LABEL{uniq 3}\dpd{2}%
On a finite set $X$ with odd cardinality, a
left unique, symmetric, and left serial relation cannot
be irreflexive, and hence not anti-transitive (\QRef{235}).
%~
On the set $X = \set{a,b,c,d,e,f}$, the relation
$R = \set{ \tpl{a,b}, \tpl{b,a}, \tpl{c,d}, \tpl{d,c},
\tpl{e,f}, \tpl{f,e} }$
satisfies all properties simultaneously.
}
\PROOF{
Since $R$ is left unique and symmetric, each $x \in X$ can be
comparable to at most one element:
$xRy \lor yRx$ and $xRz \lor zRx$ implies
$xRy \land yRx \land xRz \land zRx$,
and in turn $y=z$.
%~
Due to irreflexivity and seriality, each $x$ must be
comparable to at least one element different from $x$.
%~
From both conditions together we obtain that each $x$ is
related to exactly one different element $x'$.
%~
This is impossible if $X$ is finite and of odd cardinality.
%~
By Lem.~\REF{antitrans 1}\dpr{1}, each anti-transitive relation is
irreflexive.
\qed
}

\begin{figure}
\begin{center}
\begin{picture}(115,25)
	%\put(0,0){\makebox(0,0){$+$}}
	%\put(115,25){\makebox(0,0){$+$}}
\put(3.000,5.000){\makebox(0.000,0.000)[r]{$y$}}
\put(6.000,5.000){\line(1,0){102.000}}
\put(109.000,5.000){\makebox(0.000,0.000)[l]{\ldots}}
\put(6.000,4.000){\line(0,1){2.000}}
\put(26.000,4.000){\line(0,1){2.000}}
\put(46.000,4.000){\line(0,1){2.000}}
\put(66.000,4.000){\line(0,1){2.000}}
\put(86.000,4.000){\line(0,1){2.000}}
\put(106.000,4.000){\line(0,1){2.000}}
\put(6.000,3.000){\makebox(0.000,0.000)[t]{$0$}}
\put(26.000,3.000){\makebox(0.000,0.000)[t]{$1$}}
\put(46.000,3.000){\makebox(0.000,0.000)[t]{$2$}}
\put(66.000,3.000){\makebox(0.000,0.000)[t]{$3$}}
\put(86.000,3.000){\makebox(0.000,0.000)[t]{$4$}}
\put(106.000,3.000){\makebox(0.000,0.000)[t]{$5$}}
\put(3.000,19.000){\makebox(0.000,0.000)[r]{$x$}}
\put(6.000,19.000){\line(1,0){102.000}}
\put(109.000,19.000){\makebox(0.000,0.000)[l]{\ldots}}
\put(6.000,18.000){\line(0,1){2.000}}
\put(16.000,18.000){\line(0,1){2.000}}
\put(26.000,18.000){\line(0,1){2.000}}
\put(36.000,18.000){\line(0,1){2.000}}
\put(46.000,18.000){\line(0,1){2.000}}
\put(56.000,18.000){\line(0,1){2.000}}
\put(66.000,18.000){\line(0,1){2.000}}
\put(76.000,18.000){\line(0,1){2.000}}
\put(86.000,18.000){\line(0,1){2.000}}
\put(96.000,18.000){\line(0,1){2.000}}
\put(106.000,18.000){\line(0,1){2.000}}
\put(6.000,21.000){\makebox(0.000,0.000)[b]{$0$}}
\put(16.000,21.000){\makebox(0.000,0.000)[b]{$1$}}
\put(26.000,21.000){\makebox(0.000,0.000)[b]{$2$}}
\put(36.000,21.000){\makebox(0.000,0.000)[b]{$3$}}
\put(46.000,21.000){\makebox(0.000,0.000)[b]{$4$}}
\put(56.000,21.000){\makebox(0.000,0.000)[b]{$5$}}
\put(66.000,21.000){\makebox(0.000,0.000)[b]{$6$}}
\put(76.000,21.000){\makebox(0.000,0.000)[b]{$7$}}
\put(86.000,21.000){\makebox(0.000,0.000)[b]{$8$}}
\put(96.000,21.000){\makebox(0.000,0.000)[b]{$9$}}
\put(106.000,21.000){\makebox(0.000,0.000)[b]{$10$}}
\put(16.000,17.000){\vector(-1,-1){10.000}}
\put(26.000,17.000){\vector(0,-1){10.000}}
\put(36.000,17.000){\vector(-1,-1){10.000}}
\put(46.000,17.000){\vector(0,-1){10.000}}
\put(56.000,17.000){\vector(-1,-1){10.000}}
\put(66.000,17.000){\vector(0,-1){10.000}}
\put(76.000,17.000){\vector(-1,-1){10.000}}
\put(86.000,17.000){\vector(0,-1){10.000}}
\put(96.000,17.000){\vector(-1,-1){10.000}}
\put(106.000,17.000){\vector(0,-1){10.000}}
\end{picture}
\caption{Counter-example in Lem.~\REF{uniq 2}}
\LABEL{Counter-example in Lem.uniq 2}
\end{center}
\end{figure}

%\clearpage
\subsection{Semi-order property 1}
\LABEL{Semi-order property 1}

\LEMMA{Sufficient conditions for semi-order property 1}{%
\LABEL{semiOrd1 3a}%
A relation $R$ on a set $X$ satisfies semi-order property 1 if one
of the following sufficient conditions holds:
\begin{enumerate}
\item\LABEL{semiOrd1 3a 6}\dpd{1}%
	$R$ is transitive and incomparability-transitive (\QRef{166});
\item\LABEL{semiOrd1 3a 1}\dpd{4}%
	$R$ is semi-connex and transitive (\QRef{158});
\item\LABEL{semiOrd1 3a 2}\dpd{5}%
	$R$ is semi-connex and left Euclidean (\QRef{156});
\item\LABEL{semiOrd1 3a 3}\dpd{5}%
	$R$ is semi-connex and right Euclidean (\QRef{157});
\item\LABEL{semiOrd1 3a 4}\dpd{1}%
	$R$ is dense and transitive
	and satisfies semi-order property 2 (\QRef{233}); or
\item\LABEL{semiOrd1 3a 5}\dpd{5}%
	$R$ is quasi-reflexive and incomparability-transitive (\QRef{168}).
\end{enumerate}
}
\PROOF{
\begin{enumerate}
\item%semiOrd1 3a 6
	Let $wRx$ and $yRz$ hold for $x,y$ incomparable.
	%~
	Then $w,y$ must be comparable by incomparability-transitivity.
	%~
	However, $yRw$ would imply the contradiction $yRx$ by transitivity.
	%~
	Hence, we have $wRy$, which implies $wRz$ by transitivity.
\item%semiOrd1 3a 1
	Follows from~\REF{semiOrd1 3a 6}\dpr{1}, since $R$ is
	incomparability-transitive by Lem.~\REF{incTrans 1}\dpr{3}.
\item%semiOrd1 3a 2
	Follows from~\REF{semiOrd1 3a 1}\dpr{4},
	since semi-connex and left or right Euclidean relation is
	transitive by Lem.~\REF{eucl 6}\dpr{4}.
\item%semiOrd1 3a 3
	By~\REF{semiOrd1 3a 1}\dpr{4}
	and Lem.~\REF{eucl 6}\dpr{4} again.
\item%semiOrd1 3a 4
	Let $wRx \land yRz$ hold for some incomparable $x,y$.
	%~
	By density, we have $y R y' \land y' R z$ for some $y'$.
	%~
	By semi-order property 2, $x$ is comparable to one of $y,y',z$,
	that is, one of $y',z$.
	%~
	Now $y' R x$ would imply $yRx$ by transitivity;
	from $zRx$, we would get the same contradiction.
	%~
	But $x R y'$ implies $wRz$ by transitivity,
	and $x R z$ implies the same.
	%~
	Altogether, $wRz$ must hold.
\item%semiOrd1 3a 5
	By Lem.~\REF{incTrans 6}\dpr{2}, $R$ is connex or empty.
	%~
	Hence, $R$ satisfies semi-order property 1,
	in the first case by Lem.~\REF{semiOrd1 1a}\dpr{4},
	in the second case trivially (cf.\ \QRef{033}).
\qed
\end{enumerate}
}

\LEMMA{Sufficient for semi-order property 1 implying transitivity}{%
\LABEL{semiOrd1 12}%
An relation $R$ satisfying semi-order property 1 is
transitive if one of the following sufficient conditions is met:
\begin{enumerate}
\item\LABEL{semiOrd1 12 3}\dpd{1}%
	$R$ is irreflexive;
\item\LABEL{semiOrd1 12 1}\dpd{1}%
	$R$ is left unique (\QRef{120});
\item\LABEL{semiOrd1 12 2}\dpd{1}%
	$R$ is right unique (\QRef{121});
\item\LABEL{semiOrd1 12 5}\dpd{2}%
	$R$ is asymmetric (\QRef{123});
	or
\item\LABEL{semiOrd1 12 4}\dpd{2}%
	$R$ is anti-transitive (\QRef{122}).
\end{enumerate}
}
On the set $X = \set{ 0,1,2,3 }$,
delimiting examples of transitive relations satisfying
semi-order property 1 are the following:
\begin{itemize}
\item the universal relation is
	neither irreflexive, nor left or right unique;
\item the relation defined by $xRy :\Lra x=0$ is
	left unique, but not irreflexive;
\item the relation defined by $xRy :\Lra y=0$ is
	right unique, but not irreflexive;
\item the relation $x<y$ is irreflexive, but neither left nor right unique.
\end{itemize}
\PROOF{
\begin{enumerate}
\item%semiOrd1 12 3
	Let $xRy$ and $yRz$ hold.
	%~
	Since $y,y$ is incomparable due to irreflexivity,
	we obtain $xRz$ by semi-order property 1.
\item%semiOrd1 12 1
	If $R$ is left unique and satisfies semi-order property 1,
	then $wRx$, $x,y$ incomparable, and $yRz$ implies
	$wRz$ by semi-order property 1,
	hence $y=w$ by uniqueness, hence the contradiction $yRx$.
	%~
	Therefore, the precondition of semi-order property 1 can never be
	satisfied.
	%~
	Hence, if $aRb$ and $bRc$ holds,
	then $b$ must be comparable to itself, i.e.\ $bRb$ must hold.
	%~
	By uniqueness, we get $a=b$, that is, $aRc$ holds trivially.
\item%semiOrd1 12 2
	The proof is dual to that of case~\REF{semiOrd1 12 1}.
\item%semiOrd1 12 5
	Follows from~\REF{semiOrd1 12 3}\dpr{1},
	since $R$ is irreflexive
	by Lem.~\REFF{asym 1b}{1}\dpr{1}.
\item%semiOrd1 12 4
	Follows from~\REF{semiOrd1 12 3}\dpr{1},
	since $R$ is irreflexive
	by Lem.~\REF{antitrans 1}\dpr{1}.
	%~
	Note that in this case, $R$ is vacuously transitive, that is,
	$xRy \land yRz$ never holds, cf.\ Lem.~\REF{antitrans 4}.
\qed
\end{enumerate}
}

\LEMMA{}{%
\LABEL{semiOrd1 4}%
Let $R$ be
symmetric and satisfy semi-order property 1.
%~
Then:
\begin{enumerate}
\item\LABEL{semiOrd1 4 1}\dpd{1}%
	if $R$ is left quasi-reflexive relation,
	then $R$ is left Euclidean (\QRef{219}); and
\item\LABEL{semiOrd1 4 2}\dpd{1}%
	if $R$ is right quasi-reflexive relation,
	then $R$ is right Euclidean.
\end{enumerate}
}
\PROOF{
We show~\REF{semiOrd1 4 1}; the proof of~\REF{semiOrd1 4 2} is
similar.
%~
If $yRx$ and $zRx$ holds, then $yRy$ and $zRz$ by left quasi-reflexivity.
%~
Hence incomparability of $y$ and $z$ would contradict semi-order
property 1.
%~
By symmetry, therefore both $yRz$ and $zRy$.
\qed
}

\LEMMA{Sufficient for semi-order property 1
	implying incomparability-transitivity}{%
\LABEL{semiOrd1 6}%
A relation $R$ that satisfies semi-order property 1 is
incomparability-transitive if it satisfies one of the following
sufficient conditions:
\begin{enumerate}
\item\LABEL{semiOrd1 6 1}\dpd{1}%
	$R$ is left unique and left serial (\QRef{241});
\item\LABEL{semiOrd1 6 4}\dpd{1}%
	$R$ is right unique and right serial (\QRef{249}).
\end{enumerate}
}
\PROOF{
We show case~\REF{semiOrd1 6 4}; the other condition is proven dually.
%~
Assume for contradiction $aRb$ holds and $c$ is incomparable both to
$a$ and to $b$.
%~
By right seriality, we obtain $c R c'$.
%~
By semi-order property 1, we have $a R c'$.
%~
By right uniqueness, we get $b=c'$,
contradicting $c R c' \land \lnot cRb$.
\qed
}

\LEMMA{Incompatibilities of semi-order property 1}{%
\LABEL{semiOrd1 7}%
If $X$ has at least $2$ elements,
no relation $R$ can satisfy semi-order property 1
and one of the following conditions:
\begin{enumerate}
\item\LABEL{semiOrd1 7 1}\dpd{1}%
	$R$ is left unique and right serial (\QRef{206});
\item\LABEL{semiOrd1 7 2}\dpd{1}%
	$R$ is right unique and left serial (\QRef{194}).
\end{enumerate}
On a singleton set, the universal relation is a counter-example.
}
\PROOF{
For a $2$ element set $X$, all $4$ relations on $X$ are easily checked;
we assume in the following the $X$ has $\geq 3$ elements.
%~
We show case~\REF{semiOrd1 7 1}; the other condition is proven dually.
%~
For an arbitrary $w$, find $wRx$ by right seriality.
%~
Choose $y \in X \setminus \set{w}$.
%~
Then $yRx$ would imply the contradiction $y=w$ by left uniqueness.
%~
Moreover, $x,y$ incomparable would imply $w R y'$ by semi-order property 1,
where $y R y'$ is obtained by right seriality;
hence $w=y$ by left uniqueness, contradicting $wRx$.

So for arbitrary $w$,
and $x$ an $R$-successor of $w$,
we have that $\forall y \in X \setminus \set{w} . \; xRy$ must hold.
%~
Now let $w_1,w_2,y$ be pairwise distinct,
obtain $w_i R x_i$ by right seriality,
and $x_i R y$ by the above argument.
%~
Then $x_1 = x_2$, hence $w_1 = w_2$, both by left uniqueness;
this contradicts our assumption.
\qed
}

\LEMMA{}{%
\LABEL{semiOrd1 1a}\dpd{4}%
The following conditions are equivalent:
\begin{enumerate}
\item\LABEL{semiOrd1 1a 1}%
	$R$ is connex;
\item\LABEL{semiOrd1 1a 2}%
	$R$ is reflexive and satisfies semi-order property 1
	(\QRef{035}, \QRef{136});
\item\LABEL{semiOrd1 1a 3}%
	$R$ is reflexive and satisfies semi-order property 2
	(\QRef{069}, \QRef{175});
\item\LABEL{semiOrd1 1a 4}%
	$R$ is reflexive and semi-connex
	(\QRef{045}, \QRef{062}, \QRef{159}).
\end{enumerate}
}
\PROOF{
\begin{itemize}
\item \REF{semiOrd1 1a 1} $\Ra$ \REF{semiOrd1 1a 2}:
	%~
	If $R$ is connex, no $x,y$ are incomparable;
	by Lem.~\REF{conn 1}\dpr{2} $R$ is reflexive.
\item\REF{semiOrd1 1a 2} $\Ra$ \REF{semiOrd1 1a 1}:
	%~
	If $x,y$ were incomparable, applying semi-order property 1 to
	$xRx$ and $yRy$ would yield the contradiction $xRy$.
\item \REF{semiOrd1 1a 1} $\Ra$ \REF{semiOrd1 1a 3}:
	%~
	If $R$ is connex, then it is reflexive by
	Lem.~\REF{conn 1}\dpr{2}
	and satisfies semi-order property 2
	by Lem.~\REF{incTrans 1}\dpr{3}.
\item \REF{semiOrd1 1a 3} $\Ra$ \REF{semiOrd1 1a 1}:
	%~
	Given $w$ and $x$, apply semi-order property 2 to
	$w$ and $xRx \land xRx$.
\item \REF{semiOrd1 1a 1} $\Lra$ \REF{semiOrd1 1a 4}:
	%~
	Shown in Lem.~\REF{conn 1}\dpr{2}.
\qed
\end{itemize}
}

%\clearpage
\subsection{Semi-order property 2}
\LABEL{Semi-order property 2}

\LEMMA{}{%
\LABEL{semiOrd2 2}\dpd{1}%
For a non-empty relation satisfying semi-order property 2,
reflexivity and quasi-reflexivity are equivalent.
}
\PROOF{
If $R$ is nonempty and quasi-reflexive,
$xRy$, and hence $xRx$ holds for some $x,y$.
%~
Applying semi-order property 2 to $xRx \land xRx$ and an
arbitrary $w$ yields $wRx \lor xRw$, i.e.\ $wRw$ by
quasi-reflexivity.
%~
Hence $R$ is reflexive.
%~
The converse direction is trivial.
\qed
}

\LEMMA{}{%
\LABEL{semiOrd2 4}\dpd{4}%
A left Euclidean relation is transitive if it
satisfies semi-order property 2 (\QRef{170}).
%~
The same applies to a right Euclidean relation (\QRef{171}).
%~
On the set $X = \set{ a,b }$,
the relation $R = \set{ \tpl{a,a}, \tpl{a,b} }$
satisfies semi-order property 2 and
if left Euclidean and transitive, but not reflexive.
}
\PROOF{
Let $R$ be left Euclidean and satisfying semi-order property 2.
%~
Let $xRy$ and $yRz$ hold; we will show $xRz$.
%~
By Lem.~\REF{eucl 3}\dpr{3}, $x \in \dom(R)$ implies
$xRx$.
%~
Applying semi-order property 2 to $xRx$, $xRx$, and $z$,
we obtain $xRz$ or $zRx$.
%~
In the former case, we are done.
%~
In the latter case, we have $z \in \dom(R)$,
hence $zRx$ implies $xRz$ by Lem.~\REF{eucl 3}\dpr{3}.

The proof for right Euclideanness is similar.
\qed
}

\LEMMA{}{%
\LABEL{semiOrd2 5}\dpd{4}%
If $R$ is right Euclidean and satisfies semi-order property 2,
then it is incomparability-transitive (\QRef{178}),
and satisfies semi-order property 1 (\QRef{174}).
%~
The same applies if $R$ is left rather than right Euclidean
(\QRef{177}, \QRef{173}).
}
\PROOF{
Let $R$ be right Euclidean and satisfy semi-order property 2.
\begin{itemize}
\item
	Let $wRx$ and $yRz$ hold, and $x,y$ be incomparable;
	we show $wRz$.
	%~
	By Lem.~\REF{eucl 3}\dpr{3}, we have $xRx$;
	applying semi-order property 2 to $wRx \land xRx$ and $y$
	yields $yRw$ or $wRy$.
	%, since $x,y$ are incomparable.
	%~
	In the former case,
	$yRz$ and right Euclideanness yields $wRz$.
	%~
	In the latter case,
	$wRx$ and right Euclideanness yields $xRy$,
	contradicting $x,y$'s incomparability.
\item
	Incomparability-transitivity has been shown in
	Lem.~\REFF{incTrans 14}{6}\dpr{2}.
\item
	The proofs for left Euclideanness are similar.
\qed
\end{itemize}
}

\LEMMA{}{%
\LABEL{semiOrd2 6}\dpd{1}%
On a set $X$ of at least $5$ elements,
each relation $R$ that is left and right unique
and satisfies semi-order property 2,
needs to be transitive or left serial (\QRef{254}).
%~
By duality, such a relation also needs to be transitive or right
serial.
%~
On the $4$-element set $X = \set{a,b,c,d}$,
the relation $R = \set{ \tpl{a,b}, \tpl{b,c}, \tpl{c,d} }$
satisfies all antecedent properties, but non of the conclusion
properties.
}
\PROOF{
Let $X$ have $\geq 5$ elements,
let $R$ be a a left and right unique relation on $X$ satisfying
semi-order property 2.
%~
Assume for contradiction $aRb$ and $bRc$, but not $aRc$,
and $\lnot xRd$ for all $x \in X$.
%~
By semi-order property 2, $d$ needs to be related to at least one of
$a,b,c$, that is, $dRa$ or $dRb$ or $dRc$ holds.
%~
By left uniqueness, this implies
$dRa \land aRb \land bRc$
or $d=a \land aRb \land bRc$
or $d=b$;
the last case is impossible due to $\lnot aRd$.

In both possible cases, we have a chain $d R x_1 \land x_1 R x_2$.
%~
Now choose two distinct $y,z \in X \setminus \set{d,x_1,x_2}$.
%~
Then by semi-order property 2,
$y$ must be comparable to one of $d,x_1,x_2$.
%~
Due to left uniqueness,
we cannot have $y R x_1 \lor y R x_2$,
due to right uniqueness,
we cannot have $d R y \lor x_1 R y$.
%~
By definition of $d$, we cannot have $yRd$,
hence $x_2 R y$ must hold.
%~
However, the same arguments apply to $y$ as well, so $x_2 R y$ must
hold, too.
%~
By right uniqueness, this implies the contradiction $y=z$.
\qed
}

\LEMMA{Necessary for uniqueness and semi-order property 2}{%
\LABEL{semiOrd2 7b}%
Let $R$ on $X$ be left unique and satisfy semi-order property 2.
%~
Then $R$ is necessarily
\begin{enumerate}
\item\LABEL{semiOrd2 7b 2}\dpd{3}%
	left Euclidean or anti-transitive (\QRef{224});
\item\LABEL{semiOrd2 7b 2b}\dpd{4}%
	left quasi-reflexive or anti-transitive;
\item\LABEL{semiOrd2 7b 3}\dpd{1}%
	asymmetric or vacuously quasi-transitive (\QRef{230}); and
\item\LABEL{semiOrd2 7b 4}\dpd{1}%
	asymmetric or left serial (\QRef{242}).
\end{enumerate}
Dually, let $R$ on $X$ be right unique and satisfy semi-order property 2.
%~
Then $R$ is necessarily
\begin{enumerate}
\item right Euclidean or anti-transitive (\QRef{225});
\item right quasi-reflexive or anti-transitive;
\item asymmetric or vacuously quasi-transitive (\QRef{231}); and
\item asymmetric or right serial (\QRef{250}).
\end{enumerate}
}
\PROOF{
\begin{enumerate}
\item%semiOrd2 7b 2
	Shown in Lem.~\REFF{eucl 11b}{5}\dpr{2}.
\item%semiOrd2 7b 2b
	Follows from~\REF{semiOrd2 7b 2}\dpr{3},
	since left quasi-reflexivity is equivalent to left
	Euclideanness for left unique relations
	by Lem.~\REF{eucl 11}\dpr{1}.
\item%semiOrd2 7b 3
	Assume for contradiction $aRb$, but $bRa$,
	and $cRd \land \lnot dRc \land dRe \land \lnot eRd$,
	but $\lnot cRe \lor eRc$.
	%~
	By semi-order property 2, $e$ must be comparable to one of $a,b$;
	w.l.o.g.\ to $a$.
	%~
	If $aRe$, then $a=d$ and hence $c=b$, both by left uniqueness;
	this implies the contradiction $dRc$.
	%~
	If $eRa$, then $b=e$ and hence $a=d$ by left uniqueness,
	implying the contradiction $dRe$.
	%~
	Note that we didn't need the conclusion, $\lnot cRe \lor eRc$,
	of negated quasi-transitivity.
\item%semiOrd2 7b 4
	Assume for contradiction $aRb$ and $bRa$ holds,
	and $cRy$ doesn't hold for any $y$.
	%~
	By semi-order property 2, $c$ must be related to one of $a,b$.
	%~
	Since $cRa$ is impossible, we have w.l.o.g.\ that $aRc$ holds.
	%~
	By uniqueness, $c=b$, hence we have the contradiction $cRa$.
\end{enumerate}
The proof of the dual claims is similar.
\qed
}

\LEMMA{}{%
\LABEL{semiOrd2 8}\dpd{1}%
If $X$ has at least $3$ elements,
every left and right unique relation $R$ on $X$
that satisfies semi-order property 2
is asymmetric (\QRef{226}).
}
\PROOF{
Assume for contradiction $aRb$ and $bRa$.
%~
Let $a \neq c \neq b$.
%~
Then by semi-order property 2,
$c$ must be comparable to one of $a,b$; we assume w.l.o.g., to $a$.
%~
If $aRc$, then $b=c$ by right uniqueness,
if $cRa$ then $b=c$ by left uniqueness;
both contradicts our assumptions.
\qed
}

\LEMMA{}{%
\LABEL{semiOrd2 10}\dpd{1}%
Every left or right quasi-reflexive relation
satisfying semi-order property 2
vacuously also satisfies semi-order property1 (\QRef{176}).
}
\PROOF{
We show that the antecedent of Def.~\REFF{def}{SemiOrd1} is never
satisfied.
%~
If $wRx$ and $yRz$ then $yRy$ for a left quasi-reflexive $R$;
applying semi-order property 2 to $x$ and $yRy \land yRy$ yields that
$x,y$ can't be incomparable.
%~
For a right quasi-reflexive $R$, apply semi-order property 2 to $y$
and $xRx$.
\qed
}

\clearpage
\section{Examples}
\LABEL{Examples}

In this section, we collect those computed properties that gave
rise to single examples, rather than to general laws.
%~
Laws about the empty (Exm.~\REF{empty}) and the universal
(Exm.~\REF{univ}) relation were reported properly by our algorithm.
%~
Each other example arose since the reported law suggestion was true on a
$5$-element universe, but turned out to be false on a larger one.
%~
In most cases, we just relied on the algorithm for both claims;
in Exm.~\REF{exm 186} and~\REF{exm 252,255,257} we gave formal proofs.

\EXAMPLE{Empty relation}{%
\LABEL{empty}%
The empty relation $R = \set{}$ on a set $X$ has the following
properties:
\begin{enumerate}
\item\LABEL{empty 1}\dpd{2}%
	$R$ is co-reflexive (\QRef{002}),
	hence by Lem.~\REFF{corefl 5}{10}\dpr{1}
	also quasi-reflexive (\QRef{047});
\item\LABEL{empty 2}\dpd{1}%
	$R$ is left (\QRef{004}) and right (\QRef{007}) Euclidean;
\item\LABEL{empty 3}\dpd{1}%
	$R$ is left (\QRef{010}) and right (\QRef{013}) unique;
\item\LABEL{empty 4}\dpd{1}%
	$R$ is symmetric (\QRef{016});
\item\LABEL{empty 5}\dpd{1}%
	$R$ is anti-transitive (\QRef{019});
\item\LABEL{empty 6}\dpd{2}%
	$R$ is asymmetric (\QRef{021}),
	hence by Lem.~\REF{asym 1b}\dpr{1}
	also anti-symmetric (\QRef{052})
	and irreflexive (\QRef{036});
\item\LABEL{empty 7}\dpd{1}%
	$R$ is transitive (\QRef{030});
\item\LABEL{empty 7b}\dpd{3}%
	$R$ is quasi-transitive (\QRef{072})
	by~\REF{empty 7}\dpr{1} and Lem.~\REF{quasiTrans 4}\dpr{2};
\item\LABEL{empty 8}\dpd{1}%
	$R$ satisfies semi-order properties~1 (\QRef{033})
	and~2 (\QRef{067});
\item\LABEL{empty 8b}\dpd{1}%
	$R$ is incomparability-transitive (\QRef{063}).
\item\LABEL{empty 9}\dpd{1}%
	$R$ is dense (\QRef{079}).
\end{enumerate}
If $X$ is not empty, then
\begin{enumerate}
\setcounter{enumi}{11}%
\item\LABEL{empty 10}\dpd{1}%
	$R$ is not the universal relation (\QRef{001});
\item\LABEL{empty 11}\dpd{4}%
	$R$ is neither left (\QRef{087}) nor right (\QRef{091})
	serial,
	hence by Lem.~\REF{ser 1}\dpr{3} not reflexive (\QRef{041}).
\end{enumerate}
If $X$ has at least $2$ elements, then
\begin{enumerate}
\setcounter{enumi}{13}%
\item\LABEL{empty 12}\dpd{3}%
	$R$ is not semi-connex (\QRef{056}),
	and hence by Lem.~\REF{conn 1}\dpr{2} not connex (\QRef{023}).
\qed
\end{enumerate}
}

\EXAMPLE{Universal relation}{%
\LABEL{univ}%
On a set $X$,
the universal relation $R = X \times X$ has the following properties:
\begin{enumerate}
\item\LABEL{univ 1}\dpd{1}%
	$R$ is left (\QRef{005})
	and right (\QRef{008}) Euclidean;
\item\LABEL{univ 2}\dpd{1}%
	$R$ is symmetric (\QRef{017});
\item\LABEL{univ 3}\dpd{3}%
	$R$ is connex (\QRef{024})
	hence by Lem.~\REF{conn 1}\dpr{2}
	also semi-connex (\QRef{057});
\item\LABEL{univ 4}\dpd{3}%
	$R$ is transitive (\QRef{031}),
	hence by Lem.~\REF{quasiTrans 4}\dpr{2}
	also quasi-transitive (\QRef{073});
\item\LABEL{univ 5}\dpd{1}%
	$R$ satisfies semi-order properties~1 (\QRef{034})
	and~2 (\QRef{068});
\item\LABEL{univ 6}\dpd{4}%
	$R$ is reflexive (\QRef{042})
	hence by Lem.~\REF{quasiRefl 2}\dpr{1}
	also quasi-reflexive (\QRef{048}),
	by Lem.~\REF{ser 1}\dpr{3}
	also left (\QRef{088}),
	and right (\QRef{092}) serial,
	and by Lem.~\REFF{den 1}{1}\dpr{1}
	also dense (\QRef{080});
\item\LABEL{univ 7}\dpd{1}%
	$R$ is incomparability-transitive (\QRef{064}).
\end{enumerate}
If $X$ is not empty, then
\begin{enumerate}
\setcounter{enumi}{7}%
\item\LABEL{univ 8}\dpd{1}%
	$R$ is not the empty relation (\QRef{001});
\item\LABEL{univ 9}\dpd{2}%
	$R$ is not irreflexive (\QRef{037}),
	hence by Lem.~\REF{antitrans 1}\dpr{1}
	not anti-transitive (\QRef{020}),
	and by Lem.~\REFF{asym 1b}{1}\dpr{1}
	not asymmetric (\QRef{022}).
\end{enumerate}
If $X$ has at least $2$ elements, then
\begin{enumerate}
\setcounter{enumi}{9}%
\item\LABEL{univ 10}\dpd{1}%
	$R$ is not co-reflexive (\QRef{003});
\item\LABEL{univ 11}\dpd{1}%
	$R$ is not anti-symmetric (\QRef{053});
\item\LABEL{univ 12}\dpd{1}%
	$R$ is neither left (\QRef{011})
	nor right (\QRef{014}) unique.
\qed
\end{enumerate}
}

% qmc law suggestions that turned out to be false,
% except for small universe sets:
% 186 227 243 251 258 259 261 262 263 264 266 267 268 269 270 271 272 274

\begin{figure}
\begin{center}
\includegraphics[scale=0.5]{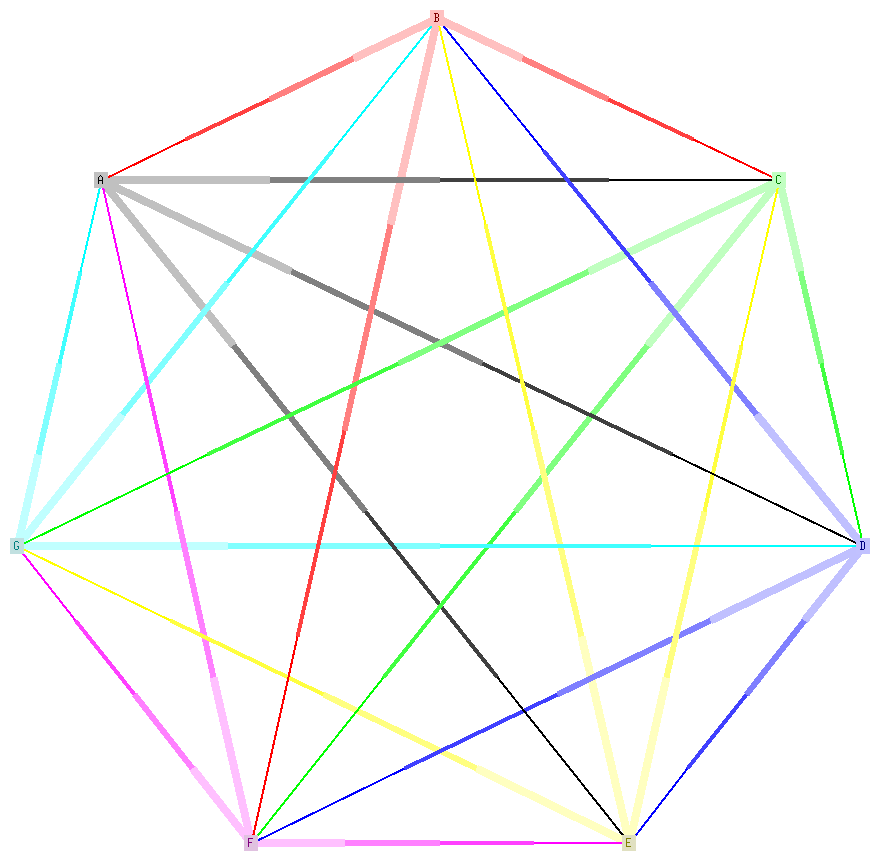}
\caption{Relation graph in Exm.~\REF{exm 186}}
%%\caption{Relation graph in Exm.~\REF{exm 186} (\QRef{186})}
\LABEL{asymmetric dense 7-element example}
\end{center}
\end{figure}

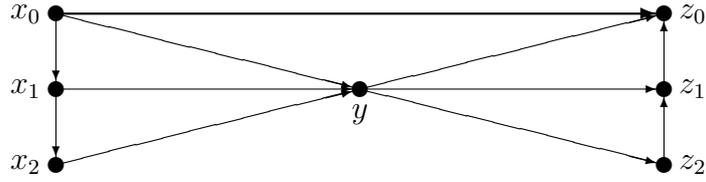
\begin{figure}
\begin{center}
\begin{picture}(100,30)
	%\put(0,0){\makebox(0,0){$+$}}
	%\put(100,30){\makebox(0,0){$+$}}
\put(10,5){\circle*{2}}
\put(10,15){\circle*{2}}
\put(10,25){\circle*{2}}
\put(50,15){\circle*{2}}
\put(90,5){\circle*{2}}
\put(90,15){\circle*{2}}
\put(90,25){\circle*{2}}
\put(50,13){\makebox(0,0)[t]{$y$}}
\put(8,5){\makebox(0,0)[r]{$x_2$}}
\put(8,15){\makebox(0,0)[r]{$x_1$}}
\put(8,25){\makebox(0,0)[r]{$x_0$}}
\put(92,5){\makebox(0,0)[l]{$z_2$}}
\put(92,15){\makebox(0,0)[l]{$z_1$}}
\put(92,25){\makebox(0,0)[l]{$z_0$}}
\put(10,15){\vector(0,-1){9}}
\put(10,25){\vector(0,-1){9}}
\put(10,25){\vector(4,-1){39}}
\put(10,5){\vector(4,1){39}}
\put(90,5){\vector(0,1){9}}
\put(90,15){\vector(0,1){9}}
\put(50,15){\vector(4,1){39}}
\put(50,15){\vector(4,-1){39}}
\put(10,15){\vector(1,0){39}}
\put(50,15){\vector(1,0){39}}
\thicklines
\put(10,25){\vector(1,0){79}}
\end{picture}
\caption{Asymmetry and density requires at least 7 distinct elements}
\LABEL{asymmetric dense minimum count}
\end{center}
\end{figure}

\EXAMPLE{}{%
\LABEL{exm 186}\dpd{2}%
A relation $R$ on a set $X$ cannot non-empty, dense, and
asymmetric if $X$ has no more than 6 elements (\QRef{186}).
%~
On a set $X$ of at least 7 elements,
these properties can be satisfied simultaneously.
%~
On the infinite set of all rational numbers,
they are satisfied e.g.\ by the usual strict ordering.

Unsatisfiability on small sets has been machine-checked.
%~
Subsequently, the following argument was found,
cf.\
Fig.~\REF{asymmetric dense minimum count}:
%~
If $x_0 R z_0$ for some elements $x_0, z_0$,
then by density
$x_0 R y \land y R z_0$,
hence
$x_0 R x_1 \land x_1 R y$
hence
$x_1 R x_2 \land x_2 R y$;
and dually
$y R z_1 \land z_1 R z_0$
and
$y R z_2 \land z_2 R z_1$.
%~
Since $R$ is asymmetric, it is also irreflexive
by Lem.~\REFF{asym 1b}{1}\dpr{1};
therefore
$x_0 \neq x_1 \neq x_2$
and
$z_0 \neq z_1 \neq z_2$
and
$x_i \neq y \neq z_j$ for all $i,j \in \set{0,1,2}$.
%~
The plain asymmetry of $R$ implies
$x_i \neq z_j$ for all $i,j \in \set{0,1,2}$
and $x_0 \neq x_2$
and $z_0 \neq z_2$.
%~
Therefore, all seven elements
$x_0, x_1, x_2, y, z_0, z_1, z_2$ are pairwise distinct;
i.e.\ $X$ has at least 7 elements.

The graph shown in
Fig.~\REF{asymmetric dense 7-element example}
shows a non-empty asymmetric dense relation on a 7-element
set; its properties have been machine-checked, too.
%~
An arrow from $x$ (light blunt end) to $y$ (dark peaked end)
indicates $xRy$.
%~
Each vertex has three outgoing edges, all sharing its color;
their opposite vertices are always connected to each other
by a directed cycle.
%~
Dually, each vertex has three incoming edges, with their
opposite vertices again connected by a directed cycle.
%~
For example, vertex $A$ points to the cycle $C \ra D \ra E \ra C$,
and is pointed to by the cycle $B \ra F \ra G \ra B$;
in terms of
Fig.~\REF{asymmetric dense minimum count},
$x_0, x_1, x_2, y, z_0, z_1, z_2$
corresponds to $B, F, G, A, E, D, C$, respectively.
%~
However,
Fig.~\REF{asymmetric dense minimum count}
matches
Fig.~\REF{asymmetric dense 7-element example}
in a multitude of other ways.
\qed
}

\begin{figure}
\begin{center}
\includegraphics[scale=0.5]{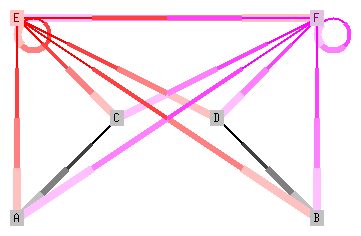}
\caption{Relation graph in Exm.~\REF{exm 227}}
%%\caption{Relation graph in Exm.~\REF{exm 227} (\QRef{227})}
\LABEL{Counter-example 227}
\end{center}
\end{figure}
%	Checking 227: Trans, SemiOrd2 ==> AntiSym || SemiOrd1
%	x\y      0 1 2 3 4 5
%		 A B C D E F
%
%	 0       . . X . X X
%	 1       . . . X X X
%	 2       . . . . X X
%	 3       . . . . X X
%	 4       . . . . X X
%	 5       . . . . X X
%
%	tr 1, s2 1, as 0, s1 0

\EXAMPLE{}{%
\LABEL{exm 227}\dpd{1}%
On a set $X$ of $6$ elements, a relation $R$ can be transitive but not
anti-symmetric, and satisfy semi-order property 2, but not 1
(\QRef{227}).
%~
An example is shown in Fig.~\REF{Counter-example 227}.
%~
Anti-symmetry is violated by $eRf \land fRe$.
%~
Semi-order property 1 is violated by $aRc$, $bRd$, $b,c$ incomparable,
but not $aRd$.
%~
On a set of $\leq 5$ elements, no relation with these properties exists.
\qed
}

\begin{figure}
\begin{center}
\includegraphics[scale=0.5]{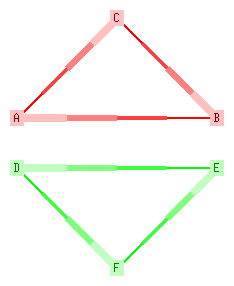}
\caption{Relation graph in Exm.~\REF{exm 243,251}}
%%\caption{Relation graph in Exm.~\REF{exm 243,251} (\QRef{243}, \QRef{251})}
\LABEL{Counter-example 251}
\end{center}
\end{figure}

%	Checking 243: RgUnique, ASym, LfSerial ==> SemiOrd2
%	x\y      0 1 2 3 4 5 6 7
%
%	 0       . X . . . . . .
%	 1       . . X . . . . .
%	 2       X . . . . . . .
%	 3       . . . . X . . .
%	 4       . . . . . X . .
%	 5       . . . . . . X .
%	 6       . . . . . . . X
%	 7       . . . X . . . .
%
%	ru 1, as 1, ls 1, s2 0

%	Checking 251: LfSerial, RgSerial, ASym ==> SemiOrd2
%	x\y      0 1 2 3 4 5
%
%	 0       . X . . . .
%	 1       . . X . . .
%	 2       X . . . . .
%	 3       . . . . X .
%	 4       . . . . . X
%	 5       . . . X . .
%
%	ls 1, rs 1, as 1, s2 0

\EXAMPLE{}{%
\LABEL{exm 243,251}\dpd{1}%
On a set $X$ of $6$ elements, a relation $R$ can be right unique, left
serial, and asymmetric but not
satisfying semi-order property 2 (\QRef{243}).
%~
An example is shown in Fig.~\REF{Counter-example 251}.
%~
Semi-order property 2 is violated since $aRb \land bRc$, but $d$ isn't
comparable to any of $a,b,c$.
%~
Since the relation shown there is also right serial, it is a
counter-example for \QRef{251}, too.
%~
On a set of $\leq 5$ elements, no relation with either property combination
exists.
\qed
}

\EXAMPLE{}{%
\LABEL{exm 252,255,257}\dpd{1}%
On a set $X$ of $\leq 5$ elements,
a relation $R$ must be anti-transitive, or anti-symmetric, or
quasi-transitive, if one of the following holds:
\begin{enumerate}
\item $R$ is left and right unique (\QRef{252}), or
\item $R$ is left unique, and not left serial (\QRef{255}).
\item $R$ is right unique, and not right serial (\QRef{257}).
\end{enumerate}
On the $6$-element set $X = \set{a,b,c,d,e,f}$,
the relations
\begin{enumerate}
\item $R_1 =
	\set{ \tpl{a,a}, \tpl{b,c}, \tpl{c,b}, \tpl{d,e}, \tpl{e,f} }$,
\item $R_2 =
	\set{ \tpl{a,a}, \tpl{b,c}, \tpl{c,b}, \tpl{c,d}, \tpl{d,e} }$,
	and
\item $R_3 =
	\set{ \tpl{a,a}, \tpl{b,c}, \tpl{c,b}, \tpl{d,c}, \tpl{e,d} }$
\end{enumerate}
satisfy all respective properties simultaneously.
%~
$R_1$ also satisfies property~2 and~3.

We show case 1.\ and 2.; the proof for case 3.\ is similar to that of
case 2.
%~
First, for both cases,
we investigate the properties of a relation $R$
that is left unique, but neither anti-transitive, nor anti-symmetric,
nor quasi-transitive.

As a counter-example to anti-transitivity,
let $a R a_2$, $a_2 R a_3$, but $a R a_3$ hold.
%~
Then $a_2 = a$ by left uniqueness;
that is, the counter-example collapses to $aRa$, we don't use $a R a_3$.

As a counter-example to anti-symmetry,
let $bRc$, but $cRb$ hold, for $b \neq c$.
%~
Then $a \neq b$, since else $cRa$ and $aRa$ would imply $a=b=c$ by
left uniqueness; similar $a \neq c$.

As a counter-example to quasi-transitivity,
let $dRe \land \lnot eRd \land eRf \land \lnot fRe$,
but $\lnot dRf \lor fRd$ hold.
%~
Then $e \neq a$, since else $e=d=a$ by left uniqueness, implying the
contradiction $\lnot aRa$.
%~
Similarly, $f \neq a$.
%~
Moreover, $e \neq b$, since else $d = c$ by left uniqueness, implying
the contradiction $eRd$.
%~
Similarly, $f \neq b$.
%~
By a symmetry argument, we also have $e \neq c \neq f$.
%~
We have $d \neq e$,
since the contrary would imply the contradiction $eRd$;
similarly $e \neq f$.
%~
And we have $d \neq f$, since else we had the contradiction $eRd$.
%~
To sum up, we have shown that the set $\set{a,b,c,e,f}$ has a
cardinality of $5$, and $\set{d,e,f}$ has $3$ elements,
but we couldn't rule out $d \in \set{a,b,c}$.
%~
Second, we distinguish the cases 1.\ and 2.\ in order to use the
additional properties of $R$ in each case:
\begin{enumerate}
\item Let $R$ be additionally right unique.
	%~
	Then $d \neq a$ since else $e=a$ by right uniqueness,
	contradicting $\lnot eRd$.
	%~
	Moreover, $d \neq b$, since else $e=c$ by right uniqueness,
	hence $eRd$, contradicting again $\lnot eRd$.
	%~
	Altogether, $a,b,c,d,e,f$ are pairwise distinct.
\item
	As a counter-example to left seriality, let $g$ be such that
	$xRg$ is false for every $x$.
	%~
	Then $g \not\in \set{a,b,c,e,f}$, since all members of that set
	have an $R$-predecessor.
	%~
	Hence $a,b,c,e,f,g$ are pairwise distinct.
\end{enumerate}
In both cases, $X$ must have at least $6$ elements.
\qed
}

\begin{figure}
\begin{center}
\includegraphics[scale=0.5]{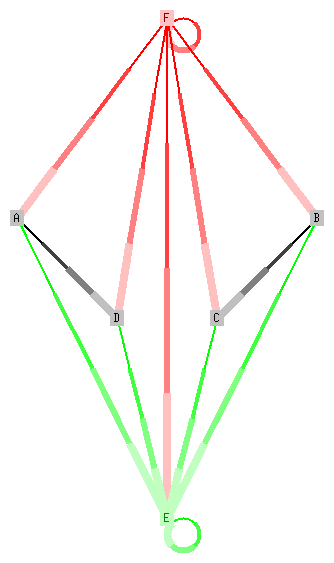}
\caption{Relation graph in Exm.~\REF{exm 258}}
%%\caption{Relation graph in Exm.~\REF{exm 258} (\QRef{258})}
\LABEL{Counter-example 258}
\end{center}
\end{figure}

%	Checking 258: Trans, SemiOrd2, LfSerial, RgSerial ==> SemiOrd1
%	x\y      0 1 2 3 4 5
%		 A B C D E F
%
%	 0       . . . . . X
%	 1       . . . . . X
%	 2       . X . . . X
%	 3       X . . . . X
%	 4       X X X X X X
%	 5       . . . . . X
%
%	tr 1, s2 1, ls 1, rs 1, s1 0

\EXAMPLE{}{%
\LABEL{exm 258}\dpd{1}%
On a set $X$ of $6$ elements, a relation $R$ can be transitive, left
and right serial,
and satisfy semi-order property 2, but not 1 (\QRef{258}).
%~
An example is shown in Fig.~\REF{Counter-example 258}.
%~
Semi-order property 1 is violated, since $dRa$, $cRb$, $a,c$
incomparable, but not $dRb$.
%~
On a set of $\leq 5$ elements, no relation with these properties exists.
\qed
}

\begin{figure}
\begin{center}
\includegraphics[scale=0.5]{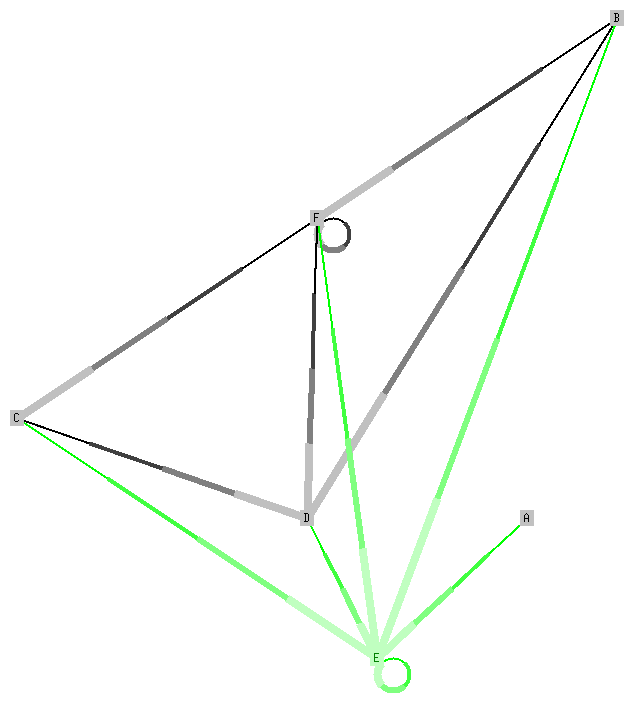}
\caption{Relation graph in Exm.~\REF{exm 259}}
%%\caption{Relation graph in Exm.~\REF{exm 259} (\QRef{259})}
\LABEL{Counter-example 259}
\end{center}
\end{figure}

%	Checking 259: SemiOrd1, AntiSym, LfSerial ==>
%	Trans || IncTrans || Dense
%	x\y      0 1 2 3 4 5
%
%	 0       . . . . . .
%	 1       . . . . . .
%	 2       . . . . . X
%	 3       . X X . . X
%	 4       X X X X X X
%	 5       . X . . . X

\EXAMPLE{}{%
\LABEL{exm 259}\dpd{1}%
On a set $X$ of $6$ elements, a relation $R$ can satisfy semi-order
property 1 and be anti-symmetric and
left serial, but neither transitive nor incomparability-transitive nor
dense (\QRef{259}).
%~
An example is shown in Fig.~\REF{Counter-example 259}.
%~
It is not transitive, since $cRf \land fRb$ but not $cRb$;
it is not incomparability-transitive, since $a,f$ and $a,b$ are
incomparable, while $fRb$;
it is not dense, since $eRa$ has no intermediate element.
%~
By inverting the arrow directions, a counter example for the dual
\QRef{260} is obtained.
%~
On a set of $\leq 5$ elements, no relation with either property combination
exists.
\qed
}

\begin{figure}
\begin{center}
\includegraphics[scale=0.5]{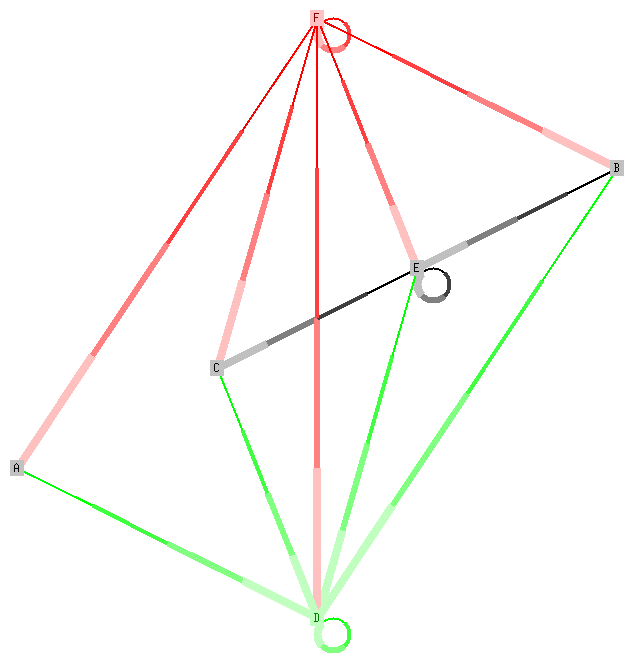}
\caption{Relation graph in Exm.~\REF{exm 261}}
%%\caption{Relation graph in Exm.~\REF{exm 261} (\QRef{261})}
\LABEL{Counter-example 261}
\end{center}
\end{figure}

%	Checking 261: SemiOrd1, AntiSym, LfSerial, RgSerial ==>
%	Trans || IncTrans
%	x\y      0 1 2 3 4 5
%
%	 0       . . . . . X
%	 1       . . . . . X
%	 2       . . . . X X
%	 3       X X X X X X
%	 4       . X . . X X
%	 5       . . . . . X

\EXAMPLE{}{%
\LABEL{exm 261}\dpd{1}%
On a set $X$ of $6$ elements, a relation $R$ can
satisfy semi-order property 1 and be anti-symmetric, left and right
serial, but neither transitive nor incomparability-transitive
(\QRef{261}).
%~
An example is shown in Fig.~\REF{Counter-example 261}.
%~
It is not transitive, since $cRe \land eRb$ but not $cRb$;
it is not incomparability-transitive, since $a,c$ and $a,e$ are
incomparable, but $cRe$.
%~
On a set of $\leq 5$ elements, no relation with these properties exists.
\qed
}

\begin{figure}
\begin{center}
\includegraphics[scale=0.5]{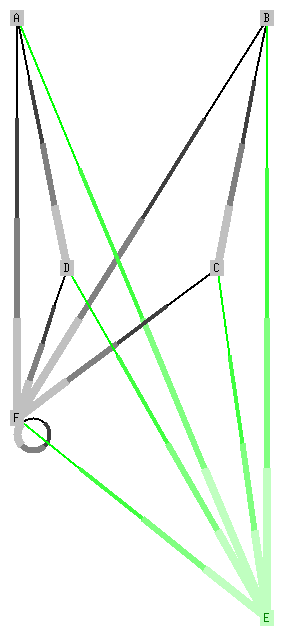}
\caption{Relation graph in Exm.~\REF{exm 262}}
%%\caption{Relation graph in Exm.~\REF{exm 262} (\QRef{262})}
\LABEL{Counter-example 262}
\end{center}
\end{figure}

%	Checking 262: Trans, SemiOrd2
%	==> ASym || SemiOrd1 || LfSerial || RgSerial
%	x\y      0 1 2 3 4 5
%
%	 0       . . . . . .
%	 1       . . . . . .
%	 2       . X . . . .
%	 3       X . . . . .
%	 4       X X X X . X
%	 5       X X X X . X
%
%	tr 1, s2 1, as 0, s1 0, ls 0, rs 0

\EXAMPLE{}{%
\LABEL{exm 262}\dpd{1}%
On a set $X$ of $6$ elements, a relation $R$ can be transitive, but
neither asymmetric nor left nor right serial,
and satisfy semi-order property 2, but not 1
(\QRef{262}).
%~
An example is shown in Fig.~\REF{Counter-example 262}.
%~
It is not asymmetric, since $fRf$;
it is not left and right serial since $e$ and $a$ has no predecessor
and successor, respectively;
it violates semi-order property 1 since $dRa$, $a,c$ are incomparable,
$cRb$, but not $dRb$.
%~
On a set of $\leq 5$ elements, no relation with these properties exists.
\qed
}

\begin{figure}
\begin{center}
\includegraphics[scale=0.5]{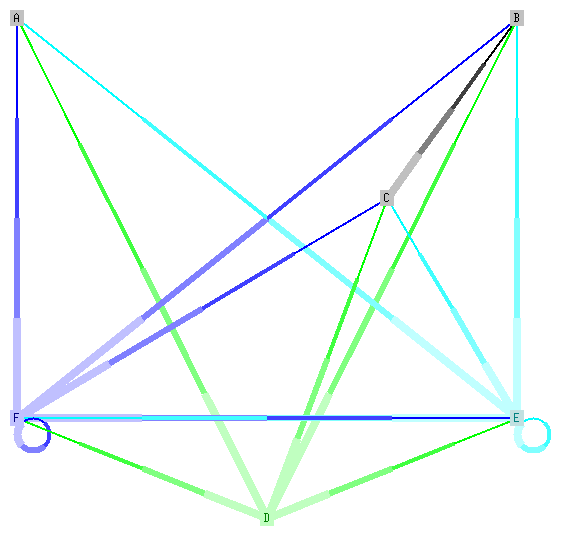}
\caption{Relation graph in Exm.~\REF{exm 263}}
%%\caption{Relation graph in Exm.~\REF{exm 263} (\QRef{263})}
\LABEL{Counter-example 263}
\end{center}
\end{figure}

%	Checking 263: Trans, SemiOrd2
%	==> AntiSym || IncTrans || LfSerial || RgSerial
%	x\y      0 1 2 3 4 5
%
%	 0       . . . . . .
%	 1       . . . . . .
%	 2       . X . . . .
%	 3       X X X . X X
%	 4       X X X . X X
%	 5       X X X . X X
%
%	tr 1, s2 1, as 0, it 0, ls 0, rs 0

\EXAMPLE{}{%
\LABEL{exm 263}\dpd{1}%
On a set $X$ of $6$ elements, a relation $R$ can satisfy semi-order
property 2 and be transitive but neither anti-symmetric nor
incomparability-transitive nor left nor right serial
(\QRef{263}).
%~
An example is shown in Fig.~\REF{Counter-example 263}.
%~
It is not anti-symmetric since $eRf$ and $fRe$;
it is not incomparability-transitive since $a,b$ and $a,c$ are
incomparable but $cRb$;
it is not left and right serial, since $d$ and $a$ has no predecessor
and successor, respectively.
%~
On a set of $\leq 5$ elements, no relation with these properties exists.
\qed
}

\begin{figure}
\begin{center}
\includegraphics[scale=0.5]{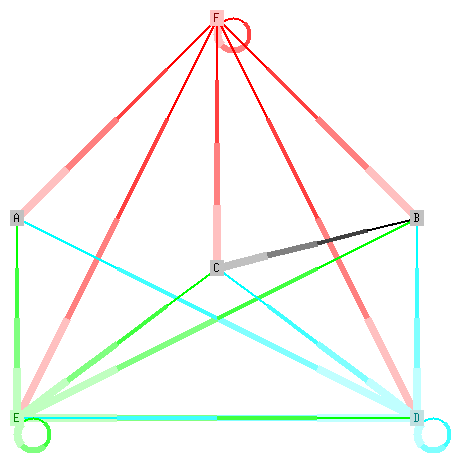}
\caption{Relation graph in Exm.~\REF{exm 264}}
%%\caption{Relation graph in Exm.~\REF{exm 264} (\QRef{264})}
\LABEL{Counter-example 264}
\end{center}
\end{figure}

%	Checking 264: Trans, SemiOrd2, LfSerial, RgSerial
%	==> AntiSym || IncTrans
%	x\y      0 1 2 3 4 5
%
%	 0       . . . . . X
%	 1       . . . . . X
%	 2       . X . . . X
%	 3       X X X X X X
%	 4       X X X X X X
%	 5       . . . . . X
%
%	tr 1, s2 1, ls 1, rs 1, as 0, it 0

\EXAMPLE{}{%
\LABEL{exm 264}\dpd{1}%
On a set $X$ of $6$ elements, a relation $R$ can satisfy semi-order
property 2 and be transitive and left and right serial, but neither
anti-symmetric nor incomparability-transitive
(\QRef{264}).
%~
An example is shown in Fig.~\REF{Counter-example 264}.
%~
It is not anti-symmetric, since $eRd$ and $dRe$;
it is not incomparability-transitive since $a,b$ and $a,c$ are
incomparable, but $cRb$.
%~
On a set of $\leq 5$ elements, no relation with these properties exists.
\qed
}

\begin{figure}
\begin{center}
\includegraphics[scale=0.5]{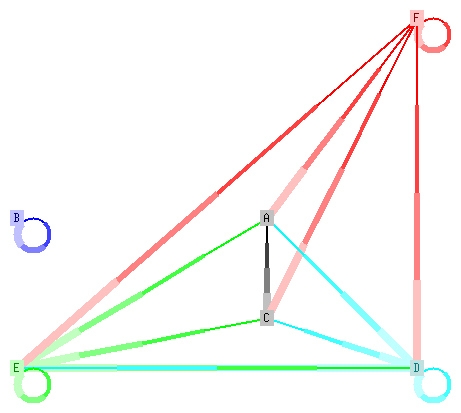}
\caption{Relation graph in Exm.~\REF{exm 266}}
%%\caption{Relation graph in Exm.~\REF{exm 266} (\QRef{266})}
\LABEL{Counter-example 266}
\end{center}
\end{figure}

%	Checking 266: Trans, LfSerial, RgSerial
%	==> SemiOrd1 || AntiSym || Dense
%	x\y      0 1 2 3 4 5
%
%	 0       . . . . . X
%	 1       . X . . . .
%	 2       X . . . . X
%	 3       X . X X X X
%	 4       X . X X X X
%	 5       . . . . . X
%
%	tr 1, ls 1, rs 1, s1 0, as 0, de 0

\EXAMPLE{}{%
\LABEL{exm 266}\dpd{1}%
On a set $X$ of $6$ elements, a relation $R$ can be transitive and
left and right serial, but neither anti-symmetric nor dense nor
satisfying semi-order property 1
(\QRef{266}).
%~
An example is shown in Fig.~\REF{Counter-example 266}.
%~
It is not anti-symmetric, since $dRe$ and $eRd$;
it is not dense, since $cRa$ has no intermediate element;
it doesn't satisfy semi-order property 1, since $bRb$, $b,c$ are
incomparable, $cRa$, but not $bRa$.
%~
On a set of $\leq 5$ elements, no relation with these properties exists.
\qed
}

\begin{figure}
\begin{center}
\includegraphics[scale=0.5]{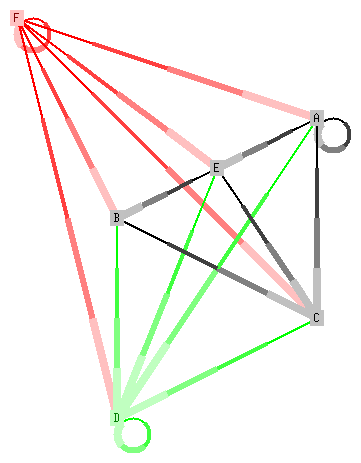}
\caption{Relation graph in Exm.~\REF{exm 267}}
%%\caption{Relation graph in Exm.~\REF{exm 267} (\QRef{267})}
\LABEL{Counter-example 267}
\end{center}
\end{figure}

%	Checking 267: SemiOrd1, AntiSym, LfSerial, RgSerial
%	==> Trans || Dense
%	x\y      0 1 2 3 4 5
%
%	 0       . . . . . X
%	 1       . . . . X X
%	 2       X X . . X X
%	 3       X X X X X X
%	 4       X . . . X X
%	 5       . . . . . X
%
%	s1 1, as 1, ls 1, rs 1, tr 0, de 0

\EXAMPLE{}{%
\LABEL{exm 267}\dpd{1}%
On a set $X$ of $6$ elements, a relation $R$ can be satisfy semi-order
property 1 and be anti-symmetric and left and right serial, but
neither transitive nor dense
(\QRef{267}).
%~
An example is shown in Fig.~\REF{Counter-example 267}.
%~
It is not transitive, since $bRe \land eRa$, but not $bRa$;
it is not dense, since $bRe$ has no intermediate element.
%~
On a set of $\leq 5$ elements, no relation with these properties exists.
\qed
}

\EXAMPLE{}{%
\LABEL{exm 268}\dpd{1}%
On a set $X$ of $\leq 5$ elements,
a transitive and left and right serial relation $R$
must be anti-symmetric, or semi-connex, or dense (\QRef{268}).
%~
On the $6$-element set
$X = \set{a,b,c,d,e,f}$,
the relation
$R = \set{ \tpl{a,a}, \tpl{a,b}, \tpl{b,a}, \tpl{b,b},
\tpl{c,c}, \tpl{c,d}, \tpl{c,e}, \tpl{c,f},
\tpl{d,e}, \tpl{d,f}, \tpl{e,f}, \tpl{f,f} }$
is a counter-example.

Both claims have been machine-checked.
%~
For the $6$-element counter-example,
$R$ isn't anti-symmetric due to $aRb \land bRa$,
not semi-connex since $a$ and $c$ are incomparable, and
not dense since $dRe$ has no intermediate element;
left and right seriality has been achieved by making the elements at
the start and at the end of each chain reflexive.
\qed
}

\begin{figure}
\begin{center}
\includegraphics[scale=0.5]{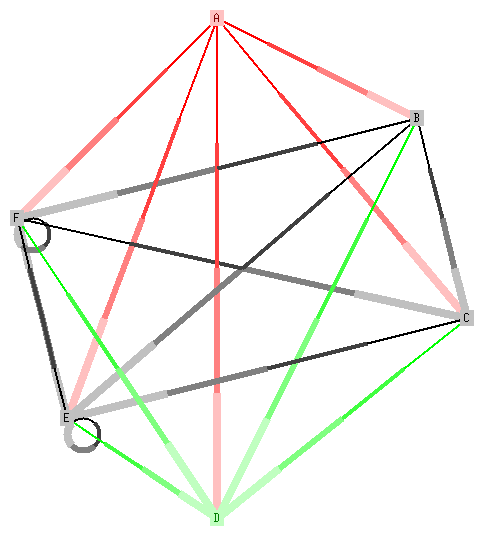}
\caption{Relation graph in Exm.~\REF{exm 269}}
%%\caption{Relation graph in Exm.~\REF{exm 269} (\QRef{269})}
\LABEL{Counter-example 269}
\end{center}
\end{figure}

%	Checking 269: SemiOrd1, SemiConnex
%	==> QuasiTrans || Dense || LfSerial || RgSerial
%	x\y      0 1 2 3 4 5
%
%	 0       . . . . . .
%	 1       X . . . . .
%	 2       X X . . . X
%	 3       X X X . X X
%	 4       X X X . X X
%	 5       X X . . X X
%
%	s1 1, sc 1, qt 0, de 0, ls 0, rs 0

\EXAMPLE{}{%
\LABEL{exm 269}\dpd{1}%
On a set $X$ of $6$ elements, a relation $R$ can satisfy semi-order
property 1 and be semi-connex, but neither quasi-transitive nor dense
nor left nor right serial
(\QRef{269}).
%~
An example is shown in Fig.~\REF{Counter-example 269}.
%~
It is not quasi-transitive, since $eRc \land \lnot cRe$
and $cRf \land \lnot fRc$, but $fRe$;
it is not dense, since $bRa$ has no intermediate element;
it is not left and right serial, since $d$ and $a$ has no predecessor
and successor, respectively.
%~
On a set of $\leq 5$ elements, no relation with these properties exists.
\qed
}

\begin{figure}
\begin{center}
\includegraphics[scale=0.5]{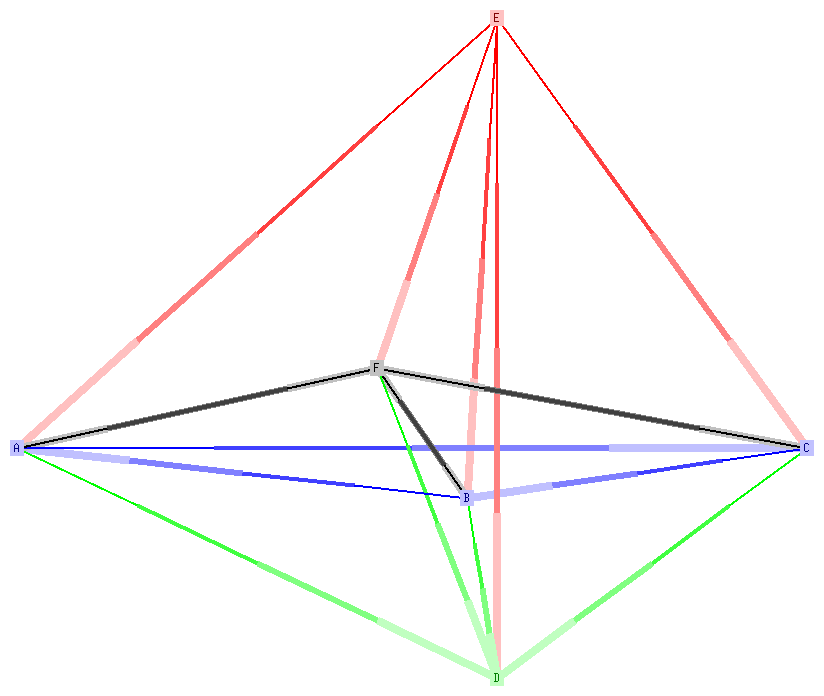}
\caption{Relation graph in Exm.~\REF{exm 270}}
%%\caption{Relation graph in Exm.~\REF{exm 270} (\QRef{270})}
\LABEL{Counter-example for Exm. exm 270}
\end{center}
\end{figure}

\EXAMPLE{}{%
\LABEL{exm 270}\dpd{1}%
On a set $X$ of $\leq 5$ elements,
an irreflexive, semi-connex, and dense relation $R$
must be quasi-transitive, or left or right serial (\QRef{270}).
%~
On a larger set, this is no longer true;
Fig.~\REF{Counter-example for Exm. exm 270}
shows a counter-example for the $6$-element set $X = \set{a,b,c,d,e,f}$,

Both claims have been machine-checked.
%~
In Fig.~\REF{Counter-example for Exm. exm 270},
the cycle $a$, $b$, $c$ (shown in blue) violates quasi-transitivity,
the minimal and maximal element $d$ and $e$ (green and red) violates
left and right seriality, respectively.
%~
Element $f$ (grey) is related to each of $a,b,c$ in both directions,
e.g.\ $aRf \land fRa$ holds, thus helping to establish density.
\qed
}

\begin{figure}
\begin{center}
\includegraphics[scale=0.5]{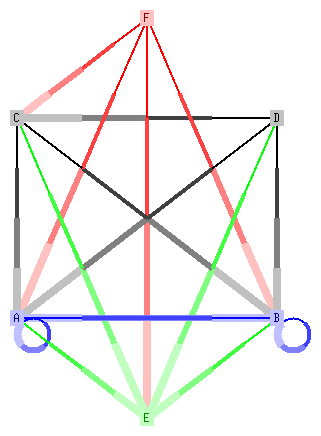}
\caption{Relation graph in Exm.~\REF{exm 271}}
%%\caption{Relation graph in Exm.~\REF{exm 271} (\QRef{271})}
\LABEL{Counter-example for Exm. exm 271}
\end{center}
\end{figure}

%	Checking 271: Trans, IncTrans ==>
%		AntiSym || SemiConnex || Dense || LfSerial || RgSerial
%
%		 D F C E A B
%	x\y      0 1 2 3 4 5
%
%	D 0      . . . . . .
%	F 1      . . . . . .
%	C 2      X X . . . .
%	E 3      X X X . X X
%	A 4      X X X . X X
%	B 5      X X X . X X

\EXAMPLE{}{%
\LABEL{exm 271}\dpd{1}%
On a set $X$ of $\leq 5$ elements,
an transitive and incomparability-transitive relation $R$
must be anti-symmetric, semi-connex, dense,
or left or right serial (\QRef{271}).
%~
On a larger set, this is no longer true;
Fig.~\REF{Counter-example for Exm. exm 271}
shows a counter-example for the $6$-element set $X = \set{a,b,c,d,e,f}$,
%~
It is not anti-symmetric, since $aRb$ and $bRa$;
it is not semi-connex, since $d,f$ are incomparable;
it is not dense, since $cRf$ has no intermediate element;
it is not left and right serial, since $e$ and $f$ has no predecessor
and successor, respectively.
\qed
}

\begin{figure}
\begin{center}
\includegraphics[scale=0.5]{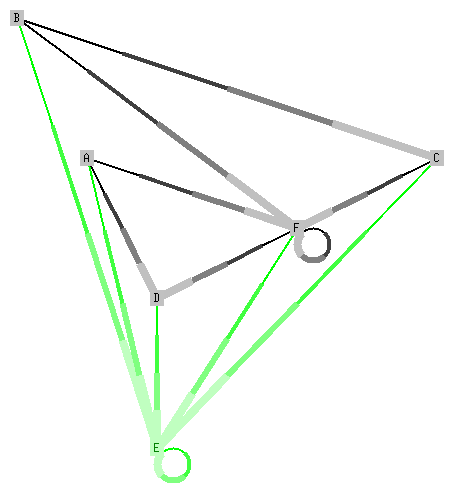}
\caption{Relation graph in Exm.~\REF{exm 272}}
%%\caption{Relation graph in Exm.~\REF{exm 272} (\QRef{272})}
\LABEL{Counter-example 272}
\end{center}
\end{figure}

%	Checking 272:  SemiOrd1, SemiOrd2 LfSerial ==>
%	IncTrans || QuasiTrans || Dense || RgSerial
%	x\y      0 1 2 3 4 5
%
%	 0       . . . . . .
%	 1       . . . . . .
%	 2       . X . . . .
%	 3       . X . . . X
%	 4       X X X X X X
%	 5       X X X . . X
%
%	s1 1, s2 1, ls 1, it 0, qt 0, de 0, rs 0

% 273 dual zu 272

\EXAMPLE{}{%
\LABEL{exm 272}\dpd{1}%
On a set $X$ of $6$ elements, a relation $R$ can satisfy semi-order
properties 1 and 2 and be left serial but neither right serial nor
incomparability-transitive nor quasi-transitive nor dense
(\QRef{272}).
%~
An example is shown in Fig.~\REF{Counter-example 272}.
%~
It is not right serial, since $b$ has no successor;
it is not incomparability-transitive, since $b,a$ and $b,d$ are
incomparable, but $dRa$;
it is not quasi-transitive,
since $dRf \land \lnot fRd$ and $fRc \land \lnot cRf$, but not $dRc$;
it is not dense, since $cRb$ has no intermediate element.
%~
By reverting the arrow directions, a counter-example for the dual
\QRef{273} is obtained.
%~
On a set of $\leq 5$ elements, no relation with either property combination
exists.
\qed
}

\begin{figure}
\begin{center}
\includegraphics[scale=0.5]{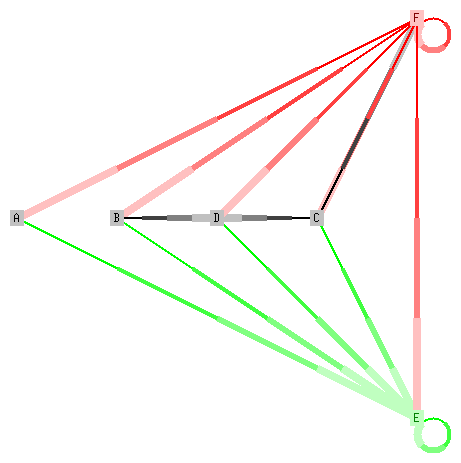}
\caption{Relation graph in Exm.~\REF{exm 274}}
%%\caption{Relation graph in Exm.~\REF{exm 274} (\QRef{274})}
\LABEL{Counter-example 274}
\end{center}
\end{figure}

% 274		 A B C D E F
%	x\y      0 1 2 3 4 5
%
%	 0       . . . . . X
%	 1       . . . . . X
%	 2       . . . . . X
%	 3       . X X . . X
%	 4       X X X X X X
%	 5       . . X . . X
%
%	s1 1, s2 1, qt 1, ls 1, rs 1, tr 0, it 0, de 0

\EXAMPLE{}{%
\LABEL{exm 274}\dpd{1}%
On a set $X$ of $6$ elements, a relation $R$ can satisfy semi-order
properties 1 and 2 and be quasi-transitive and left and right serial,
but neither transitive nor incomparability-transitive nor dense
(\QRef{274}).
%~
An example is shown in Fig.~\REF{Counter-example 274}.
%~
It is not transitive, since $cRf \land fRc$ but not $cRc$;
it is not incomparability-transitive, since
$a,b$ and $a,d$ are incomparable, but $dRb$;
it is not dense, since $dRb$ has no intermediate element.
%~
On a set of $\leq 5$ elements, no relation with these properties exists.
\qed
}

\clearpage
\section{Implementation issues}
\LABEL{Implementation issues}

In this section, we comment on some program details.
%~
The source code is available in the ancillary file 
\texttt{nonprominentProperties.c}.

\subsection{Improved relation enumeration}
\LABEL{Improved relation enumeration}

\begin{figure}
\begin{center}
$\begin{array}{l|cccc|l @{\hspace*{2cm}} l|cccc|l}
\multicolumn{1}{l}{} & a & b & c & \multicolumn{1}{c}{d} &&
\multicolumn{1}{l}{} & b & c & d & \multicolumn{1}{c}{a} & \\
\cline{2-5}
\cline{8-11}
a & 1 & 1 & 1 & 1 & \tpl{3,1} & b & 1 & 1 & . & . & \tpl{1,1}	\\
b & . & 1 & 1 & . & \tpl{1,1} & c & . & . & 1 & 1 & \tpl{2,0}	\\
c & 1 & . & . & 1 & \tpl{2,0} & d & . & 1 & . & 1 & \tpl{2,0}	\\
d & 1 & . & 1 & . & \tpl{2,0} & a & 1 & 1 & 1 & 1 & \tpl{3,1}	\\
\cline{2-5}
\cline{8-11}
\end{array}$
\caption{Example relation (left) and its normal form (right)}
\LABEL{Example relation (left) and its normal form (right)}
\end{center}
\end{figure}

We encoded a binary relation by a square array of \texttt{bool}s.
%~
Originally, we enumerated all possible assignments of such an array.
%~
However, if $R$ is a relation on a finite set $X$ of $n$ elements,
and $\pi: X \ra X$ is a permutation,
then $R'$ defined by $x R' y :\Lra \pi_x R \pi_y$
shares all properties from Def.~\REF{def} with $R$;
that is, $R'$ is reflexive iff $R$ is, etc.

In order to speed up relation enumeration, we therefore defined a
normal form for a square array as follows:
%~
To the $i$th array row, assign the pair $\tpl{c_i,d_i}$, where $c_i$ is the
number of \texttt{true} values in columns $j \neq i$, and $d_i$ is
the value of the $i$th column.
%~
An array is in normal form if
$\:{\leq}1n{ \tpl{c_\i,d_\i} }$ holds,
where ``$\leq$'' denotes the lexicographic order.

Every array can be converted into normal form by applying a
row-sorting permutation simultaneous to rows and columns.
%~
Note that the pairs are designed to be invariant under simultaneous
row and column permutation.

As an example,
the left half of
Fig.~\REF{Example relation (left) and its normal form (right)}
shows a relation on the set $X = \set{a,b,c,d}$,
and the pairs assigned to each row.
%~
For readability, we denoted the value \texttt{false} and \texttt{true}
by ``.'' and ``$1$'', respectively.
%~
The right half shows the corresponding normal form, obtained by
sorting the rows by ascending associated pairs, and permuting the
columns in the same way.

\begin{figure}
\begin{center}
\small
\begin{verbatim}
const bool allRows[72] = { 
  /*  0: ***** row position 0 ********************* */
  /*  0: group <0,0> */   0,0,0,    
  /*  3: group <0,1> */   1,0,0,
  /*  6: group <1,0> */   0,1,0, 0,0,1,
  /* 12: group <1,1> */   1,1,0, 1,0,1,
  /* 18: group <2,0> */   0,1,1,       
  /* 21: group <2,1> */   1,1,1,
  /* 24: ***** row position 1 ********************* */
  /* 24: group <0,0> */   0,0,0,       
  /* 27: group <0,1> */   0,1,0,
  /* 30: group <1,0> */   1,0,0, 0,0,1,
  /* 36: group <1,1> */   1,1,0, 0,1,1,
  /* 42: group <2,0> */   1,0,1,       
  /* 45: group <2,1> */   1,1,1,
  /* 48: ***** row position 2 ********************* */ 
  /* 48: group <0,0> */   0,0,0,       
  /* 51: group <0,1> */   0,0,1,
  /* 54: group <1,0> */   1,0,0, 0,1,0,
  /* 60: group <1,1> */   1,0,1, 0,1,1,
  /* 66: group <2,0> */   1,1,0,       
  /* 69: group <2,1> */   1,1,1,
  /* 72: ***** end ******************************** */ 
};
const int gpS[3*2+1] = {  0,  3,  6, 12, 18, 21, 24 };
int gp[card];
int rw[card];

void initFromRow(bool R[card][card]) {
  for (int r=0; r<card; ++r)
    for (int c=0; c<card; ++c)
      R[r][c] = allRows[ r*24 + rw[r] + cc ];
}

void check03(void) {
  bool R[card][card];
  for (gp[0]=0; gp[0]<card*2; ++gp[0])
    for (gp[1]=0; gp[1]<=gp[0]; ++gp[1])
      for (gp[2]=0; gp[2]<=gp[1]; ++gp[2])
        for (rw[0]=gpS[gp[0]]; rw[0]<gpS[gp[0]+1]; rw[0]+=card)
          for (rw[1]=gpS[gp[1]]; rw[1]<gpS[gp[1]+1]; rw[1]+=card)
            for (rw[2]=gpS[gp[2]]; rw[2]<gpS[gp[2]+1]; rw[2]+=card) {
              initFromRow(R);
              if (isRgEucl(R) && ! isTrans(R))
                printRel(R);
            }
}
\end{verbatim}
\caption{Improved code to search for right Euclidean non-transitive relations}
\LABEL{Improved code to search for right Euclidean non-transitive relations}
\end{center}
\end{figure}

It is sufficient to consider property combinations only for relation
arrays that are
in normal form, and this is what our improved algorithm does.
%~
Figure~\REF{Improved code to search for right Euclidean non-transitive relations}
shows the improved code to search for right Euclidean non-transitive
relations over a $3$-element universe.\footnote{
	The unimproved code was shown in 
	Fig.~\REF{Source code to search for right Euclidean non-transitive relations}
}
%~
The list \texttt{allRows} contains all possible rows for a 
$3 \times 3$ array representing a relation.
%~
The rows are grouped by assigned pair;
we have one copy for
each row position, to account for the varying column position of the
diagonal element.
%~
For example, in horizontal position $0$, the rows \texttt{1,1,0} and
\texttt{1,0,1} are assigned the pair $\tpl{1,1}$; they are found at
starting index \texttt{12} of \texttt{allRows}.
%~
The corresponding row values for horizontal position $1$ are
\texttt{1,1,0} and \texttt{0,1,1}; they start at index \texttt{36}.
%~
Procedure \texttt{check03} iterates in the loops on \texttt{gp[0]},
\texttt{gp[1]}, and
\texttt{gp[2]} over all combinations of groups
that lead to a normal form, and in the loops on \texttt{rw[0]},
\texttt{rw[1]}, and \texttt{rw[2]},
over all combinations of rows from the current groups.
%~
The lists \texttt{allRows} and \texttt{gpS} were precomputed by
another program; its source code is available in the ancillary file
\texttt{genTables.c}.

Figure~\REF{Number of relations vs. carrier set cardinality}
shows, for set cardinalities $1$ through $7$, the number of all array
assignments (column ``Unpruned''), the number of arrays in normal form
(column ``Pruned''), and the quotient of both numbers, indicating the
speed-up factor.
%~
Figure~\REF{Number of relations on a $5$-element set}
shows, for each property, the number of satisfying relations found
with  the old (column ``Old'') and with the improved
(column ``Pruned'') enumeration method.

\begin{figure}
\begin{center}
\begin{tabular}{|r|r|r|r|}
\hline
\bf Card & \bf Unpruned & \bf Pruned & \bf Ratio	\\
\hline
% 0&                                        1&               &	\\
 1&                                        2&              2&  1.0	\\
 2&                                       16&             10&  1.6	\\
 3&                                      512&            140&  3.6	\\
 4&                                   65 536&          6 170& 10.6	\\
 5&                               33 554 432&        907 452& 36.9	\\
 6&                           68 719 476 736&    460 631 444&149.1	\\
 7&                      562 949 953 421 312&827 507 617 792&680.2	\\
%8&               18 446 744 073 709 551 616&               &  \\
% 9&        2 417 851 639 229 258 349 412 352&               &  \\
%10&1 267 650 600 228 229 401 496 703 205 376&               &  \\
\hline
\end{tabular}
\caption{Number of relations vs.\ carrier set cardinality}
\LABEL{Number of relations vs. carrier set cardinality}
\end{center}
\end{figure}

\subsection{Quine-McCluskey implementation}

The procedure \texttt{computeLaws} iterates over all
relations, determining for each the set\footnote{
	encoded as bit vector, see
	Fig.~\REF{Number of relations on a $5$-element set}
}
of its properties, and counting the number of occurrences of each such
vector.
%~
After that, it calls the Quine-McCluskey implementation \texttt{qmc}
to compute all prime implicants of the non-occurring vectors.
%~
The latter procedure performs a top-down breadth-first search on the
search graph.

\begin{figure}
\begin{center}
\begin{tabular}{|r|r|l|l|l|}
\hline
\bf Count & \bf Count & \bf Name & \bf Def. & \bf Encoding	\\
\bf (Old) & \bf (Pruned) & &  & 	\\
\hline
         1 &       1 & Empty		&			& \verb|0x000001|	\\
         1 &       1 & Univ		&			& \verb|0x000002|	\\
        32 &       6 & CoRefl		& \REFF{def}{CoRefl}	& \verb|0x000004|	\\
     3 163 &     166 & LfEucl		& \REFF{def}{LfEucl}	& \verb|0x000008|	\\
     3 163 &     186 & RgEucl		& \REFF{def}{RgEucl}	& \verb|0x000010|	\\
     7 776 &     440 & LfUnique		& \REFF{def}{LfUnique}	& \verb|0x000020|	\\
     7 776 &   1 818 & RgUnique		& \REFF{def}{RgUnique}	& \verb|0x000040|	\\
    32 768 &   1 012 & Sym		& \REFF{def}{Sym}	& \verb|0x000080|	\\
    47 462 &   4 841 & AntiTrans	& \REFF{def}{AntiTrans}	& \verb|0x000100|	\\
    59 049 &   3 870 & ASym		& \REFF{def}{ASym}	& \verb|0x000200|	\\
    59 049 &   3 870 & Connex		& \REFF{def}{Connex}	& \verb|0x000400|	\\
   154 303 &   3 207 & Trans		& \REFF{def}{Trans}	& \verb|0x000800|	\\
   467 750 &  11 103 & SemiOrd1		& \REFF{def}{SemiOrd1}	& \verb|0x001000|	\\
 1 048 576 &  70 436 & Irrefl		& \REFF{def}{Irrefl}	& \verb|0x002000|	\\
 1 048 576 &  70 436 & Refl		& \REFF{def}{Refl}	& \verb|0x004000|	\\
 1 069 742 &  71 198 & QuasiRefl	& \REFF{def}{QuasiRefl}	& \verb|0x008000|	\\
 1 889 568 &  50 480 & AntiSym		& \REFF{def}{AntiSym}	& \verb|0x010000|	\\
 1 889 568 &  50 480 & SemiConnex	& \REFF{def}{SemiConnex}& \verb|0x020000|	\\
 3 756 619 & 113 142 & IncTrans		& \REFF{def}{IncTrans}	& \verb|0x040000|	\\
 4 498 393 & 144 128 & SemiOrd2		& \REFF{def}{SemiOrd2}	& \verb|0x080000|	\\
 5 531 648 & 131 994 & QuasiTrans	& \REFF{def}{QuasiTrans}& \verb|0x100000|	\\
15 339 497 & 425 854 & Dense		& \REFF{def}{Dense}	& \verb|0x200000|	\\
28 629 151 & 764 962 & LfSerial		& \REFF{def}{LfSerial}	& \verb|0x400000|	\\
28 629 151 & 817 185 & RgSerial		& \REFF{def}{RgSerial}	& \verb|0x800000|	\\
\hline
\end{tabular}
\caption{Number of relations on a $5$-element set}
\LABEL{Number of relations on a $5$-element set}
\end{center}
\end{figure}

An example graph, showing all possible prime implicants
for a Boolean function of 3 variables is given in
Fig.~\REF{Search graph for the Quine-McCluskey algorithm on 3 variables}.
%~
At each node of the search graph, the corresponding conjunction is
checked by the procedure \texttt{qmcRect}:
if no combination in its covered rectangle\footnote{
	This terminology is inspired by the Karnaugh diagram method;
	in
	Fig.~\REF{Search graph for the Quine-McCluskey algorithm on 3 variables},
	the rectangle covered by a node
	corresponds to the set of all leaves below it.
}
is ``\texttt{off}'' and
at least one is ``\texttt{on}'',\footnote{
	Since we are interested in \emph{non}-occurring vectors,
	``\texttt{on}'' corresponds to an occurrence count of zero,
	and ``\texttt{off}'' to a count $>0$.
	%~
	We encode ``\texttt{don't care}'' by a count of $-1$.
}
then it is actually a prime implicant.
%~
In that case, we output its description using
\texttt{qmcPrint},\footnote{
	In particular, we don't perform the usual search of a
	\emph{minimal}
	set of prime implicants covering all ``\texttt{on}'' vectors.
}
and set all vectors in its
covered rectangle to \texttt{don't care}.

Note that we can't perform a depth-first search: for example, if $a$
isn't a prime implicant, we can't check its child $ab$ next, since it
could satisfy the above primeness criterion, but nevertheless be
covered by a simpler prime implicant, such as $b$.

\begin{figure}
\begin{center}
	\newcommand{\node}[3]{%
		\put(#1,#2){\circle*{2}}%
		\put(#1,#2){\put(2,0){\makebox(0,0)[l]{$#3$}}}%
	}
\begin{picture}(110,60)%
	%\put(0,0){\makebox(0,0){$+$}}%
	%\put(110,60){\makebox(0,0){$+$}}%
% 0-variable conjuncts
%\put(55,60){\Line{  0,40}}%	true - a
%\put(55,60){\Line{ 22,40}}%	true - A
%\put(55,60){\Line{ 44,40}}%	true - b
%\put(55,60){\Line{ 66,40}}%	true - B
%\put(55,60){\Line{ 88,40}}%	true - c
%\put(55,60){\Line{110,40}}%	true - C
\put(55.000,60.000){\line(-3,-1){30.000}}%...
 \put(25.000,50.000){\line(-5,-2){25.000}}%	true - a
\put(55.000,60.000){\line(-5,-3){33.000}}%	true - A
%\put(55.000,60.000){\line(-5,-3){30.000}}%...
% \put(25.000,42.000){\line(-3,-2){3.000}}%	true - A
\put(55.000,60.000){\line(-3,-5){6.000}}%...
 \put(49.000,50.000){\line(-1,-2){5.000}}%	true - b
\put(55.000,60.000){\line(1,-2){5.000}}%...
 \put(60.000,50.000){\line(3,-5){6.000}}%	true - B
\put(55.000,60.000){\line(5,-3){33.000}}%	true - c
%\put(55.000,60.000){\line(3,-2){3.000}}%...
% \put(58.000,58.000){\line(5,-3){30.000}}%	true - c
\put(55.000,60.000){\line(5,-2){25.000}}%...
 \put(80.000,50.000){\line(3,-1){30.000}}%	true - C
\node{ 55}{60}{\true}%	true
% 1-variable conjuncts
%\put(  0,40){\Line{  0,20}%	a - ab
%\put(  0,40){\Line{ 10,20}%	a - aB
%\put(  0,40){\Line{ 20,20}%	a - ac
%\put(  0,40){\Line{ 30,20}%	a - aC
%\put( 22,40){\Line{ 40,20}%	A - Ab
%\put( 22,40){\Line{ 50,20}%	A - AB
%\put( 22,40){\Line{ 60,20}%	A - Ac
%\put( 22,40){\Line{ 70,20}%	A - AC
%\put( 44,40){\Line{  0,20}%	b - ab
%\put( 44,40){\Line{ 40,20}%	b - Ab
%\put( 44,40){\Line{ 80,20}%	b - bc
%\put( 44,40){\Line{ 90,20}%	b - bC
%\put( 66,40){\Line{ 10,20}%	B - aB
%\put( 66,40){\Line{ 50,20}%	B - AB
%\put( 66,40){\Line{100,20}%	B - Bc
%\put( 66,40){\Line{110,20}%	B - BC
%\put( 88,40){\Line{ 20,20}%	c - ac
%\put( 88,40){\Line{ 60,20}%	c - Ac
%\put( 88,40){\Line{ 80,20}%	c - bc
%\put( 88,40){\Line{100,20}%	c - Bc
%\put(110,40){\Line{ 30,20}%	C - aC
%\put(110,40){\Line{ 70,20}%	C - AC
%\put(110,40){\Line{ 90,20}%	C - bC
%\put(110,40){\Line{110,20}%	C - BC
\put(0.000,40.000){\line(0,-1){20.000}}%	a - ab
\put(0.000,40.000){\line(1,-2){10.000}}%	a - aB
\put(0.000,40.000){\line(1,-1){20.000}}%	a - ac
\put(0.000,40.000){\line(3,-2){30.000}}%	a - aC
\put(22.000,40.000){\line(5,-6){10.000}}%...
 \put(32.000,28.000){\line(1,-1){8.000}}%	A - Ab
\put(22.000,40.000){\line(4,-3){16.000}}%...
 \put(38.000,28.000){\line(3,-2){12.000}}%	A - AB
\put(22.000,40.000){\line(5,-3){10.000}}%...
 \put(32.000,34.000){\line(2,-1){28.000}}%	A - Ac
\put(22.000,40.000){\line(2,-1){8.000}}%...
 \put(30.000,36.000){\line(5,-2){40.000}}%	A - AC
\put(44.000,40.000){\line(-5,-2){20.000}}%...
 \put(24.000,32.000){\line(-2,-1){24.000}}%	b - ab
\put(44.000,40.000){\line(-1,-5){4.000}}%	b - Ab
\put(44.000,40.000){\line(5,-3){20.000}}%...
 \put(64.000,28.000){\line(2,-1){16.000}}%	b - bc
\put(44.000,40.000){\line(2,-1){16.000}}%...
 \put(60.000,32.000){\line(5,-2){30.000}}%	b - bC
\put(66.000,40.000){\line(-3,-1){36.000}}%...
 \put(30.000,28.000){\line(-5,-2){20.000}}%	B - aB
\put(66.000,40.000){\line(-4,-5){16.000}}%	B - AB
\put(66.000,40.000){\line(5,-3){30.000}}%...
 \put(96.000,22.000){\line(2,-1){4.000}}%	B - Bc
\put(66.000,40.000){\line(2,-1){24.000}}%...
 \put(90.000,28.000){\line(5,-2){20.000}}%	B - BC
\put(88.000,40.000){\line(-4,-1){32.000}}%...
 \put(56.000,32.000){\line(-3,-1){36.000}}%	c - ac
\put(88.000,40.000){\line(-3,-2){12.000}}%...
 \put(76.000,32.000){\line(-4,-3){16.000}}%	c - Ac
\put(88.000,40.000){\line(-2,-5){8.000}}%	c - bc
\put(88.000,40.000){\line(3,-5){12.000}}%	c - Bc
\put(110.000,40.000){\line(-4,-1){80.000}}%	C - aC
\put(110.000,40.000){\line(-2,-1){40.000}}%	C - AC
\put(110.000,40.000){\line(-1,-1){20.000}}%	C - bC
\put(110.000,40.000){\line(0,-1){20.000}}%	C - BC
\node{  0}{40}{a}%		a
\node{ 22}{40}{\overline{a}}%	A
\node{ 44}{40}{b}%		b
\node{ 66}{40}{\overline{b}}%	B
\node{ 88}{40}{c}%		c
\node{110}{40}{\overline{c}}%	C
% 2-variable conjuncts
%\put(  0,20)\Line{  0,0}%	ab - abc
%\put(  0,20)\Line{ 15,0}%	ab - abC
%\put( 10,20)\Line{ 30,0}%	aB - aBc
%\put( 10,20)\Line{ 45,0}%	aB - aBC
%\put( 20,20)\Line{  0,0}%	ac - abc
%\put( 20,20)\Line{ 30,0}%	ac - aBc
%\put( 30,20)\Line{ 15,0}%	aC - abC
%\put( 30,20)\Line{ 45,0}%	aC - aBC
%\put( 40,20)\Line{ 60,0}%	Ab - Abc
%\put( 40,20)\Line{ 75,0}%	Ab - AbC
%\put( 50,20)\Line{ 90,0}%	AB - ABc
%\put( 50,20)\Line{105,0}%	AB - ABC
%\put( 60,20)\Line{ 60,0}%	Ac - Abc
%\put( 60,20)\Line{ 90,0}%	Ac - ABc
%\put( 70,20)\Line{ 75,0}%	AC - AbC
%\put( 70,20)\Line{105,0}%	AC - ABC
%\put( 80,20)\Line{  0,0}%	bc - abc
%\put( 80,20)\Line{ 60,0}%	bc - Abc
%\put( 90,20)\Line{ 15,0}%	bC - abC
%\put( 90,20)\Line{ 75,0}%	bC - AbC
%\put(100,20)\Line{ 30,0}%	Bc - aBc
%\put(100,20)\Line{ 90,0}%	Bc - ABc
%\put(110,20)\Line{ 45,0}%	BC - aBC
%\put(110,20)\Line{105,0}%	BC - ABC
\put(0.000,20.000){\line(0,-1){20.000}}%	ab - abc
\put(0.000,20.000){\line(3,-4){15.000}}%	ab - abC
\put(10.000,20.000){\line(1,-1){20.000}}%	aB - aBc
\put(10.000,20.000){\line(5,-3){25.000}}%...
 \put(35.000,5.000){\line(2,-1){10.000}}%	aB - aBC
\put(20.000,20.000){\line(-1,-1){20.000}}%	ac - abc
\put(20.000,20.000){\line(1,-2){10.000}}%	ac - aBc
\put(30.000,20.000){\line(-3,-4){15.000}}%	aC - abC
\put(30.000,20.000){\line(3,-4){15.000}}%	aC - aBC
\put(40.000,20.000){\line(1,-1){20.000}}%	Ab - Abc
\put(40.000,20.000){\line(5,-3){25.000}}%...
 \put(65.000,5.000){\line(2,-1){10.000}}%	Ab - AbC
\put(50.000,20.000){\line(2,-1){40.000}}%	AB - ABc
\put(50.000,20.000){\line(5,-2){25.000}}%...
 \put(75.000,10.000){\line(3,-1){30.000}}%	AB - ABC
\put(60.000,20.000){\line(0,-1){20.000}}%	Ac - Abc
\put(60.000,20.000){\line(3,-2){30.000}}%	Ac - ABc
\put(70.000,20.000){\line(1,-4){5.000}}%	AC - AbC
\put(70.000,20.000){\line(5,-3){25.000}}%...
 \put(95.000,5.000){\line(2,-1){10.000}}%	AC - ABC
\put(80.000,20.000){\line(-4,-1){80.000}}%	bc - abc
\put(80.000,20.000){\line(-1,-1){20.000}}%	bc - Abc
\put(90.000,20.000){\line(-4,-1){60.000}}%...
 \put(30.000,5.000){\line(-3,-1){15.000}}%	bC - abC
\put(90.000,20.000){\line(-3,-4){15.000}}%	bC - AbC
\put(100.000,20.000){\line(-4,-1){40.000}}%...
 \put(60.000,10.000){\line(-3,-1){30.000}}%	Bc - aBc
\put(100.000,20.000){\line(-1,-2){10.000}}%	Bc - ABc
\put(110.000,20.000){\line(-4,-1){20.000}}%...
 \put(90.000,15.000){\line(-3,-1){45.000}}%	BC - aBC
\put(110.000,20.000){\line(-1,-4){5.000}}%	BC - ABC
\node{  0}{20}{ab}%				ab
\node{ 10}{20}{a\overline{b}}%			aB
\node{ 20}{20}{ac}%				ac
\node{ 30}{20}{a\overline{c}}%			aC
\node{ 40}{20}{\overline{a}b}%			Ab
\node{ 50}{20}{\overline{a}\overline{b}}%	AB
\node{ 60}{20}{\overline{a}c}%			Ac
\node{ 70}{20}{\overline{a}\overline{c}}%	AC
\node{ 80}{20}{bc}%				bc
\node{ 90}{20}{b\overline{c}}%			bC
\node{100}{20}{\overline{b}c}%			Bc
\node{110}{20}{\overline{b}\overline{c}}%	BC
% 3-variable conjuncts
\node{  0}{ 0}{abc}%					abc
\node{ 15}{ 0}{ab\overline{c}}%				abC
\node{ 30}{ 0}{a\overline{b}c}%				aBc
\node{ 45}{ 0}{a\overline{b}\overline{c}}%		aBC
\node{ 60}{ 0}{\overline{a}bc}%				Abc
\node{ 75}{ 0}{\overline{a}b\overline{c}}%		AbC
\node{ 90}{ 0}{\overline{a}\overline{b}c}%		ABc
\node{105}{ 0}{\overline{a}\overline{b}\overline{c}}%	ABC
\end{picture}
\caption{Search graph for the Quine-McCluskey algorithm on 3 variables}
\LABEL{Search graph for the Quine-McCluskey algorithm on 3 variables}
\end{center}
\end{figure}
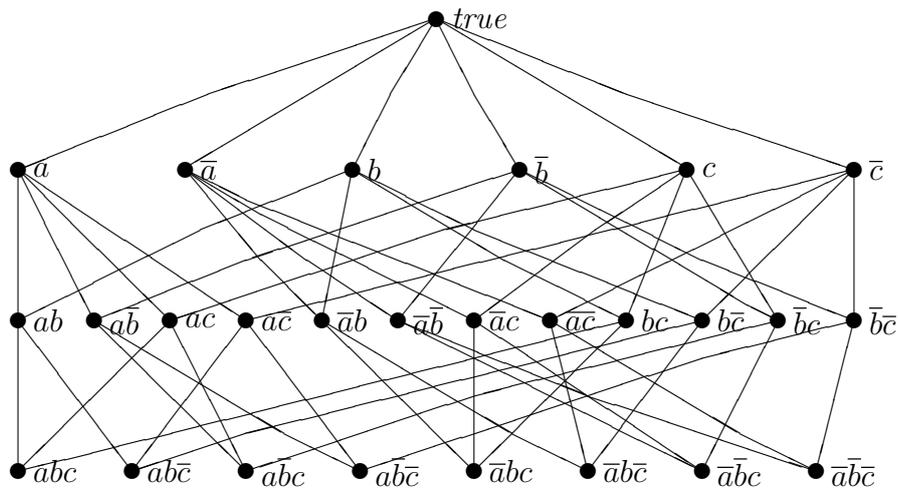

Figure~\REF{Timing for Quine-McCluskey run}
shows the timing for our Quine-McCluskey run.
%~
Enumeration\footnote{
	We still used the old enumeration scheme,
	i.e.\ we enumerated and checked also relations in
	non-normal form.
}
of all relations on a $5$-element set and checking all properties of
each relation was done after $20$ seconds wall-clock time.
%~
After that, levels $1$ to $5$ were completely checked within an hour,
but it took over a day to find the last law, on level $8$.
%~
The largest share of run time was used in looping though all possible
rectangles of a level; level $n$ has
$\left( \begin{array}{c} 24 \\ n \end{array} \right) \cdot 2^n$
rectangles.

Note that the improved relation enumeration described in
Sect.~\REF{Improved relation enumeration}
would have affected only the very first phase, which was completed in $20$
seconds, anyway.
%~
However, it would have allowed for using a $6$-element universe set in
reasonable computation time, thereby avoiding the report of all non-laws in
Sect.~\REF{Examples}, except \QRef{186} (Exm.~\REF{exm 186}), which
needs $\geq 7$ elements in the universe.

%	4/3	3.14
%	5/4	6.76
%	6/5	5.72
%	7/6	4.43
%	8/7	3.36
%	9/8	2.76

Somewhat unexpected, the $2^{25}$ relations inhabited no more than
$495$ of the $2^{24}$ possible combinations of properties.
%~
A listing of property combinations by number of satisfying relations
is provided in the ancillary file
\texttt{propertyCombinationsByCount.txt}.

\begin{figure}
\begin{center}
\begin{tabular}{|r|r|r|r|r|}
\hline
\bf Time & \multicolumn{4}{l|}{}	\\
(sec) & \multicolumn{4}{l|}{}	\\
\hline
     0 & \multicolumn{4}{l|}{Counting property combinations}	\\
    20 & \multicolumn{4}{l|}{Number of relations (Fig.~\REF{Number of relations on a $5$-element set})}	\\
    21 & \multicolumn{4}{l|}{Prime implicants (Fig.~\REF{Reported laws for level 2}--\REF{Reported laws for level 8})}	\\
\cline{2-5}
& \bf Level & \bf On & \bf Off & \bf Don't care	\\
\cline{2-5}
    21 & 1 & 16776721 & 495 &        0	\\
    21 & 2 & 16776721 & 495 &        0	\\
    27 & 3 &    32063 & 495 & 16744658	\\
    85 & 4 &      161 & 495 & 16776560	\\
   575 & 5 &       32 & 495 & 16776689	\\
  3291 & 6 &       24 & 495 & 16776697	\\
 14591 & 7 &        4 & 495 & 16776717	\\
 49134 & 8 &        1 & 495 & 16776720	\\
135925 & 9 &        0 & 495 & 16776721	\\
\hline
\end{tabular}
\caption{Timing for Quine-McCluskey run}
\LABEL{Timing for Quine-McCluskey run}
\end{center}
\end{figure}

\subsection{On finding ``nice'' laws}

It is desirable to find a set of laws as ``elegant'' as possible.
%~
While ``elegance'' is a matter of mathematicians' taste and can hardly
be rigorously defined, some criteria for it are beyond doubt.

Each single law should be as general as possible.
%~
On a technical level, this translates into two requirements:
\begin{enumerate}
\item\LABEL{nice 1}%
	Each law should consist of as few literals as possible.
\item\LABEL{nice 2}%
	Each law should use the sharpest predicates possible.
\end{enumerate}
Considering sets of laws,
\begin{enumerate}
\setcounter{enumi}{2}
\item\LABEL{nice 3}%
	some balance should be kept between
	conciseness and convenience.
\end{enumerate}

We discuss criteria~\REF{nice 1} to~\REF{nice 3} in the following.

\paragraph*{Criterion \REF{nice 1}}

Each law should consist of as few literals as possible.

That
is, it should be obtained from a rectangle as large as possible.
%~
For example, in a law
$\lnot$~LfSerial $\lor$ $\lnot$~Asym $\lor$ Irrefl,
the first literal should be omitted.

This requirement is fulfilled, since our algorithm checks
rectangles in order of decreasing size, and for every reported
rectangle prevents properly contained rectangles from being
reported also.

\paragraph*{Criterion \REF{nice 2}}

Each law should use the sharpest predicates possible.

For example, both
``LfUnique $\land$ RgUnique $\land$ SemiOrd2 $\Ra$ ASym''
(\QRef{226})
and
``LfUnique $\land$ RgUnique $\land$ SemiOrd2 $\Ra$ Irrefl''
is a law,
but the latter is redundant since it follows from the former
and the law ``ASym $\Ra$ Irrefl'' (\QRef{039}).

In a naive approach to cope with this requirement, we ordered
the properties by extension set cardinality, see
Fig.~\REF{Number of relations on a $5$-element set},\footnote{
	We used the figures from column ``Old'' for that.
	%~
	Note that column ``Pruned'' would result in a slightly
	different order.
	%~
	Also note that
	due to the nature of the pruning described in
	Sect.~\REF{Improved relation enumeration}, dual
	properties could have different extension set
	cardinalities.
}
and ensured
that the algorithm checks rectangles in order of increasing
encoding.
%~
This way, \QRef{226} corresponds to the encoding
\texttt{0x08000260},
while its weaker consequence would correspond to a larger encoding
\texttt{0x08002060} and therefore isn't found by the algorithm.

However, when a predicate occurs negated in a law,
this order doesn't lead to the desired result.
%~
For example, the law
``LfUnique $\land$ Irrefl $\Ra$ AntiTrans'' (\QRef{127})
is found,
but its weaker consequence
``LfUnique $\land$ ASym $\Ra$ AntiTrans'' (\QRef{109})
was found before it;
the former and the latter
corresponds to the encoding \texttt{0x002120} and
\texttt{0x000320}, respectively.

The latter problem is caused by our too simple enumeration
method.
%~
In procedure \texttt{qmc}, we iterate in the loop on
\texttt{mask} over all property sets of cardinality
given by the current \texttt{level}.
%~
For each such set we the iterate in the loop on \texttt{val}
over all assignments of negation symbols to the properties.
%~
For example, referring to level 2 of
Fig.~\REF{Search graph for the Quine-McCluskey algorithm on 3 variables},
\texttt{mask} may take the values $ab$, $ac$, $bc$,
and for each value $xy$, \texttt{val} may take the values
$xy$, $x\overline{y}$, $\overline{x}y$, $\overline{x}\overline{y}$.
%~
However, assuming that the order by extension set cardinalities
is
$\overline{c} < \overline{b} < \overline{a} < a < b < c$,
we should check the set $ac$ before the larger set $bc$
but after the smaller set $\overline{b}c$,
which is impossible with our simple enumeration method.

\begin{figure}
\begin{center}
\includegraphics[scale=0.5]{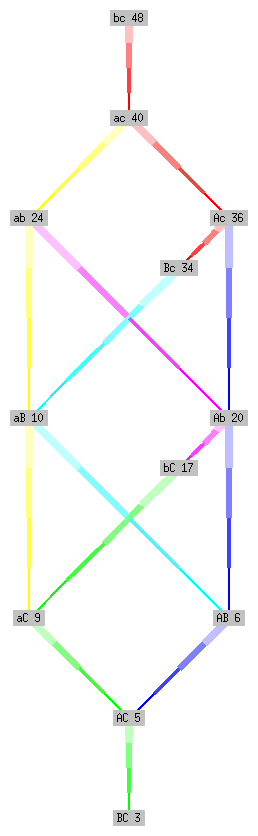}
%\hfill
%$\begin{array}[b]{|r | c@{}c@{}c | c | c@{}c@{}c |}
%\hline
%\rm wg & \multicolumn{3}{c|}{\rm rect}
%	& \texttt{mask} & \multicolumn{3}{c|}{\texttt{val}} \\\hline
%\hline
%48&            &          b &          c &\tt 011&\tt -&\tt 0&\tt 0\\\hline
%40&          a &            &          c &\tt 101&\tt 0&\tt -&\tt 0\\\hline
%36&\overline{a}&            &          c &\tt 101&\tt 1&\tt -&\tt 0\\\hline
%34&            &\overline{b}&          c &\tt 011&\tt -&\tt 1&\tt 0\\\hline
%24&          a &          b &            &\tt 110&\tt 0&\tt 0&\tt -\\\hline
%20&\overline{a}&          b &            &\tt 110&\tt 1&\tt 0&\tt -\\\hline
%17&            &          b &\overline{c}&\tt 011&\tt -&\tt 0&\tt 1\\\hline
%10&          a &\overline{b}&            &\tt 110&\tt 0&\tt 1&\tt -\\\hline
% 9&          a &            &\overline{c}&\tt 101&\tt 0&\tt -&\tt 1\\\hline
% 6&\overline{a}&\overline{b}&            &\tt 110&\tt 1&\tt 1&\tt -\\\hline
% 5&\overline{a}&            &\overline{c}&\tt 101&\tt 1&\tt -&\tt 1\\\hline
% 3&            &\overline{b}&\overline{c}&\tt 011&\tt -&\tt 1&\tt 1\\\hline
%\multicolumn{8}{c}{}	\\
%\multicolumn{8}{c}{}	\\
%\multicolumn{8}{c}{}	\\
%\multicolumn{8}{c}{}	\\
%\multicolumn{8}{c}{}	\\
%\multicolumn{8}{c}{}	\\
%\multicolumn{8}{c}{}	\\
%\multicolumn{8}{c}{}	\\
%\end{array}$
% table contents as ASCII:
%	48	.bc	011	*00
%	40	a.c	101	0*0
%	36	A.c	101	1*0
%	34	.Bc	011	*10
%	24	ab.	110	00*
%	20	Ab.	110	10*
%	17	.bC	011	*01
%	10	aB.	110	01*
%	09	a.C	101	0*1
%	06	AB.	110	11*
%	05	A.C	101	1*1
%	03	.BC	011	*11
\caption{Partial rectangle order Induced by
	%$c < \overline{a} < \overline{b} < b < a < \overline{c}$
	$\overline{c} < \overline{b} < \overline{a} < a < b < c$}
\LABEL{Partial rectangle order Induced by}
\end{center}
\end{figure}

Given the cardinality of each property's extension set
(as in
Fig.~\REF{Number of relations on a $5$-element set}), a partial
order on the intersection set is induced;
Fig.~\REF{Partial rectangle order Induced by}\footnote{
	Capital letters denote negations.
}
gives an example for the facts about two-set
intersections inferrable\footnote{
	E.g.\ $a < b$ implies $ac < bc$ by monotonicity.
}
from the above order
$\overline{c} < \overline{b} < \overline{a} < a < b < c$.

An improved approach should enumerate the rectangles in an
order that is some linearization of this inferred partial order.
%~
An efficient method to do this is still to be found.
%~
One possibility might be to assign to each property a weight,
such that increasing extension set cardinalities correspond to
increasing weights,
and to linearize the intersections in order of
increasing weight sums.
%~
Choosing powers of $2$ as weights will guarantee that all
weight sums are distinct; cf.\ the numbers in
Fig.~\REF{Partial rectangle order Induced by}.
%~
However, sorting a list of all
$\left( \begin{array}{c} 24 \\ n \end{array} \right) \cdot 2^n$
rectangles on level $n$ would definitely not be efficient.

%Red arcs indicate comparability along $c$,
%magenta: $b$ ($bc \ra ab$ induced),
%yellow: $a$,
%blue: $\overline{a}$,
%cyan: $\overline{b}$
%($\overline{a} \overline{b} \ra \overline{b} \overline{c}$ induced),
%green: $\overline{c}$.
%~
%Weighting numbers are obtained by summing literal weights,
%with
%$\overline{c}:1,
%\overline{b}:2,
%\overline{a}:4,
%a:8,
%b:16
%c:32$.

\paragraph*{Criterion \REF{nice 3}}

Considering sets of laws, some balance should be kept between
conciseness and convenience.

For example, in a textbook about
commutative groups, commutative variants
of the associativity axiom,
like $(x y) z = (y z) x$,
usually aren't explicitly mentioned, in
order to keep the presentation concise.
%~
On the other hand, while all theorems are redundant in the presence of
an axiomatization, the book will undoubtedly present some of them for
convenience.

In our setting, we have a simple formal criterion about which laws to
consider redundant:
%~
those that follow from other laws solely by propositional logic.
%~
For example, ``CoRefl$\Ra$LfEucl'' (\QRef{006}), is not considered
redundant despite the triviality of its proof, since the latter needs
to use Def.~\REFF{def}{CoRefl} and~\REFF{def}{LfEucl}.
%~
In contrast,
``LfUnique$\land$IncTrans$\Ra$ASym$\lor$LfSerial'' (\QRef{239})
is redundant since it follows from
``LfUnique$\land$SemiOrd2$\Ra$ASym$\lor$LfSerial'' (\QRef{242})
and
``SemiOrd2$\Ra$IncTrans'' (\QRef{071})
by propositional inference alone, without employing Def.~\REF{def}.

\begin{figure}
\begin{center}
\setlength{\unitlength}{1cm}
\renewcommand{\karnNumStyle}[1]{}
\renewcommand{\karnMrkStyle}[2]{%
	\put( 0.00,-0.40){\makebox(0,0)[b]{\scriptsize\sf #1}}%
}
\karnIII{%
\karnCmd{\a}{\i}{\r}%
\karnNam{as}{ir}{rf}%
\karnBox%
\karnCap%
\karnMul\i{\karnMul\r{ \karnMrk{\raisebox{4ex}{\QRef{039}}} }} {      }%
\karnMul\r{            \karnMrk{\raisebox{4ex}{\QRef{039}}}  } {\a\i  }%
\karnMul\r{            \karnMrk{\raisebox{4ex}{--        }}  } {\a    }%
\karnMul\i{\karnMul\r{ \karnMrk{\raisebox{2ex}{\QRef{044}}} }} {      }%
\karnMul\i{            \karnMrk{\raisebox{2ex}{\QRef{044}}}  } {\a    }%
\karnMul\i{            \karnMrk{\raisebox{2ex}{--        }}  } {\a  \r}%
\karnMul\a{\karnMul\r{ \karnMrk{\raisebox{0ex}{\QRef{046}}} }} {      }%
\karnMul\a{            \karnMrk{\raisebox{0ex}{\QRef{046}}}  } {  \i  }%
\karnMul\a{            \karnMrk{\raisebox{0ex}{--        }}  } {  \i\r}%
}
\caption{Redundant law suggestion example}
\LABEL{Redundant law suggestion example}
\end{center}
\end{figure}

As can be seen from the previous example, our algorithm doesn't avoid
reporting redundant law suggestions.
%~
The reason for this is that it just reports prime implicants in order
of appearance, rather than selecting a minimal covering subset of
them.
%~
The latter technique is commonly employed in proper Quine-McCluskey
implementations, it is, however, NP-complete
\CITEPAGE{p.14}{Feldman.2009}.
%~
The Karnaugh diagram in Fig.~\REF{Redundant law suggestion example}
illustrates the problem in a simplified setting ($3$ properties only)
along the example
``ASym$\Ra$Irrefl'' (\QRef{039}),
``$\lnot$Refl$\lor \lnot$ASym'' (\QRef{044}), and
``$\lnot$Refl$\lor \lnot$Irrefl'' (\QRef{046}).
%~
Rectangle \QRef{044} is inspected before \QRef{046}, so the former is
reported and its fields are set to \texttt{don't care}.
%~
When the latter is inspected, there is still one \texttt{on} field in it,
so it is reported, too.
%~
In contrast, a minimal covering subset approach wouldn't report
\QRef{044} since all its fields are covered by the union of \QRef{039}
and \QRef{046}.

\section{References}

\bibliographystyle{elsarticle-num-names}
\bibliography{lit}

\end{document}